\documentclass[a4paper]{amsart}

\usepackage{relsize}
\usepackage{amsmath}
\usepackage[dvipsnames]{xcolor} 
\usepackage{physics}
\usepackage[shortlabels]{enumitem} 
\setlist[enumerate]{topsep=0pt} % reduce space before enu
\setlist[itemize]{topsep=0pt}  
\def\<#1>{\mathinner{\left<#1\right>}}

\usepackage{hyperref}
\hypersetup{colorlinks}
\hypersetup{
linkcolor=Periwinkle!70!BlueViolet
,citecolor=Green
,filecolor=BurntOrange
,urlcolor=Blue!50!PineGreen
,menucolor=BrickRed
,runcolor=YelloOrange
,linkbordercolor=BrickRed
,citebordercolor=Green
,filebordercolor=Mulberry
,urlbordercolor=NavyBlue
,menubordercolor=BrickRed
,runbordercolor=BurntOrange
}

\usepackage{dirtytalk}

\newcommand{\R}{\mathbb{R}}
\newcommand{\Z}{\mathbb{Z}}
\newcommand{\Q}{\mathbb{Q}}
\newcommand{\F}{\mathbb{F}}
\newcommand{\E}{\mathbb{E}}
\newcommand{\bbS}{\mathbb{S}}

\usepackage[font=footnotesize]{caption}

\usepackage{graphicx}

\makeatletter
\newcommand*{\defeq}{\mathrel{\rlap{%
                     \raisebox{0.3ex}{$\m@th\cdot$}}%
                     \raisebox{-0.4ex}{$\m@th\cdot$}}%
                     =}
\makeatother

\usepackage{setspace}

\usepackage{dashrule}
\usepackage{csquotes}
\usepackage{amsthm}
\usepackage[capitalise,nameinlink]{cleveref}
\usepackage{cancel}
\usepackage{mathtools}

\usepackage{listings}

\usepackage{mathrsfs}

\DeclareMathOperator{\St}{St}    %for Steinberg modules

\DeclareMathOperator{\GL}{GL}

\DeclareMathOperator{\gr}{gr}

\DeclareMathOperator{\im}{im}

\DeclareMathOperator{\Cobar}{Cobar}

\DeclareMathOperator{\Spec}{Spec}
\DeclareMathOperator{\Sq}{Sq}

\DeclareMathOperator{\SL}{SL}

\DeclareMathOperator{\Map}{Map}

\DeclareMathOperator{\Aut}{Aut}

\DeclareMathAlphabet{\mathantt}{OT1}{antt}{li}{it}
\DeclareMathAlphabet{\mathpzc}{OT1}{pzc}{m}{it}

\usepackage{amssymb}
\usepackage{quiver}
\usepackage{fontawesome}

\DeclareMathOperator{\sk}{sk}

\newcommand{\qedd}{\hfill $\blacksquare$ \vspace{1.3em}} 
\newcommand{\colim}{\mathrm{colim}}
\newcommand{\Fun}{\cat{Fun}}

\newcommand{\Alg}{ \cat{Alg}}
\newcommand{\Coalg}{ \cat{Coalg}}
\newcommand{\bbk}{ \mathbb{k}}

\newcommand{\Spaces}{\mathscr{S}}

\newcommand{\cat}[1]{\mathrm{#1}} % special font for category
\newcommand{\fancycat}[1]{\mathscr{#1}} % special font for category
\newtheoremstyle{exampstyle}
{7pt} % Space above
{5pt} % Space below
{} % Body font
{} % Indent amount
{\bfseries } % Theorem head font
{.} % Punctuation after theorem head
{.5em} % Space after theorem head
{\thmname{#1}\thmnumber{ #2}\thmnote{ (#3)}} % Theorem head spec (can be left empty, meaning `normal')

\newtheoremstyle{thmstyle}
{7pt} % Space above
{5pt} % Space below
{\itshape} % Body font
{} % Indent amount
{\bfseries } % Theorem head font
{.} % Punctuation after theorem head
{.5em} % Space after theorem head
{\thmname{#1}\thmnumber{ #2}\thmnote{ (#3)}} % Theorem head spec (can be left empty, meaning `normal')

\newtheoremstyle{inlinestyle}
{5pt} % Space above
{5pt} % Space below
{\itshape} % Body font
{1em} % Indent amount
{ } % Theorem head font
{.} % Punctuation after theorem head
{.5em} % Space after theorem head
{<\thmname{#1}\thmnumber{ #2}\thmnote{ #3}>}

\theoremstyle{thmstyle}
\newtheorem{thm}{Theorem}[section]
\newtheorem{prop}[thm]{Proposition}
\newtheorem{lem}[thm]{Lemma}
\newtheorem{cor}[thm]{Corollary}

\newtheorem*{claim}{Claim}
\newtheorem{claimnumbered}[thm]{Claim}

\theoremstyle{exampstyle}

\newtheorem{obs}[thm]{Observation}
\newtheorem{notation}[thm]{Notation}

\newtheorem{construction}[thm]{Construction}
\newtheorem{defn}[thm]{Definition}
\newtheorem{exmp}[thm]{Example}
\newtheorem{rem}[thm]{Remark}

\usepackage{thmtools}

\declaretheoremstyle[%
spaceabove=10pt,%
spacebelow=5pt,%
headfont=\normalfont\itshape,%
notebraces={}{},
headformat = \NAME \NOTE,
postheadspace=1em,%
qed=$\blacksquare$%
]{mystyle} 
\declaretheorem[name={Proof},style=mystyle,unnumbered,]{prf}

\declaretheoremstyle[%
spaceabove=5pt,%
spacebelow=3pt,%
headfont=\normalfont\itshape,%
postheadspace=1em,%
qed=$\blacktriangleleft$%
]{mystyle2} 
\declaretheorem[name={Proof of claim},style=mystyle2,unnumbered,]{prfoc}
\DeclareMathOperator{\id}{id}

\DeclareMathOperator{\Fil}{Fil}
\DeclareMathOperator{\fil}{fil}

\DeclareMathOperator{\Tor}{Tor}
\DeclareMathOperator{\Cotor}{Cotor}

\makeatletter
\newcommand{\leqnomode}{\tagsleft@true\let\veqno\@@leqno}
\makeatother
\usepackage[style=alphabetic,sorting=debug]{biblatex}

\renewbibmacro{in:}{}
\renewbibmacro{url:}{}
\renewbibmacro{doi:}{}

\DeclareFieldFormat{labelalpha}{\thefield{entrykey}}
\DeclareFieldFormat{extraalpha}{}

\definecolor{grass}{HTML}{829356}
\definecolor{cinnamon}{HTML}{c24e0a}
\definecolor{gold}{HTML}{BCA136}
\definecolor{flame}{rgb}{0.89, 0.35, 0.13}
\definecolor{Amethyst}{HTML}{8953C6}
\definecolor{goodpink}{HTML}{e02963}

\usepackage{colortbl}

\usepackage{yhmath} % \widetriangle and \wideparen ]

\usepackage[colorinlistoftodos]{todonotes}

\usepackage{tikz}
\usepackage{tensor} % to decorate the bockstein spectral sequence
\usepackage{spectralsequences}
\usepackage{adjustbox}
\DeclareMathOperator{\UT}{UT}
\DeclareMathOperator{\ab}{ab}

\DeclareMathOperator{\CGL}{CGL}
\DeclareMathOperator{\cgl}{cgl}

\DeclareMathOperator{\Ba}{Bar}

\usepackage{scalerel}

\usepackage[ margin=8em]{geometry}

\usepackage{graphicx}
\usepackage{tikz-cd}

\usepackage[capitalise,nameinlink]{cleveref}
\usepackage{wasysym}

\usepackage[full]{textcomp} % to get the right copyright, etc.
\usepackage[p,osf]{fbb} % osf in text, tabular lining figures in math
\usepackage[scaled=.95,type1]{cabin} % sans serif in style of Gill Sans
\usepackage[varqu,varl]{zi4}% inconsolata typewriter
\usepackage[T1]{fontenc} % LY1 also works

\usepackage[libertine]{newtxmath}
\usepackage[cal=boondoxo,bb=boondox,frak=boondox]{mathalfa} 
\usepackage[mathscr]{euscript}
\title{Periodic families in the homology of $\GL_n(\F_2)$}

\author{Kelly Wang}
\begin{document}
\begin{abstract}
    We construct infinite families of nonzero classes in  $H_d(\GL_n(\F_2);\F_2)$  along lines of the form $d =\frac{2}{3}n + $(constant), thereby showing that the known slope $\frac{2}{3}$-stability for these homology groups are optimal. Using the new stability Hopf algebra  perspective of Randal-Williams,  our computations in addition recover the slope-$\frac{2}{3}$ stability for $\GL_n(\Z)$ with coefficients in $\F_2$, improve that for $\Aut(F_n)$ to $\frac{2}{3}$, and demonstrate that those slopes are optimal.
    Perhaps of independent interest, we also provide a manual for computing stability Hopf algebras over $\F_2$.
\end{abstract}

%------------------- Spacings ------------------------------
\setlength{\parindent}{0pt}
\setlength{\parskip}{1ex}
%-----------------------------------------------------------
%------------------- Section settings ----------------------
\setcounter{tocdepth}{2}
\setcounter{secnumdepth}{3}
%
%-----------------------------------------------------------

\maketitle
\renewcommand{\baselinestretch}{0.75}\normalsize
\tableofcontents
\renewcommand{\baselinestretch}{1.0}\normalsize
\makeatletter
\renewcommand\subsubsection{\@startsection{subsubsection}{3}{\z@}%
{-3.25ex\@plus -1ex \@minus .2ex}%
{-1em}%
{\normalfont\normalsize\bfseries \relscale{0.95}}}
\makeatother

\section{Introduction}
When $l$ is an odd prime, the cohomologies $H^*(\GL_n(\F_2);\F_l)$ are completely determined by Quillen in his seminal paper \autocite[]{Qui72}.  With coefficients in $\F_2$, very little is known. Quillen showed in the same paper that the  stable homology $\colim_{n\to \infty}\widetilde{H}_*(\GL_n(\F_2);\F_2)$  vanishes, and this combined with the stability result of \autocite[Theorem B]{GKRW18b} tells us the $\F_2$-homology groups exhibit a slope-$\frac{2}{3}$ vanishing line: $\widetilde{H}_d(\GL_n(\F_2);\F_2)=0$ for $d<\frac{2}{3}n-1$. One natural yet previously unresolved question is whether the $\frac{2}{3}$ slope can be improved. In this paper show that the answer is negative by exhibiting infinite families of non-zero classes along lines of slope $\frac{2}{3}$.
\begin{thm}[\cref{thm: absolute clases}]\label{intro thm: abs}
    For every $i\geqslant 0$, there are non-zero classes
    \begin{align*}
u_{i0}& \in H_{8i+1}(\GL_{12i+2}(\F_2);\F_2)    \\
u_{i1}&\in H_{8i+2}(\GL_{12i+4}(\F_2);\F_2)\\
u_{i2}&\in H_{8i+3}(\GL_{12i+6}(\F_2);\F_2)\\
s_{i+1}&\in  H_{8i+7}(\GL_{12i+11}(\F_2);\F_2).
    \end{align*}
    $u_{00}\in H_1(\GL_2(\F_2);\F_2)$ is represented by the matrix $\smqty(  0& 1 \\ 1&0 )$, and  $u_{i1} = u_{00}u_{i0}$, $u_{i2} = u_{00}^2u_{i0}$ under the algebra structure of $H_*(\GL_*(\F_2);\F_2)$.
 \begin{figure}[h]
     \centering
     \adjustbox{scale=0.83}{
        \NewSseqGroup\mygroupnew {} {
\class(2,1)
\class(4,2)
\class(6,3)
\class(11,7)
\structline[orange](2,1)(4,2)
\structline[orange](4,2)(6,3)
}
\begin{sseqpage}[ classes = fill,  xscale = 0.6,  yscale = 0.5, grid = go, grid color =gray!20, xrange = {0}{12}, yrange = {0}{7}, tick step = 3,
major tick step = 3,
minor tick step = 1, x tick handler = {
\ifnum#1 = 0\relax
12i
\else
\ifnum#1 = 12\relax
% \vphantom is fragile so we \protect it
\protect\vphantom{2}12i+#1
\else
\ifnum#1 = 6\relax
% \vphantom is fragile so we \protect it
\protect\vphantom{2}12i+#1
\else
\relax 
\fi
\fi
\fi
},  y tick handler = {
\ifnum#1 = 0\relax
8i
\else
\relax 8i+#1\fi
}]
\mygroupnew(0,0)
\classoptions["{u}_{i0}" right](2,1)
\classoptions["{u}_{i1}" right](4,2)
\classoptions["{u}_{i2}" right](6,3)
\classoptions["{s}_{i+1}" right](11,7)
\end{sseqpage}}
     \caption{The non-zero classes of \cref{intro thm: abs}; $i\geqslant 0$. The orange line indicates an $u_{00}$-extension.}
     \label{fig:absolute classes}
 \end{figure}
\end{thm}

In particular, if $\lambda>0$ is such that there exists some constant $\mu$ with  $\widetilde{H}_d(\GL_n(\F_2);\F_2)=0$ whenever $d<\lambda n+\mu$, then one must have $\lambda\leqslant\frac{2}{3}$. It is worth noting that Lahtinen-Sprehn have constructed families of  non-zero cohomology classes by very different means \autocite[Theorem 1]{LS18}. For instance, they constructed 
for every $n\geqslant 1$, a non-zero cohomology classes in $H^{2^{n}-1}(\GL_{2^n}(\F_2);\F_2)$. However the classes they constructed all lie above the slope-$1$ line $d= n-2$.

Our investigation lies within the general framework of homological stability. Given a sequence of groups \[G_0\xrightarrow[]{i_0}G_1\xrightarrow[]{i_1}\cdots\] and a fixed abelian group of coefficient $M$, one can ask the following question. {\itshape Does there exists a function $\alpha(n)$ with $\alpha(n)\to \infty$ as $n\to \infty$, such that $H_d(i_n;M)$ is an isomorphism whenever $d<\alpha(n)$? }
If such a function of the form $\alpha(n) = \lambda n +\mu$ exists, we say that the $G_n$ have a \textit{slope $\lambda$ stability} with coefficients in $M$. We would like to remark that exhibiting stability and proving non-stability are two tasks very distinct in flavour. Moreover, a wide range of techniques have been developed for tackling the former, for instance Quillen's spectral sequence method \autocite[1974-1; pp.1-15]{QuiN}, \autocite[]{SW20},  and Galatius-Kupers-Randal-Williams's $\E_k$-cell machine \autocite[]{GKRW18a}. The latter task, however, has not received much attention. We hope this document contributes to this direction.  The calculations going into slope-$\frac{2}{3}$ stability result in \autocite[]{GKRW18b} suggests\footnote{More precisely, they essentially had $\varphi_*$ on a  filtration quotient; c.f. the class $\tau^4$ in \autocite[\S6.4]{GKRW18b}.} a  $(12,8)$ self-map
\stepcounter{thm}
\[\label{eq: self map on rel homology}\leqnomode\tag{\thethm}\varphi_{*}:H_{d-8}(\GL_{g-12}(\F_2),\GL_{g-12-1}(\F_2);\F_2)\to H_d(\GL_g(\F_2),\GL_{g-1}(\F_2);\F_2).\]
So to preclude any improved slope it is enough to show that $\varphi_*$ is non-nilpotent. 

\begin{thm} [\cref{prop: BGL/sigma classes}]\label{intro thm: rel}
The classes \[\alpha_{0j}\in H_j(\GL_{2j}(\F_2),\GL_{2j-1}(\F_2);\F_2),\quad \gamma_{0j}\in H_{j+2}(\GL_{2j+3}(\F_2),\GL_{2j+2}(\F_2);\F_2),\quad j= 0,1,2\]  of \cref{prop: BGL/sigma classes} are such that 
 $\alpha_{ij}\defeq \varphi_{*}^i(\alpha_{0j})$ and $\gamma_{ij}\defeq \varphi^i_{*}(\gamma_{0j})$ are non-zero for all $i\geqslant 0$.
\end{thm}

Let's  briefly comment on the situation for general linear groups of finite field other than $\F_2$. Let $p$ be a prime and $q>2$ be a $p$-power. Away from the characteristic,    
the cohomologies $H^*(\GL_n(\F_q);\F_l)$, $l$ a  prime distinct from $p$, have been completely computed by Quillen \autocite[Theorem 4]{Qui72}. With coefficients in $\F_p$ The stable homology  is trivial \autocite[Corollary 2]{Qui72} and the homology groups enjoy a slope $1$ stability \autocite[Theorem A]{GKRW18b}. The classes constructed in \autocite[]{LS18} shows that they cannot have beyond slope $2$ stability. One reason that we do not investigate the sharpness of the slope in these cases is that the methods in \autocite{RW25} only deals with stability with slope less than $1$.

Using the technologies developed in \autocite[]{RW25}, we are able to propagate the ingredients of the above theorem to consequences for general linear groups over the integers and automorphism groups of free groups, via the homomorphisms\stepcounter{thm}
\[\tag{\thethm}\Aut(F_n)\to \GL_n(\Z)\to \GL_n(\F_2).\label{eq: aut to gl}\leqnomode\]
\begin{thm}[\cref{sec: stability for Z}]\label{intor thm: cor}
    \begin{enumerate}[(i),before*=\leavevmode\vspace{-0.17\baselineskip}]
        \item   $H_d(\GL_n(\Z),\GL_{n-1}(\Z);\F_2)=0$ for $d<\frac{2}{3}(n-1)$, and this slope is optimal.
        \item $H_d(\Aut(F_n),\Aut(F_{n-1});\F_2)=0$ for $d<\frac{2}{3}(n-1)$, and this slope is optimal.
    \end{enumerate}
    In particular, the Bockstein long exact sequence associated $\Z\to \Z\to \Z/2$ shows these relative groups with integer coefficients cannot have a vanishing line with slope $>\frac{2}{3}$. 
\end{thm}
To the best of our knowledge,  the sharpness of $\frac{2}{3}$ in both cases was not known previously. The vanishing line  for $\GL_n(\Z)$  have been established in \autocite[Theorem C]{KMP22}\footnote{The statement there contains a misprint, the inequalities should be strict.}.  For $\Aut(F_n)$ the best known vanishing slope was $\frac{1}{2}$, coming from that for integer coefficients  \autocite[Theorem G]{WRW17}, \autocite[Theorem 18.3]{GKRW18a}.
In fact we will establish the \say{slope-$\frac{2}{3}$ bound}  for general linear groups over a slightly more general class of rings.
\begin{thm}[\cref{prop: GLnZ and AutFn}]
    Let $\fancycat{O}$ be Dedekind domain that has $\F_2$ a a quotient. Then $H_d(\GL_n(\fancycat{O});\F_2)$ cannot have stability of slope $>\frac{2}{3}$.
\end{thm}

The currently known best stability slope for general linear goups a general Dedekind domain is $\frac{1}{2}$ due to  van der Kallen \autocite[Theorem 4.11]{vdK80} and later Randal-Williams \autocite[Theorem 1.1]{RW24}.

The moral behind deducing \cref{intor thm: cor} from the calculations for \cref{intro thm: rel} is that the map on \textit{stability Hopf algebras} induced by the homomorphisms (\ref{eq: aut to gl}) is an isomorphism in gradings at most $5$. Suppose $\bbk$ is a field and $R$ is an $\E_\infty$-$\bbk$-module satisfying some mild condition, then a central result of  Randal-Williams \autocite[]{RW25} is that one can extract a Hopf algebra over $\bbk$ which \say{carries the same amount of stability information as $R$}. We interpret (higher) homological stability of $R$ as vanishings of cofibre of (iterated) self-maps on $R$. A more detailed exposition can be found in \S\ref{sec: tools}.

\subsection{Outline}

In \S\ref{sec: group homology F2} we describe an additive basis for $ H_d(\GL_n(\F_2);\F_2)$ in the range $n\leqslant 6, d\leqslant 4$, in terms of Dyer-Lashof operations on the class $\sigma\in H_0(\GL_1(\F_2);\F_2)$ and any choice of basis $\nu_1,\nu_2$ for $ H_3(\GL_3(\F_2);\F_2)$. The operations on the homology groups comes from the realisation of those groups as homotopy groups of an $\E_\infty$-algebra $\CGL(\F_2)$ in $\F_2$-modules. However, the arguments here only use formal properties of those operations, and hence the section is written in a way  such that minimal $\E_\infty$-knowledge is needed. 
In \S\ref{sec: group homology Z} we provide the analogous description of  $H_d(\GL_n(\Z);\F_2)$ in the narrower range of $n\leqslant 3$, $d\leqslant 2$. 

\S\ref{sec: notations} consists of three parts. \S\ref{sec: background machineries} is a recollection of $\infty$-categories, filtered objects, and the diagonal $t$-structure on $\Z$-graded modules. In \S\ref{sec: E_infty algs} we describe some assumptions that we often put on our $\E_\infty$-algebras. In \S\ref{subsubsec: bar and cobar} we review the Bar-Cobar construction and define a \textit{bar operation} on homotopy groups of $\E_n$-algebras. We then study how this operation interact with Dyer-Lashof operations and the Browder bracket. This part together with \S\ref{sec: tools} provides a set of tools for computing the \textit{stability Hopf algebra} $\Delta_A$ of an $\E_n$-algebra $A$ in $\F_2$-modules. This paper is essentially an application of the new stability Hopf algebra viewpoint of homological stability developed in \autocite[]{RW25}. In \S\ref{sec: tools} review the relevant concept and describe (in most cases) how to obtain $\Delta_A$ from an $\E_n$-cell structure on $A$.

In \S\ref{sec: building the detector} we build all the pieces needed to prove the theorems in the introduction. Those pieces are then assembled in \S\ref{sec: applications of the detector}. The strategy is simple:  \S\ref{sec: group homology F2} provides us with an initial segment of a minimal $\E_\infty$-cell structure of $\CGL(\F_2)$, using which we can compute the stability Hopf algebra $\Delta_{\CGL(\F_2)}$ up to grading $5$. We then check by hand that there is a Hopf algebra map $\Delta_{\CGL(\F_2)}$ to $A(1)_*$, the dual of the sub Hopf algebra of the Steenrod algebra generated by $\Sq^1$ and $\Sq^2$. Taking $\Cobar$ this becomes a map of $\E_\infty$-algebras $\cgl(\F_2)\to a(1)$. The algebra $\cgl(\F_2)$ is a truncated version of $\CGL(\F_2)$ that contains the same amount of stability information and it comes with  an $\E_\infty$-map $\CGL(\F_2)\to \cgl(\F_2)$. The key player is the composite map   $\CGL(\F_2)\to a(1)$.  By computing $\pi_{*,*}(a(1))$ via a May-type spectral sequence, one sees that $\pi_{*,*}(a(1)/\sigma)$ is $(12,8)$-periodic. Our next task is to lift this periodicity to $\CGL(\F_2)/\sigma$. In fact, a big portion of the work lies in building an $(12,8)$-endomorphism $\varphi$ on $\CGL(\F_2)/\sigma$ that realises the periodicity on $\pi_{*,*}(a(1)/\sigma)$. As insinuated by the notation, the map $\varphi_*$  (\ref{eq: self map on rel homology}) is induced by $\varphi$. 

Finally, 
\S\ref{sec: stability for Z} is devoted to stability consequences of $\CGL(\fancycat{O})$ and $\Aut(F_n)$, as explained at the end of introduction.

\subsection{Acknowledgements} I am deeply grateful to my PhD supervisor Oscar Randal-Williams  for his generous guidance, constant encouragement, and countless invaluable discussions. I would like to sincerely thank Jan Steinebrunner for many insightful conversations. Many thanks to Jan, Oscar, and Robin Stoll, whose comments on an earlier draft have substantially improved this document.

\section{Low degree group homology computations: $\GL_g(\F_2)$}\label{sec: group homology F2}

The homology groups $H_{g,d} \defeq H_d(\GL_g(\F_2);\F_2)$ are that of an $\mathbb{N}$-graded $\E_\infty$-space, so are equipped with a multiplication and Dyer-Lashof operations $Q^1, Q^2,Q^3,\dots$ making the bigraded vector space $H_{*,*}$ into a \textit{$W_\infty$-algebra} \autocite[Section 16]{GKRW18a}. This means that the multiplication and operations satisfies the following properties
\begin{enumerate}[(a)]
    \item The $Q^s$ are linear.
    \item If a class $x$ has homological degree $<s$, then $Q^s(x)=0$.
    \item The $Q^s$ satisfy the Adem relations; see \autocite[Section 16.2.2. (c')]{GKRW18a}.
    \item If a class $x$ has homological degree $s$, then $Q^s(x) =x^2$.
    \item $Q^s(1)=0$, where $1\in H_{*,*}$ is the unit.
    \item The $Q^s$ satisfy the Cartan formula; see  \autocite[Section 16.2.2. (f')]{GKRW18a}.
\end{enumerate}
In this section we analyse the $W_\infty$-structure on  $H_{*,*}$ in low degrees. We formulate the answer in terms of operations on the classes $\sigma,\nu_1,\nu_2$, where $\sigma\in H_{1,0}$ is the generator, and $\nu_1,\nu_2$ is any fixed basis of $H_{3,3}$ (the fact that $H_{3,3}$ is $2$-dimensional will be explained at the start of \S\ref{sec: grouphom 1}).
\begin{thm}\label{thm: additive basis}
    \cref{fig: additive basis} describes a basis for $H_{g,d}$ in  the range  
 $g\leqslant 6$, $d\leqslant 4$.
    \begin{figure}[h]
         \centering
\begin{tikzcd}[sep=small]
	& {} & {} \\
	4 &&& {Q^4(\sigma)} & {\sigma Q^4(\sigma)} & \begin{array}{c} Q^2(\sigma)^2, \\\sigma^2Q^4(\sigma) \end{array} & {\sigma Q^2(\sigma)^2} & {Q^1(\sigma)^2Q^2(\sigma)} \\
	3 &&& {Q^3(\sigma)} & \begin{array}{c} \nu_1,\\\nu_2 \end{array} & \begin{array}{c} Q^{2,1}(\sigma),\\Q^1(\sigma)Q^2(\sigma) \end{array} & {\sigma Q^{2,1}(\sigma)} & {Q^1(\sigma)^3} \\
	2 &&& {Q^2(\sigma)} & {\sigma Q^2(\sigma)} & {Q^1(\sigma)^2} &&& {} \\
	1 &&& {Q^1(\sigma)} &&&&& {} \\
	0 & 1 & \sigma & {\sigma^2} & {\sigma^3} & {\sigma^4} & {\sigma^5} & {\sigma^6} & {} \\
	{d/g} & 0 & 1 & 2 & 3 & 4 & 5 & 6
	\arrow[draw={rgb,255:red,87;green,87;blue,87}, dotted, no head, from=2-1, to=2-4]
	\arrow[draw={rgb,255:red,87;green,87;blue,87}, dotted, no head, from=2-4, to=2-5]
	\arrow[draw={rgb,255:red,87;green,87;blue,87}, dotted, no head, from=2-4, to=3-4]
	\arrow[draw={rgb,255:red,87;green,87;blue,87}, dotted, no head, from=2-5, to=2-6]
	\arrow[draw={rgb,255:red,87;green,87;blue,87}, dotted, no head, from=2-5, to=3-5]
	\arrow[draw={rgb,255:red,87;green,87;blue,87}, dotted, no head, from=2-6, to=2-7]
	\arrow[draw={rgb,255:red,87;green,87;blue,87}, dotted, no head, from=2-6, to=3-6]
	\arrow[draw={rgb,255:red,87;green,87;blue,87}, dotted, no head, from=2-7, to=2-8]
	\arrow[draw={rgb,255:red,87;green,87;blue,87}, dotted, no head, from=2-7, to=3-7]
	\arrow[draw={rgb,255:red,87;green,87;blue,87}, dotted, no head, from=2-8, to=3-8]
	\arrow[draw={rgb,255:red,87;green,87;blue,87}, dotted, no head, from=3-1, to=3-4]
	\arrow[draw={rgb,255:red,87;green,87;blue,87}, dotted, no head, from=3-4, to=3-5]
	\arrow[draw={rgb,255:red,87;green,87;blue,87}, dotted, no head, from=3-4, to=4-4]
	\arrow[draw={rgb,255:red,87;green,87;blue,87}, dotted, no head, from=3-5, to=3-6]
	\arrow[draw={rgb,255:red,87;green,87;blue,87}, dotted, no head, from=3-5, to=4-5]
	\arrow[draw={rgb,255:red,87;green,87;blue,87}, dotted, no head, from=3-6, to=3-7]
	\arrow[draw={rgb,255:red,87;green,87;blue,87}, dotted, no head, from=3-6, to=4-6]
	\arrow[draw={rgb,255:red,87;green,87;blue,87}, dotted, no head, from=3-7, to=3-8]
	\arrow[draw={rgb,255:red,87;green,87;blue,87}, dotted, no head, from=3-7, to=6-7]
	\arrow[draw={rgb,255:red,87;green,87;blue,87}, dotted, no head, from=3-8, to=6-8]
	\arrow[draw={rgb,255:red,87;green,87;blue,87}, dotted, no head, from=4-1, to=4-4]
	\arrow[draw={rgb,255:red,87;green,87;blue,87}, dotted, no head, from=4-4, to=4-5]
	\arrow[draw={rgb,255:red,87;green,87;blue,87}, dotted, no head, from=4-4, to=5-4]
	\arrow[draw={rgb,255:red,87;green,87;blue,87}, dotted, no head, from=4-5, to=4-6]
	\arrow[draw={rgb,255:red,87;green,87;blue,87}, dotted, no head, from=4-5, to=6-5]
	\arrow[draw={rgb,255:red,87;green,87;blue,87}, between={0}{0.8}, dotted, no head, from=4-6, to=4-9]
	\arrow[draw={rgb,255:red,87;green,87;blue,87}, dotted, no head, from=4-6, to=6-6]
	\arrow[draw={rgb,255:red,87;green,87;blue,87}, dotted, no head, from=5-1, to=5-4]
	\arrow[draw={rgb,255:red,87;green,87;blue,87}, between={0}{0.9}, dotted, no head, from=5-4, to=5-9]
	\arrow[draw={rgb,255:red,87;green,87;blue,87}, dotted, no head, from=5-4, to=6-4]
	\arrow[draw={rgb,255:red,87;green,87;blue,87}, dotted, no head, from=6-1, to=6-2]
	\arrow[color={rgb,255:red,87;green,87;blue,87}, between={0}{0.9}, from=6-2, to=1-2]
	\arrow[color={rgb,255:red,87;green,87;blue,87}, no head, from=6-2, to=6-3]
	\arrow[draw={rgb,255:red,87;green,87;blue,87}, between={0}{0.9}, dotted, no head, from=6-3, to=1-3]
	\arrow[draw={rgb,255:red,87;green,87;blue,87}, no head, from=6-3, to=6-4]
	\arrow[draw={rgb,255:red,87;green,87;blue,87}, dotted, no head, from=6-3, to=7-3]
	\arrow[draw={rgb,255:red,87;green,87;blue,87}, no head, from=6-4, to=6-5]
	\arrow[draw={rgb,255:red,87;green,87;blue,87}, dotted, no head, from=6-4, to=7-4]
	\arrow[draw={rgb,255:red,87;green,87;blue,87}, no head, from=6-5, to=6-6]
	\arrow[draw={rgb,255:red,87;green,87;blue,87}, dotted, no head, from=6-5, to=7-5]
	\arrow[draw={rgb,255:red,87;green,87;blue,87}, no head, from=6-6, to=6-7]
	\arrow[draw={rgb,255:red,87;green,87;blue,87}, dotted, no head, from=6-6, to=7-6]
	\arrow[draw={rgb,255:red,87;green,87;blue,87}, no head, from=6-7, to=6-8]
	\arrow[draw={rgb,255:red,87;green,87;blue,87}, dotted, no head, from=6-7, to=7-7]
	\arrow[color={rgb,255:red,87;green,87;blue,87}, between={0}{0.4}, from=6-8, to=6-9]
	\arrow[draw={rgb,255:red,87;green,87;blue,87}, dotted, no head, from=7-8, to=6-8]
\end{tikzcd}

         \caption{Vector space basis for $H_d(\GL_g)$.}
         \label{fig: additive basis}
     \end{figure}
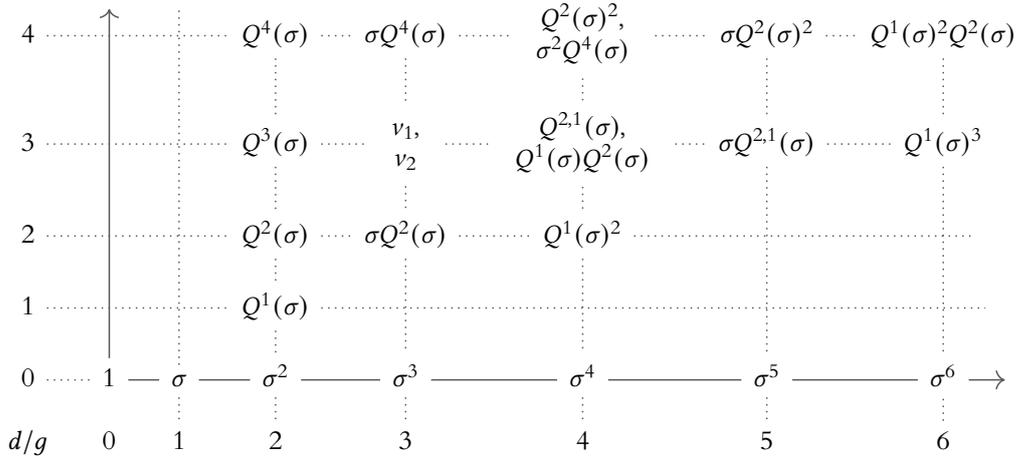
\end{thm}

Since determining \textit{all} the $W_\infty$-relations between those basis in this range seems to be difficult and is not essential for our purpose,  we restrict to the smaller range $g\leqslant6,d\leqslant 3$, where the relations can be summarised as follows.
\begin{thm}\label{thm: Winfty structure}
    \begin{enumerate}[(i),before*=\leavevmode\vspace{-0.17\baselineskip}]
            \item Let $W_\infty[\sigma,\nu_1,\nu_2]$ be the free $W_\infty$-algebra generated by symbols $\sigma,\nu_1,\nu_2$ of bidegree $(1,0),(3,3),(3,3)$ respectively. Let $I$ be $W_\infty$-ideal\footnote{Given a subset $S$ of a $W_\infty$-algebra $A$, we define the \textit{$W_\infty$-ideal} $(S)_{W_\infty}$ generated by $S$ to be the smallest subset of $A$ containing $S$ that is closed under multiplication and Dyer-Lashof operations. In particular, $(S)_{W_{\infty}}$ is an algebra ideal of $A$. The algebra quotient $A/(S)_{W_\infty}$ inherits a $W_\infty$-structure from $A$.} generated by the following list of elements
        \stepcounter{thm}
            \[\label{eq: classes that are zero}\leqnomode\tag{\thethm}\sigma Q^1(\sigma),\quad \sigma Q^3(\sigma),\quad \sigma^2 Q^2(\sigma),\quad \sigma\nu_1,\quad \sigma\nu_2.\]
            Then sending $\sigma,\nu_1,\nu_2$ to the classes of the same name in $H_{*,*}$ defines a $W_\infty$-map
            \[\varphi: W_{\infty}[\sigma, \nu_1,\nu_2]/I\to H_{*,*}.\]
            \item The map $\varphi$ is a bijection in the range $\{(g,d)\mid g\leqslant6,d\leqslant 3\}$.
    \end{enumerate}
\end{thm}

Finally we record some identities in  homological degree $4$. These are not needed for the computations in \S\ref{sec: building the detector} and \S\ref{sec: applications of the detector}, but are perhaps of independent interest.
\begin{cor}\label{cor: identities}

\begin{enumerate}[(i),before*=\leavevmode\vspace{-0.17\baselineskip}]
    \item $Q^1(\sigma)^2 Q^2(\sigma)=Q^1(\sigma) Q^{2,1}(\sigma)$ in $ H_{6,4}$.
    \item $Q^1(\sigma)Q^3(\sigma) = \sigma^2 Q^4(\sigma) $ in $H_{4,4}$.
 \item $\sigma^3 Q^4(\sigma)=0 $ in $H_{5,4}$.
\end{enumerate}    
\end{cor}

Throughout this section $\GL_g\defeq \GL_g(\F_2)$, and (co)homologies are taken with $\F_2$ coefficients.

\subsection{Proof of \cref{thm: additive basis} }\label{sec: grouphom 1} 
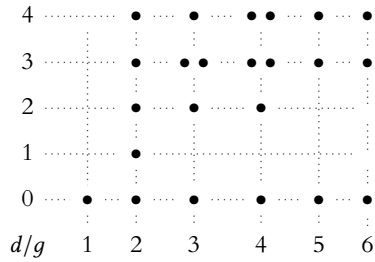
\begin{figure}[h]
    \centering

\adjustbox{scale =0.9,center}{\begin{tikzcd}[sep = tiny]
	4 & {} & \bullet & \bullet & {\bullet\;\bullet} & \bullet & \bullet \\
	3 && \bullet & {\bullet\;\bullet} & {\bullet\;\bullet} & \bullet & \bullet \\
	2 && \bullet & \bullet & \bullet && {} \\
	1 && \bullet &&&& {} \\
	0 & \bullet & \bullet & \bullet & \bullet & \bullet & \bullet \\
	{d/g} & 1 & 2 & 3 & 4 & 5 & 6
	\arrow[draw={rgb,255:red,87;green,87;blue,87}, dotted, no head, from=1-1, to=1-3]
	\arrow[draw={rgb,255:red,87;green,87;blue,87}, dotted, no head, from=1-2, to=5-2]
	\arrow[draw={rgb,255:red,87;green,87;blue,87}, dotted, no head, from=1-3, to=1-4]
	\arrow[draw={rgb,255:red,87;green,87;blue,87}, dotted, no head, from=1-3, to=2-3]
	\arrow[draw={rgb,255:red,87;green,87;blue,87}, dotted, no head, from=1-4, to=1-5]
	\arrow[draw={rgb,255:red,87;green,87;blue,87}, dotted, no head, from=1-4, to=2-4]
	\arrow[draw={rgb,255:red,87;green,87;blue,87}, dotted, no head, from=1-5, to=1-6]
	\arrow[draw={rgb,255:red,87;green,87;blue,87}, dotted, no head, from=1-5, to=2-5]
	\arrow[draw={rgb,255:red,87;green,87;blue,87}, dotted, no head, from=1-6, to=1-7]
	\arrow[draw={rgb,255:red,87;green,87;blue,87}, dotted, no head, from=1-6, to=2-6]
	\arrow[draw={rgb,255:red,87;green,87;blue,87}, dotted, no head, from=1-7, to=2-7]
	\arrow[draw={rgb,255:red,87;green,87;blue,87}, dotted, no head, from=2-1, to=2-3]
	\arrow[draw={rgb,255:red,87;green,87;blue,87}, dotted, no head, from=2-3, to=2-4]
	\arrow[draw={rgb,255:red,87;green,87;blue,87}, dotted, no head, from=2-3, to=3-3]
	\arrow[draw={rgb,255:red,87;green,87;blue,87}, dotted, no head, from=2-4, to=2-5]
	\arrow[draw={rgb,255:red,87;green,87;blue,87}, dotted, no head, from=2-4, to=3-4]
	\arrow[draw={rgb,255:red,87;green,87;blue,87}, dotted, no head, from=2-5, to=2-6]
	\arrow[draw={rgb,255:red,87;green,87;blue,87}, dotted, no head, from=2-5, to=3-5]
	\arrow[draw={rgb,255:red,87;green,87;blue,87}, dotted, no head, from=2-6, to=2-7]
	\arrow[draw={rgb,255:red,87;green,87;blue,87}, dotted, no head, from=2-6, to=5-6]
	\arrow[draw={rgb,255:red,87;green,87;blue,87}, dotted, no head, from=2-7, to=3-7]
	\arrow[draw={rgb,255:red,87;green,87;blue,87}, dotted, no head, from=3-1, to=3-3]
	\arrow[draw={rgb,255:red,87;green,87;blue,87}, dotted, no head, from=3-3, to=3-4]
	\arrow[draw={rgb,255:red,87;green,87;blue,87}, dotted, no head, from=3-3, to=4-3]
	\arrow[draw={rgb,255:red,87;green,87;blue,87}, dotted, no head, from=3-4, to=3-5]
	\arrow[draw={rgb,255:red,87;green,87;blue,87}, dotted, no head, from=3-4, to=5-4]
	\arrow[draw={rgb,255:red,87;green,87;blue,87}, dotted, no head, from=3-5, to=3-7]
	\arrow[draw={rgb,255:red,87;green,87;blue,87}, dotted, no head, from=3-5, to=5-5]
	\arrow[draw={rgb,255:red,87;green,87;blue,87}, dotted, no head, from=3-7, to=4-7]
	\arrow[draw={rgb,255:red,87;green,87;blue,87}, dotted, no head, from=4-1, to=4-3]
	\arrow[draw={rgb,255:red,87;green,87;blue,87}, dotted, no head, from=4-3, to=4-7]
	\arrow[draw={rgb,255:red,87;green,87;blue,87}, dotted, no head, from=4-3, to=5-3]
	\arrow[draw={rgb,255:red,87;green,87;blue,87}, dotted, no head, from=4-7, to=5-7]
	\arrow[draw={rgb,255:red,87;green,87;blue,87}, dotted, no head, from=5-1, to=5-2]
	\arrow[draw={rgb,255:red,87;green,87;blue,87}, dotted, no head, from=5-2, to=5-3]
	\arrow[draw={rgb,255:red,87;green,87;blue,87}, dotted, no head, from=5-2, to=6-2]
	\arrow[draw={rgb,255:red,87;green,87;blue,87}, dotted, no head, from=5-3, to=5-4]
	\arrow[draw={rgb,255:red,87;green,87;blue,87}, dotted, no head, from=5-3, to=6-3]
	\arrow[draw={rgb,255:red,87;green,87;blue,87}, dotted, no head, from=5-4, to=5-5]
	\arrow[draw={rgb,255:red,87;green,87;blue,87}, dotted, no head, from=5-4, to=6-4]
	\arrow[draw={rgb,255:red,87;green,87;blue,87}, dotted, no head, from=5-5, to=5-6]
	\arrow[draw={rgb,255:red,87;green,87;blue,87}, dotted, no head, from=5-5, to=6-5]
	\arrow[draw={rgb,255:red,87;green,87;blue,87}, dotted, no head, from=5-6, to=5-7]
	\arrow[draw={rgb,255:red,87;green,87;blue,87}, dotted, no head, from=5-6, to=6-6]
	\arrow[draw={rgb,255:red,87;green,87;blue,87}, dotted, no head, from=6-7, to=5-7]
\end{tikzcd}}
    \caption{$\dim H_d(\GL_g)=$ number of dots at $(g,d)$. }
    \label{fig: dimension}
\end{figure}

For $g=1$, there is only $\sigma$, the generator of $H_0$ of the trivial group.
The $g=2$ column was computed in \autocite[Lemma 6.3(ii)]{GKRW18b}; we will recall the argument below for completeness. For the $g\geqslant 3$ columns, we rely on information of the  dimension of $H_{g,d}$ in this range, displayed in  \cref{fig: dimension}. This is obtained using the \autocite[]{HAP} package of  \autocite[]{GAP4}.

\begin{rem}
The $g=3$ part of  \cref{fig: dimension} can be obtained from the existing computations of $H^*(\GL_3)$ (\autocite[Theorem 11.7]{BC87}, \autocite[1971-11: August 18, pp.63-66]{QuiN}). For the $g=4$ part, one can use the classical computation of $H^*(A_8)$ and the Jordan-Dickson isomorphism $\GL_4\cong A_8$, our reference for these are \autocite[Corollary 6.5 and pp.211]{AM04}.  The $g= 5,6$ part have not seem to be recorded in the literature.
\end{rem}

\subsubsection{$g=2$.}\label{sec: group homology g=2}
\begin{prop}\label{prop: description of GL_2}
    For $i\geqslant 0$, $H_{2,i}$ is $1$-dimensional spanned by $Q^i(\sigma)$.
\end{prop}
\begin{prf}
    The grading-$2$ piece of the $\E_\infty$-map from the free $E_\infty$-space on a point in grading $1$ to the $\E_\infty$-space associated to the $\GL_g$ is given by including $\Sigma_2\hookrightarrow\GL_2$ as permutation matrices. This is a homology isomorphism because $\GL_2$ can be decomposed as a semidirect product of image of $\Sigma_2$ with a cyclic subgroup of order $3$, and $\widetilde{H}_*(C_3)=0$. By the description of homology groups of free algebras, we know $H_i(\Sigma_2) $ is generated by $Q^i(\sigma)$. 
\end{prf}

\subsubsection{The map $\GL_3\to \GL_4$.}\label{sec: GL3 to GL4}Our  goal is the following proposition, which we will use later to describe the $g=3$ and $g=4$ columns of \cref{fig: additive basis}.
\begin{prop}\label{prop: GL3 to GL4}
The stabilisation homomorphism $H_i(\GL_3)\to H_i(\GL_4)$ is
\begin{enumerate}[(i)]
    \item  the zero map when $i=2,3$, and 
    \item  non-zero (and hence injective) when $i=4$.
\end{enumerate} 
\end{prop}
The strategy is to study the linear dual, i.e. the restriction map $H^i(\GL_4)\to H^i(\GL_3)$. In \autocite[]{Pha87} Pham showed that the the restriction map $H^*(\GL_4)\to \bigoplus_{i=1}^5H^*(A_i)$ to a certain collection of elementary abelian subgroups $A_1,\dots, A_5$ (described in (\ref{eq: Ai}) below) is injective \autocite[(3.1)]{Pha87}, and they described a set of algebra generators of $H^*(\GL_4)$ in terms of their image in the $A_i$. We will describe elementary abelian subgroups $B_1,B_2$ of $\GL_3$ with analogous property. The $B_i$ are also arranged so that the image of $B_1$ (resp. $B_2$) in $\GL_4$ will be a subgroup of $A_1$ (resp. $A_2$), thus the vanishing of Pham's classes in $\GL_3$ can be explicitly checked on the $B_i$.

\begin{notation}
\begin{enumerate}[before*=\leavevmode\vspace{-0.17\baselineskip}]
    \item   For $n>1$ and $1\leqslant i,j\leqslant n$, $i\neq j$, let $a_{ij}^{(n)}\in \GL_n$ denote the matrix with $1$'s along the diagonal, $1$ in position $(i,j)$, and $0$ elsewhere. The image of $a_{ij}^{(n)}$ in $\GL_{n+1}$ is $a_{ij}^{(n+1)}$, so we can write $a_{ij}$ for $a_{ij}^{(n)}$. 
    \item Suppose $S$ is a collection of $a_{ij}$ in $\GL_n$ such that  the subgroup $A$ generated by $S$ is elementary abelian, then the $a_{ij}$ form a basis of $H_1(A) =A$ and we define $\alpha_{ij}\in \mathrm{Hom}_{\F_2}(A,\F_2) =H^1(A)$ to be the dual basis.
    \item Let $\UT_n\leqslant \GL_n$  denote the subgroups of upper triangular matrices. Since $\UT_n$ is a Sylow $2$-subgroup of $\GL_n$, the restriction map  on cohomology is injective.
\end{enumerate}
  
\end{notation}

\begin{prop}\label{prop: UT3}
    Consider the elementary abelian subgroups \[B_1\defeq \<a_{12}, a_{13}>,\quad B_2\defeq \<a_{13}, a_{23}>\] of $\UT_3$. Then
    \begin{enumerate}[(i)]
        \item The restriction homomorphism $r: H^*(\UT_3)\to H^*(B_1)\oplus  H^*(B_2)$ is injective. In particular, so is the restriction homomorphism $H^*(\GL_3)\to H^*(B_1)\oplus H^*(B_2)$.
        \item Let $\rho_3$ denote the $7$-dimensional real permutation representation of $\GL_3$ given by the natural action of $\GL_3$ on the set $\F_2^{3}\setminus\{0\}$. Let $b = w_2(\rho_3)|_{\UT_3}\in H^2(\UT_3)$ be the restriction of the second Stiefel-Whitney class of $\rho_3$. Then $r(b) = (\alpha_{13}^2+\alpha_{12}^2,\alpha_{13}^2+\alpha_{23}^2)$.
        \item There are classes  $a,c\in H^1(\UT_3)$ satisfying $r(a) =(0, \alpha_{23})$, $r(c) = (\alpha_{12},0)$. By (i), $a,c$ are characterised by their image under $r$ and we have $ac = 0$. 
        \item The map of graded rings $\F_2[a,b,c]/(ac)\to H^*(\UT_3)$ is an isomorphism.
    \end{enumerate}
\end{prop}

\begin{prf}
Let $\overline{a}_{12},\overline{a}_{23}$ be the image of $a_{12},a_{23}$ in the $\UT_3/\<a_{13}>$. Then $\UT_3/\<a_{13}>$ is elementary abelian,  generated by $\overline{a}_{12},\overline{a}_{23}$. Also $B_2/\<a_{13}> =\<\overline{a_{23}}>$, $B_1/\<a_{13}> = \<\overline{a_{12}}>$, and we have the following diagram of extensions\stepcounter{thm}
    \[\begin{tikzcd}
	{(\thethm_2)} & 0 & {\<a_{13}>} & {B_2} & {\<\overline{a_{23}}>} & 0 \\
	{(\thethm)} & 0 & {\<a_{13}>} & {\UT_3} & {\<\overline{a_{12}},\overline{a_{23}}>} & 0 \\
	{(\thethm_1)} & 0 & {\<a_{13}>} & {B_1} & {\<\overline{a_{12}}>} & 0
	\arrow[from=1-2, to=1-3]
	\arrow[from=1-3, to=1-4]
	\arrow[equals, from=1-3, to=2-3]
	\arrow[from=1-4, to=1-5]
	\arrow["{f_2}"', hook, from=1-4, to=2-4]
	\arrow[from=1-5, to=1-6]
	\arrow["{g_2}", hook, from=1-5, to=2-5]
	\arrow[from=2-2, to=2-3]
	\arrow[from=2-3, to=2-4]
	\arrow[from=2-4, to=2-5]
	\arrow[from=2-5, to=2-6]
	\arrow[from=3-2, to=3-3]
	\arrow[equals, from=3-3, to=2-3]
	\arrow[from=3-3, to=3-4]
	\arrow["{f_1}", hook, from=3-4, to=2-4]
	\arrow[from=3-4, to=3-5]
	\arrow["{g_1}"', hook, from=3-5, to=2-5]
	\arrow[from=3-5, to=3-6]
\end{tikzcd}\]
induced by inclusion of subgroups. The $E_2$-page of the associated spectral sequences within $\{(p,q)\mid 0\leqslant p+q\leqslant 2\}$ is as follows.
\begin{center}
\adjustbox{scale=0.9}{\begin{tikzcd}[sep=small]
	{\alpha_{13}^2} &&& {} & {\alpha_{13}^2} &&& {} & {\alpha_{13}^2} \\
	{\alpha_{13}} & {\alpha_{13}\overline{\alpha_{23}}} &&& {\alpha_{13}} & {\alpha_{13}\overline{\alpha_{12}},\alpha_{13}\overline{\alpha_{23}}} &&& {\alpha_{13}} & {\alpha_{13}\overline{\alpha_{12}}} \\
	1 & {\overline{\alpha_{23}}} & {\overline{\alpha_{23}}^2} && 1 & {\overline{\alpha_{12}},\overline{\alpha_{23}}} & {\overline{\alpha_{12}}^2,\overline{\alpha_{23}}^2,\overline{\alpha_{12}}\overline{\alpha_{23}}} && 1 & {\overline{\alpha_{12}}} & {\overline{\alpha_{12}}^2} \\
	& {E_2^{**}(B_2)} && {} && {E_2^{**}(\UT_3)} && {} && {E_2^{**}(B_1)}
	\arrow[dotted, no head, from=1-1, to=2-1]
	\arrow[dotted, no head, from=1-5, to=2-5]
	\arrow[no head, from=1-8, to=4-8]
	\arrow[dotted, no head, from=1-9, to=2-9]
	\arrow[dotted, no head, from=2-1, to=2-2]
	\arrow[dotted, no head, from=2-1, to=3-1]
	\arrow[dotted, no head, from=2-2, to=3-2]
	\arrow[dotted, no head, from=2-5, to=2-6]
	\arrow["{d^2}"{description}, from=2-5, to=3-7]
	\arrow[dotted, no head, from=2-6, to=3-6]
	\arrow[dotted, no head, from=2-9, to=2-10]
	\arrow[dotted, no head, from=2-9, to=3-9]
	\arrow[dotted, no head, from=2-10, to=3-10]
	\arrow[dotted, no head, from=3-1, to=3-2]
	\arrow[dotted, no head, from=3-2, to=3-3]
	\arrow[dotted, no head, from=3-5, to=2-5]
	\arrow[dotted, no head, from=3-5, to=3-6]
	\arrow[dotted, no head, from=3-6, to=3-7]
	\arrow[dotted, no head, from=3-9, to=3-10]
	\arrow[dotted, no head, from=3-10, to=3-11]
	\arrow[no head, from=4-4, to=1-4]
\end{tikzcd}}
\end{center}
where $\overline{\alpha_{12}}, \overline{\alpha_{23}}\in H^1(\<\overline{a_{12}},\overline{a_{13}}>)$ is defined to be the dual basis of $\overline{a}_{12},\overline{a_{13}}$, and similarly $\overline{\alpha_{23}}\in H_1(\<\overline{\alpha}_{23}>)$, $\overline{\alpha_{12}}\in H_1(\<\overline{a_{12}}>)$. Moreover, the map $f_2^*:E_2^{**}(\UT_3)\to E_2^{**}(B_2)$ sends $\overline{\alpha_{12}}$ to $0$ and the map $f_1^*:E_2^{**}(\UT_3)\to E^2_{**}(B_1)$ sends $\overline{\alpha_{23}}$ to $0$. The $d^2$-differential is described by the following
\begin{claim}
    $d^2(\alpha_{13}) =\overline{\alpha_{12}}\overline{\alpha_{23}}\in E_2^{3,0}(\UT_3)$.
\end{claim}

\begin{prfoc} The class $d^2(\alpha_{13})\in H_3(\<\overline{a_{12}},\overline{a_{23}}>)$ classifies the central extension ($\thethm$). If $y$ is one of $\overline{a_{12}}$, $\overline{a_{23}}$, or $\overline{a_{12}}\overline{a_{23}}$, then the restriction  of  $d^2(\alpha_{13})$ to $\<\overline{y}>$ corresponds to the {pullback extension}  below.
\[\begin{tikzcd}[cramped]
	0 & {\<a_{13}>} & {\UT_3} & {\<\overline{a_{12}},\overline{a_{23}}>} & 0 \\
	0 & {\<a_{13}>} & {H} & {\<y>} & 0
	\arrow[from=1-1, to=1-2]
	\arrow[from=1-2, to=1-3]
	\arrow[from=1-3, to=1-4]
	\arrow[from=1-4, to=1-5]
	\arrow[from=2-1, to=2-2]
	\arrow[equals, from=2-2, to=1-2]
	\arrow[from=2-2, to=2-3]
	\arrow[from=2-3, to=1-3]
	\arrow["\lrcorner"{anchor=center, pos=0.125, rotate=90}, draw=none, from=2-3, to=1-4]
	\arrow[from=2-3, to=2-4]
	\arrow[from=2-4, to=1-4]
	\arrow[from=2-4, to=2-5]
\end{tikzcd}\]
When $y = \overline{a_{12}}$ (resp. $\overline{a_{23}}$) this extension is ($\thethm_1$) (resp. $(\thethm_2)$) which splits. When $y = \overline{a_{12}a_{23}}$ the extension does not split, for $H\simeq C_4$. As a consequence, if we write $d^2(\alpha_{12}) = \lambda_1\overline{\alpha_{12}}^2 + \lambda_2\overline{\alpha_{23}}^2 + \lambda_3 \overline{\alpha_{12}}\overline{\alpha_{23}}$ then we must have  $\lambda_1 =\lambda_2 =0$ and $\lambda_3\neq 0$.\end{prfoc}

Now $E_3^{*,*}(\UT_3) =\F_2[\overline{\alpha_{12}},\overline{\alpha_{23}}]/(\overline{\alpha_{12}\alpha_{23}})\otimes \F_2[\alpha_{13}^2]$ and the $d^3$-differential is determined by $d^3(\alpha_{13}^2)$. One can write $d^3(\alpha_{13}^2) = \lambda_1\overline{\alpha_{12}}^3 + \lambda_2\overline{\alpha_{23}}^2$ for some $\lambda_1,\lambda_2\in \F_2$. By naturality we have $d^3(\alpha_{13}^2) =\lambda_i\overline{\alpha_{i3}}^2\in E_3^{3,0}(B_i)$. Since ($\thethm_i$) splits as a direct product, the spectral sequence collapses at $E_2$. Thus $\lambda_1 =\lambda_2=0$. Consequently, $\{E_r(\UT_3)\}_r$ collapses at $E_3$ and $(f_1^*,f_2^*):E_\infty^{*,*}(\UT_3)\to E_\infty^{*,*}(B_1)\oplus E_\infty^{*,*}(B_2)$ is 
\[\overline{r}:\F_2[\overline{\alpha_{12}},\overline{\alpha_{23}}]/(\overline{\alpha_{12}\alpha_{23}})\otimes \F_2[\alpha_{13}^2]\to \F_2[\overline{\alpha_{12}},\alpha_{13}]\oplus\F_2[\overline{\alpha_{23}},\alpha_{13}]. \]
Since $\overline{r}$ is the map on filtration quotients of a bounded-below filtration on $r$, the fact that $\overline{r}$ is injective implies $r$ is so too, proving (i).

For (ii), it as computed in \autocite[the proof of Lemma 6.3. (ii)]{GKRW18b} that the $\rho_3|_{\<a_{12}>}$ consists of $2$ copies of the sign representation and $5$ copies of the trivial representation, so $b|_{\<a_{12}>} = \alpha_{12}^2$. Since the matrices $a_{13}$ and $a_{23}$ are conjugate to $a_{12}$ in $\GL_3$ and $b$ is restricted from $\GL_3$, the restriction of $b$ to $\<a_{13}>$ and $\<a_{23}>$ are also non-zero, so we must have $b|_{\<a_{13}>} = \alpha_{13}^2$ and $b|_{\<a_{23}>} = \alpha_{23}^2$. Now restricting $r(b) =(\lambda_{1}\alpha_{13}^2 +\lambda_2\alpha_{12}^2, \mu_1\alpha_{13}^2 +\mu_2\alpha_{23}^2)$ for some coefficients $\lambda_i,\mu_j$ to those order-2 subgroups we conclude $\lambda_1 =\lambda_2 =\mu_1 =\mu_2=1$. 

For (iii), we identify $E_\infty^{1,0}(\UT_3)$ as the bottom piece of the filtration on $H^*(\UT_3)$ and take $a,c$ to be $\overline{\alpha_{23}},\overline{\alpha_{12}}$.

Finally, for (iv) it suffices to check the composite \[\F_2[a,b,c]/(ac)\to H^*(\UT_3)\to E_\infty^{*,*}(\UT_3) = \F_2[\overline{\alpha_{12}},\overline{\alpha_{23}}]/(\overline{\alpha_{12}\alpha_{23}})\otimes \F_2[\alpha_{13}^2]\] is an isomorphism. This is indeed the case as $a$ (resp. $c$) is sent  to $\overline{\alpha_{23}}$ (resp. $\overline{\alpha_{12}}$) by (iii) and $b$ is sent to $\alpha_{13}^2 + \overline{\alpha_{12}}^2 + \overline{\alpha_{23}}^2$ by (ii).
\end{prf} 

\begin{prf}[of \cref{prop: GL3 to GL4}]
    Consider the following elementary abelian subgroups of $\GL_4$: \stepcounter{thm}
    \[\tag{\thethm}\label{eq: Ai}\leqnomode A_1\defeq \<a_{12},a_{13},a_{14}>, \quad A_2\defeq \<a_{23}, a_{13}, a_{24},a_{14}>,\quad A_3\defeq \<a_{34},a_{24},a_{14}>,\quad A_4\defeq \<a_{12},a_{34},a_{14}>,\]\[A_5\defeq \<a_1\defeq a_{12}a_{34}, a_2\defeq a_{13}a_{24},a_3\defeq a_{14}>.\] Pham showed that the $A_i$ detect zero for $\GL_4$ \autocite[Lemma 3.1]{Pha87} and they  constructed algebra generators\footnote{There is a typo in \autocite[Theorem A]{Pha87}: we should have $\abs{v_{14}}=4$ (instead of $3$), c.f. the restriction  of $v_{14}$ to the $A_i$ on loc.cit. page 58. There is also a misprint on the definition of $a_1$ (called \say{$c_1$} there).} $w_1,\dots, w_{10}$ of the subalgebra $H^*(\GL_4)$ of $H^*(\UT_4)$ by specifying their image in the $A_i$. Since the other $w_j$ lives beyond degree $4$, will only use the restriction of $w_1,w_2,w_3,w_4$ to $A_1$ and $A_2$, as recorded below. \rowcolors{0}{}{gray!5}
\[\begin{tabular}{l|c|c}
element / its restriction to  & $A_1$ & $A_2$\\
   $w_1$ &  $0$ &  $\alpha_{23}\alpha_{14}+\alpha_{13}\alpha_{24}$ \\
    $w_2$  &$0$ & $\alpha_{23}^2\alpha_{14}+\alpha_{23}\alpha_{14}^2+\alpha_{13}^2\alpha_{24}+\alpha_{13}\alpha_{24}^2$\\
     $w_3$ & $0$ & $0$\\
     $w_4$  & $x_1$  &$x_2$
\end{tabular}\]
where
\begin{align*}
    x_1&=\alpha_{14}(\alpha_{14}+\alpha_{13})(\alpha_{14}+\alpha_{12})(\alpha_{14}+\alpha_{13}+\alpha_{12}) +\alpha_{12}^4 + (\alpha_{13}^2+\alpha_{13}\alpha_{12})^2\\
    x_2&= \alpha_{14}(\alpha_{14}+\alpha_{13})(\alpha_{14}+\alpha_{24})(\alpha_{14}+\alpha_{13}+\alpha_{24}+\alpha_{23}) + (\alpha_{23}\alpha_{14}+\alpha_{13}\alpha_{24})((\alpha_{13}^2+\alpha_{13}\alpha_{23})+(\alpha_{24}^2+\alpha_{24}\alpha_{23}) + \alpha_{23}^2) \\
    &+ (\alpha_{13}^2+\alpha_{13}\alpha_{23})^2  + (\alpha_{24}^2+\alpha_{24} \alpha_{23})^2 +\alpha_{23}^4.
\end{align*}
In particular, the images of the $w_i$ in $B_1$ and $B_2$ are as follows.
\begin{enumerate}[(1)]
    \item For $1\leqslant i\leqslant 3$,  we have $w_i|_{B_1} =w_i|_{B_2}= 0$.
    \item $w_4|_{B_1} =\alpha_{12}^4 + \alpha_{12}^2\alpha_{13}^2+\alpha_{13}^4  $ and $w_4|_{B_2} =\alpha_{13}^4 + \alpha_{13}^2\alpha_{23}^2+\alpha_{23}^4$; both are non-zero.
\end{enumerate}
Now we are ready to prove part (i) of the proposition.  Since $\abs{w_1} =2$, $\abs{w_2}=\abs{w_3}=3$ and all other $w_j$ have degree at least $4$, we must have $H^2(\GL_4) =\F_2\{w_1\}$ and $H^3(\GL_4)=\F_2\{w_2\}\oplus \F_2\{w_3\}$. By \cref{prop: UT3} (i), (1) implies  $w_i|_{\GL_3}=0$ for $i =1,2,3$. 

Part (ii) follows from (2) as the class $w_4$ lives in degree $4$. 
\end{prf}

\subsubsection{$g=3$.} Let $a,b,c$ be the classes from \cref{prop: UT3}. Then

\begin{prop}\label{prop: subalg GL3}
    As a subalgebra of $H^*(\UT_3)$, $H^*(\GL_3)$ is generated by $b$, $a(b+a^2)$ and $c(b+c^2)$
\end{prop}
\begin{prf}
    We import the fact that  $H^*(\GL_3)$ is generated as an algebra in degrees $2$ and $3$. Indeed, the algebra structure $H^*(\GL_3) = \F_2[\xi,\gamma_1,\gamma_2]/(\gamma_1\gamma_2)$ with $\abs{\xi}=2$, $\abs{\gamma_1}=\abs{\gamma_2}=3$ has been computed several places: for instance \autocite[Theorem 11.7]{BC87} and \autocite[1971-11: August 18, pp.63-66]{QuiN}. Using this we are reduced to checking
    \begin{enumerate}
        \item $H^2(\GL_3) = \F_2\{b\}$.
        \item $H^3(\GL_3) = \F_2\{a(b+a^2)\}\oplus \F_2\{c(b+c^2)\}$.
    \end{enumerate}
    Since the subgroups $\<a_{12}>,\<a_{13}>$ and $\<a_{23}>$ are conjugate in $\GL_3$, the restriction of a class in $H^*(\GL_3)$ to one of these subgroups is zero iff the restriction to all of them is. Now given 
    \[x\defeq \lambda_1a^2 + \lambda_2c^2 + \lambda_3b\in H^2(\UT_3) = \F_2\{a^2\}\oplus \F_2\{c^2\}\oplus \F_2\{b\},\quad \lambda_1,\lambda_2,\lambda_3\in \F_2,\]
    we have 
    \[x|_{\<a_{12}>} = (\lambda_2+\lambda_3)\alpha_{12}^2,\quad x|_{\<a_{13}>} = \lambda_3 \alpha_{13}^2,\quad x|_{\<a_{23}>} = (\lambda_1+\lambda_3)\alpha_{23}^2\]
    by \cref{prop: UT3}. So if $x\in H^2(\GL_3)$ then either $\lambda_1 =\lambda_2 =\lambda_3 =0$ or $\lambda_1=\lambda_2=0, \lambda_3=1$. As we know $H^2(\GL_3)$ is $1$-dimensional, the  proves (1).

    For (ii), first note that $B\defeq \{ab, bc,a^3,c^3\}\subset H^3(\UT_3)$ is a basis. By \cref{prop: UT3}(iii), every element in $B$ restricts to zero in $\<a_{13}>$, and hence every element in $B$ restricts to zero in $\<a_{12}>$ and $\<a_{23}>$ as well. Given 
    \[y \defeq \lambda_1ab + \lambda_2 bc + \lambda_3a^3 + \lambda_4c^3\in H^3(\UT_3),\]
    we have \[y|_{\<a_{12}>} =(\lambda_2+\lambda_4)\alpha_{12}^3 ,\quad y|_{\<a_{23}>} = (\lambda_1+\lambda_3)\alpha_{23}^3.\]
    Thus $H^3(\GL_3)\subseteq S\defeq  \{y\text{ as above}\mid \lambda_2=\lambda_4,\lambda_1=\lambda_3\}$. But both $H^3(\GL_3)$ and $S$ have four elements, so in fact $H^3(\GL_3) = S$, proving (ii).
\end{prf}

\begin{cor}\label{cor: g=3} The map
$H_i(\GL_2)\to H_i(\GL_3) $ is an isomorphism when $i$ is even, and zero  when $i$ is odd. In particular:
    \begin{enumerate}[(i)]
    \item For every even $i\geqslant 2$, $H_{3,i}$ is $1$-dimensional and spanned by $\sigma Q^i(\sigma)$. 
    \item For every odd $i\geqslant 1$, $\sigma Q^i(\sigma)=0\in H_{3,i}$. 
    \end{enumerate}
\end{cor}
\begin{prf}
   Suppose $i =2k>0$. Using the algebra structure of \cref{prop: subalg GL3} one sees $H^{2k} (\GL_3) = \F_2\{b^k\}$.  In particular $H_{3,2k}$ is $1$-dimensional. Since $H_{2,2k} = \F_2\{Q^1(\sigma)\}$, by dualising it suffices to show the restriction map $H^{2k}(\GL_3)\to H^{2k}(\GL_2)$ is non-zero, which we do by further restricting to the subgroup $\<a_{12}>$. By \cref{prop: UT3}(ii),  $b^{2k}|_{\<a_{12}>} =\alpha_{12}^{4k}$, which is non-zero.

   Suppose $i=2k+1>0$. By dualising it is enough to show $H^{2k+1}(\GL_3)\to H^{2k+1}(\GL_2)$ is the zero map. When $i=3$ this is just because $H^3(\GL_3)$ is the zero group, so we assume $i>3$. Then every non-zero class $x\in H^{2k+1}(\GL_3)$ is divisible by at least one of $a(b+a^2)$ and $c(b+c^2)$, so $x|_{\<a_{13}>}=0$ by \cref{prop: UT3}(iii). Since $\<a_{13}>$ is conjugate to $\<a_{12}>$ in $\GL_3$ we have $x|_{\<a_{12}>}=0$. As the restriction $H^*(\GL_3)\to H^*(\<a_{12}>)$ is an isomorphism, this completes the proof.
\end{prf}

\subsubsection{$g=4$.}

Using $\dim H_{4,2}=1$, $\dim H_{4,3} = \dim H_{4,4}=2$ from (\cref{fig: dimension}), we reformulate the $g=4$ column of \cref{thm: additive basis} as the following
\begin{prop}\label{prop: pi4*}
\begin{enumerate}[(i),before*=\leavevmode\vspace{-0.17\baselineskip}]
        \item $Q^1(\sigma)^2\neq0 $  in $H_{4,2}$;
        \item $Q^{2,1}(\sigma)$ and $Q^1(\sigma) Q^2(\sigma)$ are non-zero and distinct in $H_{4,3}$;
        \item $Q^2(\sigma)^2$ and $\sigma^2 Q^4(\sigma)$ are  non-zero and distinct in $H_{4,4}$. 
    \end{enumerate}
\end{prop}
\begin{prf}
Item (i) is  \autocite[Lemma 6.3 (iii)]{GKRW18b}.

For (ii), we consider the linear dual $\Sq^i_*:H_{*,*}\to H_{*,*-i}$  of the Steenrod operation. The interaction between the dual Steenrod operation and the Dyer-Lashof operations are described by the Nishida relation \autocite[6]{CLM76}. Here are a few instances of the relation that we will use.
\begin{enumerate}[(N1)]
   \item $\Sq^1_*Q^1=0$.
    \item $\Sq^1_*Q^2 = Q^1$.
 \item $\Sq^2_*Q^4 = Q^2+Q^3\Sq_1^*$.

\end{enumerate}
By (N2)
\[\Sq^1_* (Q^{2,1}(\sigma) )= Q^1(Q^1(\sigma)) = Q^1(\sigma)^2.\]
On the other hand \[\Sq^1_*(Q^1(\sigma)Q^2(\sigma))=Q^1(\sigma) \Sq^1_*(Q^2(\sigma))+\Sq^1_*( Q^1(\sigma)) Q^2(\sigma) = Q^1(\sigma)^2\]
where the first equality is by the Cartan formula and the second equality is by  (N2) and (N1).  Since $Q^1(\sigma)^2\in H_{4,2}$ is non-zero, neither are $Q^{2,1}(\sigma)$ and $Q^1(\sigma)Q^2(\sigma)$. Those two classes are distinct because $\sigma Q^1(\sigma)Q^2(\sigma)=0$ while we will see in \cref{prop: g=5}(i) below that $\sigma Q^{2,1}(\sigma)\neq 0$.

For (iii), the class $\sigma^2Q^4(\sigma)$ is non-zero because $\sigma:H_{3,4}\to H_{4,4}$ is injective by \cref{prop: GL3 to GL4}(ii) and $\sigma Q^4(\sigma)$ is non-zero by \cref{cor: g=3}(iii). To see $Q^2(\sigma)^2$ is non-zero, one computes
    \[ \Sq^2_*(Q^2(\sigma)^2)=Q^2(\sigma) \Sq^2_*(Q^2(\sigma)) + (\Sq^1_*(Q^2(\sigma)))^2 +  \Sq^2_*(Q^2(\sigma)) Q^2(\sigma) 
    = Q^1(\sigma)^2.\]
    The first equality is by the  Cartan formula, the second equality uses (N2) in the middle term. To see that  $\sigma^2Q^4(\sigma)$ and $ Q^2(\sigma)^2$ are distinct, we evaluate $\Sq^2_*$ on them. 
    \[\Sq^2_*(\sigma^2 Q^4(\sigma)) = \sigma^2 \Sq^2_*(Q^4(\sigma)) + \Sq^1_*(\sigma^2)\Sq^1_*(Q^4(\sigma)) + \Sq^2_*(\sigma^2)Q^4(\sigma) = \sigma^2 Q^2(\sigma)\]
    by the Cartan formula and then (N3). By \cref{prop: GL3 to GL4}(i), the class $\sigma^2 Q^2(\sigma)$ is zero.  On the other hand 
    $\Sq^2_*(Q^2(\sigma)^2) =(\Sq^1_*Q^2(\sigma))^2 = Q^1(\sigma)^2\ne 0.$
    \end{prf}

\subsubsection{$g=5$.}\label{sec: group homology g=5} We need to show 
\begin{prop}\label{prop: g=5}
    \begin{enumerate}[(i),before*=\leavevmode\vspace{-0.17\baselineskip}]
        \item $\sigma Q^{2,1}(\sigma)\in H_{5,3}$ is non-zero.
        \item $\sigma Q^2(\sigma)^2\in H_{5,4}$ is non-zero.
    \end{enumerate}
\end{prop}
\begin{prf}
Part (i) follows immediately from \cref{prop: group homology g=6}(i) below, which says the stabilisation of this class is non-zero. 

Part (ii) can be obtained using the Stiefel-Whitney class trick of \autocite[Lemma 6.3]{GKRW18b}. Let $\rho_5 $
be the $32$-dimensional real permutation representation given by the natural action of $\GL_5$ on the set $\F_2^{5}$. Let $i$ be the composite $ \<a_{12},a_{34}>\hookrightarrow\GL_4\to \GL_5$. Then
\[i^*\rho_5 =(\rho_4|_{\<a_{12},a_{34}>})^{\oplus2}\]
where $\rho_4$ is the $16$-dimensional representation obtained via the $\GL_4$-action on $\F_2^4$. The first few  Stiefel-Whitney class of $\rho_4|_{\<a_{12},a_{34}>}$ 
 was computed in  \autocite[proof of 6.3 (iv)]{GKRW18b}\footnote{The $\rho_4$ there differs by the one here by a copy of trivial representation, in particular the Stiefel-Whitney class is unaffected}, giving us
 \[w(i^*\rho_5) =w(\rho_4|_{\<a_{12},a_{34}>})^2 =1+\alpha_{12}^2\alpha_{34}^2 + (\text{classes of degree}\geqslant 6).\]
 In particular $i^*w_4(\rho_5) = \alpha_{12}^2 \alpha_{34}^2$, which pairs to $1$ with $Q^2(\sigma)\times Q^2(\sigma)$. Consequently, $w_4(\rho_5)$ pairs to $1$ with $i_*(Q^2(\sigma)\times Q^2(\sigma))=\sigma Q^2(\sigma)^2$.
\end{prf}

\subsubsection{$g=6$.}\label{sec: group homology g=6} By the dimension information it's enough to check
\begin{prop}\label{prop: group homology g=6}
    \begin{enumerate}[(i),before*=\leavevmode\vspace{-0.17\baselineskip}]
    \item $\sigma^2 Q^{2,1}(\sigma)\in H_{6,3}$ is non-zero.
    \item $Q^1(\sigma)Q^{2,1}(\sigma)\in H_{6,4}$ is non-zero.
\end{enumerate}
\end{prop}
\begin{prf}
    Item (i) is explained in \autocite[Lemma 6.7 and the following paragraph]{GKRW18b}. 

    For (ii) we  will show $\Sq^1_*$ of this class is non-zero. 
    \[\Sq^1_*(Q^1(\sigma)Q^{2,1}(\sigma)) =  \Sq^1_*(Q^1(\sigma))Q^{2,1}(\sigma) + Q^1(\sigma)\Sq^1_*(Q^{2,1}(\sigma)) = Q^1(\sigma)Q^{1,1}(\sigma) = Q^1(\sigma)^3.\]
    Applying $Q^2$ to the identity $\sigma Q^1(\sigma)=0$ of \cref{cor: g=3}(ii) gives us $Q^1(\sigma)^3 = \sigma^2Q^{2,1}(\sigma)$, which is non-zero by part (i).
\end{prf}

\subsection{Proof of  \cref{thm: Winfty structure} and \cref{cor: identities}}\label{sec: grouphom 2}

\begin{prf}[of \cref{thm: Winfty structure}]
 Suppose $X$ is a natural number graded $E_\infty$-space, then given any subset $B\subseteq H_{*,*}(X)$ there is an induced $W_\infty$-map $W_\infty[B]\to H_{*,*}(X)$ that is identity on $B$. For a subset $S\subseteq W_\infty[B]$, this map factors through $W_\infty[B]/(S)_{W_\infty}$ precisely when the image in $H_{*,*}(X)$ of every element in $S$ vanishes. So to prove (i) it suffices to check the classes in \[S\defeq \{\sigma Q^1(\sigma),\sigma Q^3(\sigma), \sigma^2Q^2(\sigma),\sigma\nu_1,\sigma\nu_2\}\]
 are zero in $H_{*,*}$.  For $\sigma Q^1(\sigma)$ this follows from \cref{cor: g=3}(iii); and the vanishing of the rest of $S$ is a consequence of \cref{prop: GL3 to GL4}.

We now prove (ii). Let $B = \{\sigma, \nu_1,\nu_2\}$, 
\begin{align*}
    QB&\defeq \{Q^I(b)\mid I\text{ normal for }b\}\cap W_{\infty}[\sigma,\nu_1,\nu_2]_{\leqslant 6,\leqslant 3 },\\
    QS&\defeq \{Q^J(s)\mid J\text{ normal for }s\}\cap W_{\infty}[\sigma,\nu_1,\nu_2]_{\leqslant 6,\leqslant 3},
\end{align*}
where a sequence is said to be normal if is admissible and satisfies certain excess condition; see \S\ref{sec: normal sequence} for details. Then the algebra $\F_2[QB]/(QS)$ agrees with $W_\infty[B]/I$ in bidegrees $(g,d)$ with $g\leqslant 6, d\leqslant 3$, and hence it's enough to show that $\overline{\varphi}:\F_2[QB]/(QS)\to W_\infty[B]/I\to H_{*,*}$ is a bijection in this range. The additive basis in \cref{thm: additive basis} shows that $\overline{\varphi}$ is surjective in the proposed range, so we only need to check that the dimensions of $\F_2[QB]/(QS)$ in each bidegree  agrees with that of $H_{*,*}$. Using \cref{lem: subsript normal seq}, one  obtains
\begin{align*}
    QB&= B\cup\{Q^1(\sigma), Q^2(\sigma),Q^3(\sigma),Q^{2,1}(\sigma)\}\\
    QS& = S\cup \{Q^2(\sigma Q^1(\sigma)) \} = S\cup \{Q^1(\sigma)^3+\sigma^2 Q^{2,1}(\sigma)\}.
\end{align*}
Then one write down an additive basis of $\F_2[QB]/(QS)_{\leqslant 6,\leqslant 3}$ row-by-row, manually.
\end{prf}

\begin{prf}[of \cref{cor: identities}]

(i)\; Apply $Q^3$ to the identity $\sigma Q^1(\sigma)=0$ from \cref{cor: g=3}(ii) gives us $\sigma^2 Q^{3,1}(\sigma)+Q^1(\sigma)Q^{2,1}(\sigma)+Q^2(\sigma)Q^1(\sigma)^2=0$. But $Q^{3,1}=0$ by the Adem relation.

(ii)\; Put $x = \sigma^2Q^4(\sigma)$ and $y = Q^2(\sigma)^2$.
Using the basis of $H_{4,4}$ provided by \cref{thm: additive basis} one has  \[Q^1(\sigma)Q^3(\sigma) = \lambda_1x +\lambda_2 y\] for some $\lambda_1,\lambda_2\in \F_2$. We will determine $\lambda_2$ by evaluating $\Sq^2_*$ on both sides. By the Cartan formula
\[\Sq^2_*(Q^1(\sigma)Q^3(\sigma)) =Q^1(\sigma)\Sq_*^2(Q^3(\sigma)) + \Sq^1_*(Q^1(\sigma))\Sq^1_*(Q^3(\sigma)) + \Sq^2_*(Q^1(\sigma))Q^3(\sigma).\]
The first, second, and respectively third term vanishes by the Nishida relation $\Sq^2_*Q^3 =Q^2\Sq^1_*$, $\Sq^1_*Q^3 = 0$, and respectively $\Sq^2_*Q^1 = Q^1\Sq_*^1$. On the other hand, we have computed that $\Sq^2_*(x)=0$ and $\Sq^2_*(y)\neq 0$ in the proof of \cref{prop: pi4*}(iii). So $\lambda_2=0$. We will show $\lambda_1=1$, i.e. that $Q^1(\sigma)Q^3(\sigma)$ is non-zero, by showing $Q^1(\sigma)Q^3(\sigma)$ pairs to $1$ with the class $w_4\in H^4(\GL_4)$ of \autocite[]{Pha87}. By definition \[w_4 =v_{14}+z_1(v_{12}^2+v_{13}+v_{34}^2 + v_{24}+v_{23}^2) + v_{12}v_{34}(v_{12}^2 + v_{13} + v_{34}^2 + v_{24}) + v_{12}^4 + v_{13}^2 + v_{34}^4 + v_{24}^2 +v_{23}^4\in H^*(\UT_4),\]
where the $v_{ij}$ and $z_1$  in $H^*(\UT_4)$ are characterised by their restriction to the elementary abelian subgroups $A_1,\dots, A_5$ described in (\ref{eq: Ai}). Using the table of these restrictions in \autocite[57]{Pha87}, one can check that the restriction of $w_4|_{A_4}$ along $\<a_{12},a_{34}>\leqslant A_4$ is
\[w_4|_{\<a_{12},a_{34}>}=\alpha_{12}^3\alpha_{34}+\alpha_{12}\alpha_{34}^3 + \alpha_{12}^4+\alpha_{34}^4,\]
which pairs to $1$ with the class $Q^1(\sigma)\times Q^3(\sigma)$ of $\<a_{12}, a_{34}>$, where $Q^i(\sigma)$ is the generator of $H_i(\Sigma_2)$. So $w_4$ pairs to $1$ with $Q^1(\sigma)Q^3(\sigma)$ as promised.

(iii) is obtained by stabilising (ii) once and using $\sigma Q^1(\sigma)=0$. 
\end{prf}

\section{Low degree group homology computations: $\GL_g(\Z)$}\label{sec: group homology Z}
Similar to the previous section, we describe a basis for $H_{g,d}\defeq H_d(\GL_g(\Z);\F_2)$ in low degrees in terms of Dyer-Lashof operations and products. Let $\sigma \in H_{1,0}$, $\theta_1\in H_{1,1}$, $\theta_2\in H_{1,2}$ be the generator of the respective homology group. Then our result can be summarised as follows. 

\begin{prop}\label{prop: dyer-lashof basis for BGL(Z)}
\cref{fig: basis for Z} describes a basis for $H_{g,d}$ in the range $g\leqslant3$, $d\leqslant2 $.
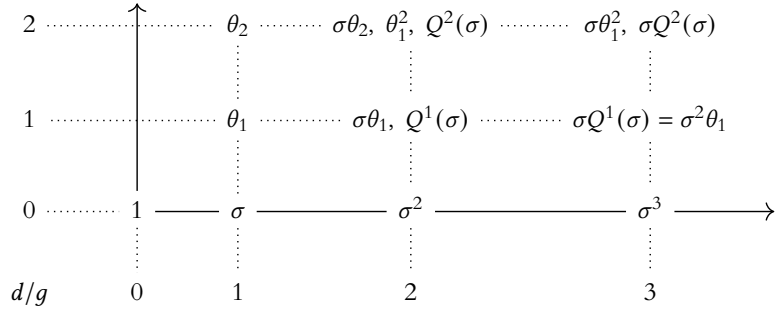
\begin{figure}[h]
    \centering
        \adjustbox{scale=0.95}{\begin{tikzcd}
	& { } \\
	2 && {\theta_2} & {\sigma\theta_2,\;\theta_1^2,\;Q^2(\sigma)} & {\sigma\theta_1^2,\;\sigma Q^2(\sigma)} \\
	1 && {\theta_1} & {\sigma\theta_1,\;Q^1(\sigma)} & {\sigma Q^1(\sigma)= \sigma^2\theta_1} \\
	0 & 1 & \sigma & {\sigma^2} & {\sigma^3} & {} \\
	d/g & 0 & 1 & 2 & 3
	\arrow[dotted, no head, from=2-1, to=2-3]
	\arrow[dotted, no head, from=2-3, to=2-4]
	\arrow[dotted, no head, from=2-3, to=3-3]
	\arrow[dotted, no head, from=2-4, to=2-5]
	\arrow[dotted, no head, from=2-4, to=3-4]
	\arrow[dotted, no head, from=2-5, to=3-5]
	\arrow[dotted, no head, from=3-3, to=3-1]
	\arrow[dotted, no head, from=3-3, to=3-4]
	\arrow[dotted, no head, from=3-3, to=4-3]
	\arrow[dotted, no head, from=3-4, to=3-5]
	\arrow[dotted, no head, from=3-4, to=4-4]
	\arrow[dotted, no head, from=3-5, to=4-5]
	\arrow[dotted, no head, from=4-1, to=4-2]
	\arrow[between={0}{0.8}, from=4-2, to=1-2]
	\arrow[no head, from=4-2, to=4-3]
	\arrow[dotted, no head, from=4-2, to=5-2]
	\arrow[no head, from=4-3, to=4-4]
	\arrow[dotted, no head, from=4-3, to=5-3]
	\arrow[no head, from=4-4, to=4-5]
	\arrow[dotted, no head, from=4-4, to=5-4]
	\arrow[between={0}{0.8}, from=4-5, to=4-6]
	\arrow[dotted, no head, from=4-5, to=5-5]
\end{tikzcd}}
    \caption{Basis for $H_d(\GL_g(\Z);\F_2)$.}
    \label{fig: basis for Z}
\end{figure}
\end{prop}
This proposition is a combination of  \cref{prop: dimension GL(Z)} and \cref{lem: computation BGLZ} below.

\begin{prop}\label{prop: dimension GL(Z)}
\cref{fig:dim for GLZ} describe $\dim H_{g,d}$ in the range $g\leqslant3$, $d\leqslant 2$.
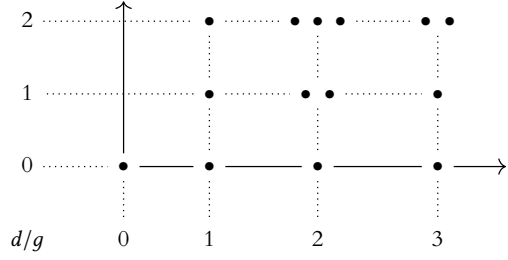
\begin{figure}[h]
    \centering
     \adjustbox{scale=0.85}{
\begin{tikzcd}
	& { } \\
	2 && \bullet & {\bullet\;\bullet\;\bullet} & {\bullet\;\;\bullet} \\
	1 && \bullet & {\bullet\;\;\bullet} & \bullet \\
	0 & \bullet & \bullet & \bullet & \bullet & {} \\
	d/g& 0 & 1 & 2 & 3
	\arrow[dotted, no head, from=2-1, to=2-3]
	\arrow[dotted, no head, from=2-3, to=2-4]
	\arrow[dotted, no head, from=2-3, to=3-3]
	\arrow[dotted, no head, from=2-4, to=2-5]
	\arrow[dotted, no head, from=2-4, to=3-4]
	\arrow[dotted, no head, from=2-5, to=3-5]
	\arrow[dotted, no head, from=3-3, to=3-1]
	\arrow[dotted, no head, from=3-3, to=3-4]
	\arrow[dotted, no head, from=3-3, to=4-3]
	\arrow[dotted, no head, from=3-4, to=3-5]
	\arrow[dotted, no head, from=3-4, to=4-4]
	\arrow[dotted, no head, from=3-5, to=4-5]
	\arrow[dotted, no head, from=4-1, to=4-2]
	\arrow[between={0}{0.8}, from=4-2, to=1-2]
	\arrow[no head, from=4-2, to=4-3]
	\arrow[dotted, no head, from=4-2, to=5-2]
	\arrow[no head, from=4-3, to=4-4]
	\arrow[dotted, no head, from=4-3, to=5-3]
	\arrow[no head, from=4-4, to=4-5]
	\arrow[dotted, no head, from=4-4, to=5-4]
	\arrow[between={0}{0.8}, from=4-5, to=4-6]
	\arrow[dotted, no head, from=4-5, to=5-5]
\end{tikzcd}}
    \caption{$\dim H_d(\GL_g(Z);\F_2)=$ number of dots at $(g,d)$. }
    \label{fig:dim for GLZ}
\end{figure}
\end{prop}
\begin{prf} The $g=1$ column is obtained by noting $\GL_1(\Z) = \Z^\times\simeq C_2$. 

For the remainder of the $d=1$ row, we use the identification $H_{g,1} = \GL_g(\Z)^{\mathrm{ab}}\otimes\F_2$. For $H_{2,1}$, consider the homomorphism $\mu:\GL_2(\Z)\to \Z/2$ defined by $\smqty(a& b\\c& d)\mapsto bc\;\mathrm{mod} 2$.
\begin{claimnumbered}\label{claim: det times mu}
    $\det\times\mu :\GL_2(\Z)^{\mathrm{ab}}\to \Z^{\times}\times \Z/2$ is an isomorphism.
\end{claimnumbered}
\begin{prfoc}
   We have a presentation
    \[\GL_2(\Z) = \<s,t,r\mid s^3 =t^2, t^4=1,r^2=1,rsr^{-1} =t^{-1}s^2t, rtr^{-1} =t^{-1}>\]
    \[s = \mqty(0 & 1 \\-1 & 1),\quad t= \mqty(0 & 1 \\-1 &0),\quad r = \mqty(-1 & 0 \\0 & 1)\]
    obtained from classical presentation of $\SL_2(\Z) = \<s,t\mid s^3 =t^2, t^4=1>$, which we learnt from \autocite[A]{KT08}, together with the semidirect product decomposition $\GL(\Z)\cong \SL_2(\Z)\rtimes C_2$, where the generator of $C_2$ acts by conjugation by $r$. Then 
    \[\GL_2(\Z)^{\mathrm{ab}} =\<\overline{t},\overline{r}\mid \overline{t}^2=1, \overline{r}^2=1, \overline{t}\overline{r}=\overline{r}\overline{t} >\]
    and one has $(\det\times \mu)(t) =(1, -1)$ and $(\det\times\mu)(r) =(-1,0)$. So the $\det\times \mu$ is surjective, but both the source and target are of the same order.
\end{prfoc}
When $g=3$, one has $\GL_3(\Z) = \<- I_3>\times \SL_3(\Z)$ as $I_3$ is central. As $\SL_3(\Z)$ is perfect (\autocite[4.3.9.Theorem and 1.2.15]{HO89}), the abelianisation of $\GL_3(\Z)$ is $C_2$ generated by the image of $-I_3$.

For $H_{2,2}$, suppose we have the following
\begin{claim}
    $H_2(\GL_2(\Z);\Z) = \Z/2$.
\end{claim}
Then  by the universal coefficient theorem, $H_2(\GL_2(\Z);\F_2)$ is (non-canonically) the direct sum of  $H_2(\GL_2(\Z);\Z)\otimes \F_2  = \F_2$ with $\Tor(\GL_1(\Z)^{\mathrm{ab}},\F_2) = \F_2\oplus \F_2$.

\begin{prfoc}
Consider the Hochschild-Serre spectral sequence associated to the extension
\[\SL_2(\Z)\to \GL_2(\Z)\xrightarrow[]{\det}\Z^{\times}.\]
The first few $E^2_{p,q} =H_p(\Z^\times ;H_q(\SL_2(\Z);\Z)) $ are as below.  The only entry that requires justification is  $E^2_{1,1} = H_1(\Z^{\times};H_1(\SL_2(\Z);\Z)) = \Z/2$. For this, we use the presentation of $\GL_2(\Z)$ and $\SL_2(\Z)$ above to get a description of $H_1(\SL_2(\Z);\Z)$ as a $\Z^{\times}$-module:  $H_1(\SL_2(\Z);\Z) =\<\overline{s}^3\mid \overline{s}^3 =1>\times \<\overline{t}\mid \overline{t}^4=1>$ with $\Z^{\times}$ action given by $(-1).\overline{s} = \overline{s}^{-1}$, $(-1).\overline{t} = \overline{t}^{-1}$.  \begin{center}
    \adjustbox{scale=0.85}{ \begin{sseqpage}[classes = { draw = none }, grid = crossword,homological Serre grading ]
     \class["\Z"](0,0) \class["\Z/2"](1,0) \class["0"](2,0)
     \class["\Z/2"](0,1) \class["\Z/2"](1,1)\class["\Z/2"](2,1)
     \class["0"](0,2)\class["\Z/2"](3,0)\d2(3,0)
 \end{sseqpage}}    \end{center}

The only possibly non-trivial differential leaving or entering an $E^2_{p,q}$ with $p+q=2$ is one drawn in the picture, namely $d^2:E^2_{3,0}\to E^2_{1,1}$. This differential is either zero or an isomorphism, and as a result $H_2(\GL_2(\Z);\Z)$ is either $\Z/2$ or $0$. We can exclude the latter case by considering the Mayer-Vietoris sequence
\[\cdots \to H_2(D_4;\Z)\to H_2(D_8;\Z)\oplus H_2(D_{12};\Z)\to H_2(\GL_2(\Z);\Z) \to \cdots \]
associated to the amalgamated free product decomposition\footnote{This decomposition is well-known, see for instance \autocite[Section 6]{BT16}. It can also be obtained from the presentation of $\GL_2(\Z)$ described above.}
\[\GL_2(\Z) \cong D_8\ast_{D_4} D_6.\]
Since $H_2(D_4;\Z) = \Z/2$, $H_2(D_8)\oplus H_2(D_{12}) =\Z/2\oplus \Z/2$, we cannot have $H_2(\GL_2(\Z);\Z)$ being zero.
\end{prfoc}

Finally, Soulé has completely computed  the integral cohomology of $\GL_3(\Z)$ \autocite[14]{Sou78}. We obtain the proposition for $H_{3,2}$ and  by using $H^2(\GL_3(\Z);\Z)=\Z/2$ and $H^3(\GL_3(\Z);\Z) = \Z/2$ together with the universal coefficient theorem.
 \end{prf}

\begin{lem}\label{lem: computation BGLZ}
    \begin{enumerate}[(i),before*=\leavevmode\vspace{-0.17\baselineskip}]
        \item $\sigma\theta_1,Q^1(\sigma)\in H_{2,1}$ are non-zero and distinct.
        \item $\sigma\theta_2, Q^2(\sigma)\in H_{2,2}$ are  non-zero and distinct.
        \item $\sigma\theta_1^2, \sigma Q^2(\sigma)\in H_{3,2}$ are non-zero and distinct. In particular $\theta_1^2\in H_{2,2}$ is non-zero.
        \item $\sigma Q^1(\sigma) =\sigma^2 \theta_1\in H_{3,1}$ and both are non-zero. 
        \item $\sigma\theta_2,\theta_1^2$ and $Q^2(\sigma)\in H_{2,2}$ are linearly independent.
    \end{enumerate}
\end{lem}

\begin{prf}
 Let $r_2:H_{g,*}\to H_*(\GL_g(\F_2);\F_2)$ be the reduction mod $2$ map, then $r_2(Q^1(\sigma)), r_2(Q^2(\sigma)), r_2(\sigma Q^2(\sigma))$  are all non-zero (\cref{prop: description of GL_2} and \cref{cor: g=3}(i)). Thus $Q^1(\sigma), Q^2(\sigma),\sigma Q^2(\sigma)\in H_{*,*}$ are non-zero. We now prove the remaining part of the lemma.

    (i)\;$\sigma\theta_1$ is represented by the matrix $\smqty(-1 & 0 \\ 0 & 1)$ and $Q^1(\sigma)$ by $\smqty(-1 & 0 \\0  &-1)$. So \[(\det\times \mu)(\sigma\theta_1) = (-1, 0)\quad\text{and}\quad  (\det\times\mu)(Q^1(\sigma)) = (1,0),\]where $\det\times \mu:H_{2,1}\to \Z^{\times}\times \Z/2$ is the homomorphism from \cref{claim: det times mu}. This shows that $\sigma\theta_1$ is non-zero and $\sigma\theta_1\neq Q^1(\sigma)$.

    (ii)\; Since $H_2(\GL_1(\F_2);\F_2)=0$, the class $r_2(\sigma\theta_2) = \sigma r_2(\theta_2)$ is zero, and hence $\sigma\theta_2\neq Q^2(\sigma)$.  For the non-vanishing of $\sigma \theta_2$, we will show that $\sigma:H_{1,2}\to H_{2,2}$ is non-zero. By the universal coefficient theorem it is enough to show that $\Tor(\GL_1(\Z)^{\mathrm{ab}},\Z/2)\to \Tor(\GL_2(\Z)^{\mathrm{ab}},\Z/2)$ is non-zero. By the proof of \cref{claim: det times mu}, the map on abelianisation sends the generator of $\GL_1(\Z)^{\mathrm{ab}}\simeq \Z/2$ to $\overline{r}$ in $\GL_2(\Z)^{\mathrm{ab}}$, which is $\Z/2\times \Z/2$ generated by $\overline{t}$ and $\overline{r}$.

    (iii)\; Again $r_2(\theta_1^2)=0$ shows that $\sigma \theta_1^2\neq \sigma Q^2(\sigma)$. To see $\sigma \theta_1^2$ is non-zero, consider the real $3$-dimensional representation $\rho:\GL_3(\Z)\to \GL_3(\R)$ which is the inclusion. Let $H\simeq C_2\times C_2$ be the subgroup of $\GL_2(\Z)$ generated by $a_1\defeq \smqty(-1 & 0 \\0 &1)$ and $a_2\defeq \smqty(1 &0 \\ 0 & -1)$, and $i:H\to \GL_2(\Z)\to \GL_3(\Z)$.   Then   
    \[i^*\rho\simeq (\mathrm{sign}\boxtimes \mathrm{triv})\oplus (\mathrm{triv}\boxtimes \mathrm{sign})\oplus (\mathrm{triv}\boxtimes \mathrm{triv}).\]
    In particular, its total Stiefel-Whitney class is 
    \[w(i^*\rho) = (1+\alpha_1)(1+\alpha_2) =1+\alpha_1+\alpha_2 +\alpha_1\alpha_2,\]
    where $\alpha_1,\alpha_2\in H^1(H;\F_2)$ is the dual basis of $a_1,a_2\in H_1(H;\F_2) = H\otimes \F_2$. We see that $w_2(i^*\rho)$ pairs to $1$ with $\theta_1\times \theta_1$, which means that $w_2(\rho)$ pairs to $1$ with $i_*(\theta_1\times\theta_1) = \sigma\theta_1^2$.

    (iv)\; We know $\dim H_{3,1}=1$, so we only need to check $\sigma^2\theta_1\neq 0$.  Indeed, $\sigma^2\theta_1$ is represented by $\smqty(-1 & 0  &0 \\0 & 1 & 0\\0  & 0 & 1)$, which has determinant $1$.

    (v)\; Given a general dependency equation  $ \lambda_1\sigma\theta_2+\lambda_2\theta_1^2 + \lambda_3Q^2(\sigma)=0$, $\lambda_i\in \F_2$,  we can evaluate $r_2$ on it to get  $\lambda_3=0$. Since $\sigma\times \theta_2\in H_*(H;\F_2)$ pairs to zero with $w_2(\rho)$, we have $0= \<w_2(\rho),\lambda_1\sigma\theta_2+\lambda_2\theta_1^2> =\lambda_2$. Since $\sigma\theta_2$ is non-zero, $\lambda_1=0$ as well.
\end{prf}
\section{Backgrounds}\label{sec: notations}
In  \S\ref{sec: background machineries} and \S\ref{sec: E_infty algs} we establish the general framework that we work in and set up some notations. This part  is included mainly for the sake of technical completeness. In  \S\ref{subsubsec: bar and cobar} we review the classical bar and cobar construction and defines the \textit{bar operation} $\overline{x}$ on homotopy classes $x$, which is essentially the effect of $\Ba$ on the attaching map of a cell. This part is crucial for \S\ref{sec: tools}, which concerns computation of stability Hopf algebra over $\F_2$. We will end this section by describing how $\overline{(-)}$ interacts with $\E_n$-operations, and establish a  Hurewicz theorem in this context which will be used later. Nevertheless, readers interested in computations should be able to proceed directly to \S\ref{sec: building the detector}-\ref{sec: applications of the detector}, and refer back to this section \S\ref{sec: tools} only when needed.

\subsection{Conventions}\label{sec: background machineries}

We work within the framework as set up in \autocite[]{HA}. In particular, we work with the Lurie model of $\E_n$-operads and their algebras. Our main $\E_\infty$-algebra $\CGL=\CGL(\F_2)$ (defined in \cref{const: algebra from symmetric monoidal groupoid}) lives in the (derived) category $\fancycat{D}(\F_2)^{\mathrm{gr}}\defeq \Fun(\Z,\fancycat{D}(\F_2))$ of graded $\F_2$-modules. 
For $\bbk$ a (discrete) commutative ring, $\fancycat{D}(\bbk)$ is the stable, presentably symmetric monoidal category obtained by 
inverting quasi-isomorphisms in the Hovey model category of chain complexes in $\bbk$-modules \autocite[\S1.3.5]{HA}. Alternatively,  one can view $\fancycat{D}(\bbk)$ as the modules over the $\E_\infty$-ring $H\bbk$ \autocite[\nopp7.1.1.16]{HA}. 
Here $\Z$ means the 1-category with objects indexed by the integers and no non-identity morphisms, and with the symmetric monoidal structure given by addition. Day convolution (c.f. \autocite[\S3]{BS23},  \autocite[]{Nik16}) equips $\fancycat{D}(\bbk)^{\mathrm{gr}}$ with a presentably symmetric monoidal structure. We will often equip functor categories between symmetric monoidal categories with  Day convolution without explicitly saying so.

\subsubsection{Bigraded homotopy groups.}\label{sec: bigraded homotopy groups} The \textit{$d$-sphere} in $\fancycat{D}(\bbk)$, $d$ an integer, is $\bbS^{d}\defeq \Sigma^d \mathbb{1}$,  where $\mathbb{1}$ is the monoidal unit.  The \textit{bigraded sphere} $\bbS^{g,d}$ in $\fancycat{D}(\bbk)^{\Z}$, $g,d\in \Z$, has  $\bbS^d$ in grading $g$ and the zero object in other gradings. The \textit{bigraded homotopy group} of an object $X$ in  $ \fancycat{D}(\bbk)^{\Z}$ is $\pi_{g,d}(X)\defeq \pi_d(X(g)) = \pi_0\Map(\bbS^{g,d},X)$. It computes the degree $d$ homology of the chain complex  $X(g)$.

\subsubsection{Filtered objects.} Let $\Z_{\leqslant}$ denote the poset of integers with its usual order, equipped with the symmetric monoidal structure given by addition. For us, a filtered objects will have ascending filtrations, i.e. they will be objects of the  stable, presentably symmetric monoidal category  $\Fil(\fancycat{D}(\bbk)^{\Z})\defeq \Fun(\Z_{\leqslant}, \fancycat{D}(\bbk)^{\Z}) $. The filtered homotopy group $\pi_{g,d,f}(\widetilde{X})$ is defined analogously \ref{sec: bigraded homotopy groups} using the filtered spheres $\bbS^{g,d,f}$, a filtration of $\bbS^{g,d}$ that is zero in filtration below $f$ and $\bbS^{g,d}$ in filtrations above and equal $f$. 

The recent {synthetic spectra revolution} exhibits the computational power of encoding spectral sequence yoga in the language of  filtered objects, which we shall make use of freely.  A practical guide to filtered objects adapted to our situation is \autocite[\S2.7]{RW25}. We recommend the original \autocite[Appendix A]{BHS22} for a more in-depth analysis.
We shall follow the set up in \autocite[\S2.7]{RW25}, but for convenience, let us recall the indexing convention of the Bockstein spectral sequence. Given a filtered object $\widetilde{X}$, one can extract a doubly filtered object $\beta\otimes \widetilde{X}$ that is the \say{$\tau$-adic filtration} on $\widetilde{X}$, with $\tau$ put in filtration $-1$. The associated spectral sequence has signature
\[\tensor*[^{\mathrm{BS}}]{E}{^1_{g,d,f,p}} =\tensor*[^{\mathrm{BS}}]{E}{^1_{g,d,f,p}}(\widetilde{X}) \cong \begin{cases*}
    \tau^{-p}\pi_{g,d,f+p}(C\tau\otimes  \widetilde{X} )&, $p\leqslant 0$\\
     0 &, $p>0$
\end{cases*}\implies \pi_{g,d,f}(\widetilde{X})\]
and differentials $\tensor*[^{\mathrm{BS}}]{d}{^r}:\tensor*[^{\mathrm{BS}}]{E}{^1_{g,d,f,p}}\to \tensor*[^{\mathrm{BS}}]{E}{^1_{g,d-1,f,p-r}}$. In this $\bbk[\tau]$-module description of the $E^1$-page,  $\tau$ has filtration $(0,0,1,-1)$. 
We denote the spectral sequence associate to $\widetilde{X}$ by $\{E^r_{*,*,*}(\widetilde{X}),d^r\}_r$.  The Omnibus theorem \autocite[Theorem 9.19]{BHS22} tells us that a $d^r$-differential $d^r(x) = y$ in $E^r(\widetilde{X})$ corresponds to a differential $\tensor*[^{\mathrm{BS}}]{d}{^r}(x)= \tau^r y$ in $\tensor*[^{\mathrm{BS}}]{E}{^r}(\widetilde{X})$.

\subsubsection{The diagonal $t$-structure \autocite[\S2.8]{RW25}.} 
We write $\tau^{\mathrm{diag}}_{\leqslant 0}:\fancycat{D}(\bbk)^{\Z}\to (\fancycat{D}(\bbk)^{\Z})^{\heartsuit}$ for the truncation functor with respect to the diagonal $t$-structure $t^{\mathrm{diag}}$ . We recall that $X$ is a connective (resp. coconnective) objects of $t^{\mathrm{diag}}$ iff each $\Sigma^gX(g)$ is connective (resp. coconnective) with respect to Postnikov $t$-structure on $\fancycat{D}(\bbk)$. Also   $t^{\mathrm{diag}}$ is compatible with the monoidal structure. When restricted to the connective objects, $\tau^{\mathrm{diag}}_{\leqslant 0}$ is strong symmetric monoidal, so its right adjoint $\iota$ is lax. When $\bbk$ is a field, $\iota$ is in fact strong symmetric monoidal \autocite[Warning 2.4]{RW25}. In particular $\iota$ sends $\E_1$-coalgebras in the heart to $\E_1$-coalgebras in  connective objects.

\subsection{$\E_\infty$-algebras and their properties}\label{sec: E_infty algs}
For  a commutative ring $R$,  there is an $\E_{\infty}$-algebra $\CGL(R)$ in natural number graded spaces whose homotopy type is the disjoint union of the classifying spaces $B\GL_g(R)$. The $\E_\infty$-structure is, roughly speaking, induced by block sum of matrices. We will be interested in its $\F_2$-linearised version $\CGL(R)_{\F_2}$, obtained by applying the $\F_2$-chain functor  $(-)_{\F_2}:\Alg_{\E_\infty}(\Spaces^{\Z})\to \Alg_{\E_\infty}(\fancycat{D}(\F_2)^{\Z})$. The honest construction is provided in \cref{const: algebra from symmetric monoidal groupoid}, but for  understanding the computations in \S\ref{sec: building the detector} and \S\ref{sec: applications of the detector}, the readers are invited to  take on  the following important properties $\CGL(R)_{\F_2}$  and skip the official definition.
    \begin{enumerate}[(i)]
     \item A map of commutative rings $R\to S$ induces an $\E_\infty$-map $\CGL(R)_{\F_2}\to\CGL(S)_{\F_2}$. 
    \item $\pi_{g,d}(\CGL(R)_{\F_2})$ computes the group homology $ H_d(\GL_g(R);\F_2)$. In particular, there is a canonical class $\sigma\in \pi_{1,0}(\CGL(R)_{\F_2})$, namely the generator of $H_0$. In particular, this class is compatible with maps in (i).
    \item The homotopy groups  $\pi_{g,d}(\CGL(R)_{\F_2}/\sigma)$ of the cofibre of multiplication by $\sigma$ computes the relative homology groups $H_d(\GL_g(R),\GL_{g-1}(R);\F_2)$.
\end{enumerate}

\begin{construction}\label{const: algebra from symmetric monoidal groupoid} Given a symmetric monoidal groupoid $\fancycat{G}$ together with a symmetric monoidal functor $r:\fancycat{G}\to \Z$, there is a symmetric monoidal $r_!:\Fun(\fancycat{G}, \fancycat{S})\to \Fun(\Z, \fancycat{S}) = \Spaces^{\mathrm{gr}}$  which is objectwise 
left Kan extension along $r$ \autocite[Proposition 3.6]{BS23}. The constant functor $\underline{\ast}$ to the terminal object $\ast\in \fancycat{S}$ has a canonical $\E_\infty$-algebra structure. Following  \autocite[\S16.1]{GKRW18a}, we set \[\mathrm{B}\fancycat{G}\defeq r_!(\underline{\ast})\in \Alg_{\E_\infty}(\Spaces^{\Z})\]
to be the \textit{$\E_\infty$-algebra associated to $(\fancycat{G},r)$}.   For every commutative ring $R$ with connected $\Spec (R)$,  taking rank defines a symmetric monoidal functor from  the symmetric monoidal 1-groupoid whose objects are finitely generated projective $R$-modules, to the integers. We define $\CGL(R)$ to be the associated $\E_\infty$-algebra in graded spaces. Finally, let $(-)_{\F_2}:\Spaces^{\Z} \to \fancycat{D}(\F_2)^{\Z}$ be the $\F_2$-chain functor\footnote{Let $\bbk$ be a commutative ring. Since $\fancycat{D}(\bbk)$ is presentably symmetric monoidal, the universal property of Day convolution says that restricting along Yoneda induces an equivalence $\Fun^{\otimes,L}(\Spaces, \fancycat{D}(\bbk))\simeq \Fun^{\otimes}(\ast, \fancycat{D}(\bbk))$. Let $C_*(-;\bbk)$ be the symmetric monoidal colimit preserving functor corresponding to the monoidal unit on right hand side. Note on object $(X_{\bbk})(n) = \colim_{X(n)}\bbk$.}. The functor $(-)_{\F_2}$ is symmetric  monoidal\footnote{by \autocite[Proposition 3.3]{BS23}.}, in particular it lifts to a functor between $\E_\infty$-algebras.
\textit{Since this document only concerns $\F_2$ homologies, from now on $\CGL(R)$ will mean $\CGL(R)_{\F_2}$.}
\end{construction}

To study homological stability of an $\E_\infty$-algebra using the methods of \autocite[]{RW25}, we need some mild assumptions on the algebra.
\begin{itemize}[$\ast$]
 \item \textbf{Connected algebras.} We say $A\in \Alg_{\E_\infty}(\fancycat{D}(\bbk)^{\Z})$ is \textit{connected} if  $\pi_{*,*}(A)$ supported in $\{(g,d)\mid g,d\geqslant 0\}$ and the unit $\mathbb{1}\to R$ is an isomorphism on $\pi_{0,0}$. In particular, a connected algebra is canonically augmented.
\item \textbf{(SCE).}  An augmented $\E_\infty$ algebra $\epsilon:A\to \mathbb{1}$ in $\fancycat{D}(\bbk)^{\Z}$  provides a $R,R$-bimodule structure on $\mathbb{1}$. We say $\epsilon:A\to \mathbb{1}$ satisfies the \textit{standard connectivity estimate (SCE)} if $\mathbb{1}\otimes_A\mathbb{1}$ is connective with respect to the $t^{\mathrm{diag}}$, i.e. $\pi_{g,d}(\mathbb{1}\otimes_R\mathbb{1})=0$ whenever $d<g$. 
\end{itemize}

We will compute the bar construction on an $\E_\infty$-algebra $A$ by finding a $\E_\infty$-cell structure on $A$. Given $x\in \pi_{g,d}(A)$, we write $A/\!/_{\E_\infty}x$ for the cofibre in $\Alg_{\E_\infty}(\fancycat{D}(\bbk)^{\Z})$ of $\E_\infty(\bbS^{g,d}\{x\})\to A$. Often we write $A/\!/^{g,d+1}_{\E_\infty}x$  to remind ourself the dimension of the cell attached.

\subsection{Bar and cobar.}\label{subsubsec: bar and cobar} For $0< n\leqslant\infty$, there is an adjunction \stepcounter{thm}
\[\leqnomode\tag{\thethm}\label{eq: bar-cobar}\Ba:\Alg^{\mathrm{aug}}_{\E_n}\left(\fancycat{D}(\bbk)^{\Z}\right)\begin{tikzcd}
     {} \arrow[r,yshift = 0.4ex] &\arrow[l,yshift =-0.4ex]{}
\end{tikzcd} \Coalg^{\mathrm{aug}}_{\E_1}\left(\Alg^{\mathrm{aug}}_{\E_{n-1}}\left(\fancycat{D}(\bbk)^{\Z}\right)\right):\Cobar
\]
obtained by applying the Bar-Cobar adjunction of \autocite[\nopp5.2.2.19]{HA} to the symmetric monoidal category of augmented $\E_{n-1}$-algebras, together with the Dunn-Lurie additivity \autocite[\nopp4.8.5.20]{HA}. Our  $n$ will mostly be $\infty$, except in \S\ref{sec: tools} where we study the effect of taking $\Ba$ to the $\E_n$-operations for general $n$. 

For convenience, we might refer the target of $\Ba$ as augmented $\E_{1,n-1}$-\textit{bialgebras}. We say an augmented $\E_n$-algebra $\epsilon:A\to \mathbb{1}$ is \textit{connected} if its augmentation ideal (i.e. the fibre of $\epsilon$ in graded $\bbk$-modules) has homotopy groups supported in $\{(g,d)\mid g,d\geqslant 0\}\setminus \{(0,0)\}$. We say an augmented bialgebra is \textit{simply-connected} if the augmentation ideal has homotopy group supported in $\{(g,d)\mid g,d\geqslant 0\}\setminus \{(0,0),(0,1)\}$. 
We will occasionally use \autocite[Proposition 2.2]{RW25}\footnote{The cited result is the $n=1$ case. But  as both connectedness and simply-connectedness are properties on the underlying object, the adjunction indeed restricts to the claimed subcategories. We are left to check the the unit and  the counit are equivalences. Again those can be checked on the underlying. }, recalled below; versions of the result are well-known.

\begin{prop}\label{prop: bar cobar restrict to equiv} The adjunction 
(\ref{eq: bar-cobar}) restricts to an equivalence between the full subcategories of connected augmented $\E_n$-algebras and simply-connected $\E_{1,n-1}$-bialgebras. \hfill $\blacksquare$
\end{prop}

\subsubsection{The bar operation}\label{sec: the bar operation} We would like to understand the effect of bar construction on an $\E_n$-cell structure. This can be summarised by \cref{def: bar operation} below: the bar construction turns an $\E_n$-cell attached along a homotopy class $x$ to an $\E_{n-1}$-cell attached along $\overline{x}$. Let $n\geqslant1$ and $A$ be an augmented $\E_n$-algebra, we write $\overline{A}$ for the augmentation ideal of $A$, i.e. the fibre of the augmentation map $\epsilon:A\to \mathbb{1}$. Then $A\cong \mathbb{1}\oplus \overline{A}$ and we can view the homotopy groups of $\overline{A}$ as a summand of that of $A$. The homotopy type of $\Ba A$ can be computed as a colimit of a simplicial diagram: There
is an equivalence $\colim_{\Delta^{\mathrm{op}}}A^{\otimes\bullet}\to \Ba A$ functorial in $A$. Here $A^{\otimes \bullet}$
is a simplicial object  with inner face maps built using multiplications and outer face maps given by collapsing along the augmentation map $\epsilon$. Consider the  the colimit of the truncated diagram:
 \stepcounter{thm}
\[\tag{\thethm}\label{eq: skeleton}\leqnomode\sk_i(A)\defeq \colim_{{\Delta}^{\mathrm{op}}_{\leqslant i}}A^{\otimes\bullet},\quad i\geqslant 0.\] Evidently  $\sk_0(A)=\mathbb{1}$. There is a natural identification \[\sk_1(A)\cong \mathbb{1}\oplus \Sigma\overline{A}\] as follows. We can map the truncated diagram  to the constant diagram with value $\mathbb{1}$ via the augmentation map. Since $\fancycat{D}(\bbk)^{\Z}$ is stable, on colimit this splits as the colimit of the constant diagram in $\mathbb{1}$, with that of the pointwise fibre $\begin{tikzcd}[cramped]
    0 \arrow[r] & \overline{A}\arrow[l, yshift =1ex ] \arrow[l, yshift =-1ex ]
\end{tikzcd}$, in which $\begin{tikzcd}
    0  & \overline{A}\arrow[l, yshift =0.7ex ] \arrow[l, yshift =-0.7ex ]
\end{tikzcd}$ is final.
In particular we have the following map \stepcounter{thm}
\[\label{not: inclusion of 1-skeleton}\tag{\thethm}i_{\overline{A}}:\Sigma\overline{A}\hookrightarrow \mathbb{1}\oplus\Sigma\overline{A}\to \Ba A\leqnomode.\]

\begin{defn}\label{def: bar operation}
    Let $A$ be an augmented $\E_n$-algebra, and $x\in \pi_{g,d}(\overline{A})$, the \textit{bar operation on $x$} is 
    \[\overline{x}\defeq (i_{\overline{A}})_*(\Sigma x)\in \pi_{g,d+1}(\Ba A).\]
\end{defn}

The following alternative description of the bar operation will be used later in \S\ref{sec: tools}.

\begin{lem}
Let $A$ and $x$ be as in \cref{def: bar operation}. Then  $\overline{x}$ is the homotopy class\footnote{Here we are using the identification $\pi_{g,d+1}(-) = \pi_0\Map_{\Alg_{\E_{n-1}}}(\E_{n-1}(\bbS^{g,d+1}),-)$ provided by the free-forgetful adjunction of $\E_{n-1}$-algebras. } in $\E_{n-1}$-algebras of $\Ba$ of any $\E_n$ representative $\E_n(\bbS^{g,d})\to A$ of $x$ precomposed with the $\E_{n-1}$-equivalence 
    \stepcounter{thm}
\[\leqnomode\tag{\thethm}\Ba(\E_n(-))\cong \E_{n-1}(\Sigma-)\label{eq: bar free suspension}\] via \autocite[\nopp  5.2.2.13]{HA}.  For instance, if  $A = \E_n(\bbS^{g,d}\{x\})$, then $\Ba (A) \cong \E_{n-1}(\bbS^{g,d+1}\{\overline{x}\})$.
\end{lem}
\begin{prf}
    Since both constructions are natural in $A$, we may assume $A = \E_n(\bbS^{g,d}x)$. Since both operation are compatible with forgetting part of the $\E_n$-structure, we may further assume $n=1$. In this case, the equivalence (\ref{eq: bar free suspension})  we use  to form $\overline{x}$ is the composite\stepcounter{thm}
    \[\label{eq: suspension to bar on free alg}\tag{\thethm}\mathbb{1}\sqcup_{\bbS^{g,d}\{x\}}\mathbb{1}\to \mathbb{1}\sqcup_{A}\mathbb{1}\to \Ba A,\leqnomode\]
    where the first map is induced by the unit $\bbS^{g,d}\{x\}\to A$ of the free-forgetful adjunction and the second  is $i_A$  (\autocite[\nopp5.2.2.14]{HA}). We can take the $\E_1$ representative of $x\in \pi_{*,*}(A)$ to be the the identity map, so both $i_{\overline{A}*}(x)$ and $\overline{x}$ are represented by the (\ref{eq: suspension to bar on free alg}) above.
\end{prf}

A key property  of the bar operation is that it kills decomposables:

\begin{lem}\label{lem: bar kills decomposables}
    Let $A$ and $x$ be as above. If $x$ factors through the multiplication map $\overline{A}\otimes \overline{A}\to \overline{A}$, then $\overline{x}=0$. 
\end{lem}

\begin{prf}
Let $\overline{\sk_i(A)}$ denote the summand in the natural splitting $\sk_i(A) = \mathbb{1}\oplus \overline{\sk_i(A)}$, and
let $F$ denote the fibre of $\overline{\sk_1(A)}\to \overline{\sk_2(A)}$.    We will see in \cref{lem: admitting qd implies null} later that the map $F\to \overline{\sk_1}(A)$ is equivalent to the multiplication $\overline{A}^{\otimes 2}\to \overline{A}$. Since $i_{\overline{A}}$ factors through $\overline{\sk_1(A)}\to \overline{\sk_2(A)}$, the result follows.
\end{prf}

We now describe how to \say{commute} the bar operation pass the $\E_n$-operations. \label{sec: En operations}The homotopy groups of an $\E_n$-algebra in $\fancycat{D}(\F_2)^{\Z}$, $n\geqslant 1$,  has a sequence of natural operations (the Dyer-Lashof operations) \[Q_s^{\E_n}:\pi_{*,*}(-)\Rightarrow\pi_{2*,2*+s} (-),\quad 0\leqslant s<n\] starting with $Q_0^{\E_n}(x) = x^2$. 
When $n<\infty$, there is also the Browder bracket
\[\lambda_{\E_n} \defeq  [-,-]_{\E_k}: \pi_{*,*}(-)\otimes\pi_{*',*'}(-)\Rightarrow \pi_{*+*', *+*'+n-1}(-)\]

\begin{thm}\label{thm: bar and operations}
    Suppose $A\in \Alg^{\mathrm{aug}}_{\E_n}(\fancycat{D}(\F_2)^{\Z})$, $n\geqslant 2$. Then\footnote{Over $\F_2$ the lower indexed Dyer-Lashof operations are related to the upper-indexed ones by the formula $Q_s(x) = Q^{\abs{x}+s}(x)$.}
    \begin{enumerate}[(i)]
  \item For every $x\in \pi_{*,*}(A)$, $\overline{Q^{\E_n}_s(x)} = Q_{s-1}^{\E_{n-1}}(\overline{x})\in\pi_{*,*}(\Ba A)$ for every $0\leqslant s<n$. 
     This includes the top operation $Q_{n-1}^{\E_n}$, often denoted by $\xi$. We set $Q_{-1}^{\E_{n-1}}$ to be zero. 
     \item For pair of classes $a,b\in \pi_{*,*}(A)$,  $\overline{[a,b]_{\E_n}} = [\overline{a},\overline{b}]_{\E_{n-1}}\in \pi_{*,*}(\Ba A)$. Note that the $\E_1$-bracket is the commutator: $[x,y]_{\E_1} = xy-yx$ \autocite[p.125, Theorem 1.2 (2)]{CLM76}. 
\end{enumerate}
\begin{prf}
By naturality it suffices to consider the case
\begin{enumerate}[(i)]
    \item  $A = \E_n(\bbS^{g,d}x)$ with $g,d$ arbitrary integers;
    \item $A= \E_n(\bbS^{g,d}a\oplus \bbS^{g',d'}b)$ with $g,g',d,d'$ arbitrary integers. 
\end{enumerate}

    (i) The $s=0$ part is a special case of \cref{lem: bar kills decomposables}, so we may assume $s>0$. Suppose first $g\neq 0$. Since $\Ba A\simeq \E_{n-1}(\bbS^{g,d+1}\overline{a})$,  the class $Q_{s-1}^{\E_{n-1}}(\overline{a})$ is non-zero and span $\pi_{2g,2d+s+1}(\Ba A)$. So we just need to show $\overline{y}$  with $y\defeq 
        Q_s^{\E_{n}}(a)$ is also non-zero, and, since $y\in \pi_{g,2d+s}(A)$ is  non-zero, it's enough to show $\pi_{g,2d+s+1}(i_A)$ is injective. Consider the cofibre sequences \[\sk_i(A)\to \sk_{i+1}(A)\to \Sigma^{i+1}\overline{A}^{\otimes i+1}\] for           $i\geqslant 0$. When  $g>0$ (resp. $g<0$), $\Sigma^{i+1}\overline{A}^{\otimes i+1}$ is supported in grading $\geqslant (i+1)g$ (resp. $\leqslant (i+1)g$). In particular $\sk_i(A)\to \sk_{i+1}(A)$ is an equivalence in grading $2g$ when $i\geqslant 2$. On the other hand, $\pi_{2g,*}(\Sigma^2\overline{A})$ is 1-dimensional concentrated in degree $2d+2$, and the boundary map \[\pi_{2g,2d+2}(\Sigma^2\overline{A})\to \pi_{2g,2d+1}(\sk_1(A)) = \F_2\] is surjective by the long exact sequence and the fact that $\pi_{2g,2d+1}(\sk_2(A)) = \pi_{2g,2d+1}(\Ba A)=0$. So in fact $\pi_{2g,D}(\sk_1(A)\to \sk_2(A))$ is an isomorphism for $D\geqslant 2d+2$. We always have $2d+s+1\geqslant 2d+2$.

        The $g=0$ case follows from the $g\neq 0$ case by the change of groupoid trick: Fix any $g'\neq 0$ and take $A' = \E_n(\bbS^{g', d}x)$. Then $A$ is equivalent to left Kan extend $A'$ along the composite of symmetric monoidal functors $f:\Z\to \{0\}\hookrightarrow \Z$. Since both the domain and the target groupoid have no non-identity morphisms, the homotopy groups of the left Kan extension admits a natural decomposition 
        \[\pi_{h,d}(f_!-) = \bigoplus_{\substack{g\in \Z\\f(g)= h}}\pi_{g,d}(-)\]
        which is respected by the Dyer-Lashof operations and the bar operations.

        (ii) can be deduced from (i) and the fact that the top Dyer-Lashof is linear up to the Browder bracket \autocite[16.2.2 (d)]{GKRW18a}. I.e. for $A = \E_n(\bbS^{g,d}a\oplus \bbS^{g',d'}b)$, we have
\begin{align*}
Q^{\E_{n-1}}_{n-1}(\overline{a}) + Q^{\E_{n-1}}_{n-1}(\overline{b})+ [\overline{a},\overline{b}]_{\E_{n-1}}&=
    Q^{\E_{n-1}}_{n-1}(\overline{a}+\overline{b})\\
    &=\overline{Q^{\E_n}_n(a+b)}\\
    &=\overline{Q^{\E_n}_n(a)+Q^{\E_n}_n(b)+[a,b]_{\E_{n}}}\\
    &=Q^{\E_{n-1}}_{n-1}(\overline{a}) +  Q^{\E_{n-1}}_{n-1}(\overline{b})+\overline{[a,b]_{\E_{n}}}.
\end{align*}
\end{prf}
\end{thm}

\subsubsection{A Hurewicz theorem} In this part we set up the Hurewicz theorem\footnote{The reason for not just citing \autocite[Corollary 11.12]{GKRW18a} is our concern of comparing models of bar construction and $\E_1$-indecomposables. The cited Hurewicz map is the map from the algebra to its derived $\E_1$-indecomposables, and we would have to use \autocite[Theorem 13.7]{GKRW18a} to relate the indecomposables to the once de-suspended bar construction.  It is not so clear that the two bar constructions agree, even though we believe they do.} in our framework.
 Let $R\to S$ be  a map of connected $\E_n$-algebras in $\fancycat{D}(\F_2)^{\mathrm{gr}}$, $n\geqslant 1$.  Let 
 \[i_{S,R}:\Sigma S/R\to \Ba(S)/\Ba(R)\]
 be the map on cofibres induced by $i_{\overline{R}},i_{\overline{S}}$ of  (\ref{not: inclusion of 1-skeleton}).
The following lemma is  the counterpart of \autocite[Corollary 11.12]{GKRW18a}.
\begin{thm}\label{lem: hurewicz}In the situation above,
    for integers $g,d$, if  $\pi_{g',d'}(S,R)=0$  whenever $g'\leqslant g,d'\leqslant d$, and  $(g',d')\neq (g,d)$, then
    \[\pi_{g,i+1}(i_{S,R}):\pi_{g,i+1}(\Sigma S,\Sigma R)\cong \pi_{g,i}(R,S) \to \pi_{g,i+1}(\mathrm{Bar}(S),\mathrm{Bar}(R))\]
    is an isomorphism for $i=d$, and an epimorphism for $i=d+1$. 
\end{thm}
\begin{prf} 
 Let $C_i$ denote the cofibre $\sk_i(S)/\sk_i(R)$, then there is a commuting diagram           
\[\begin{tikzcd}
	{\sk_{i-1}(R)} & {\sk_i(R)} & {\Sigma^i\overline{R}^{\otimes i}} \\
	{\sk_{i-1}(S)} & {\sk_i(S)} & {\Sigma^i\overline{S}^{\otimes i}} \\
	{C_{i-1}} & {C_i} & {\Sigma^i(\overline{S}/\overline{R})^{\otimes i}}
	\arrow[from=1-1, to=1-2]
	\arrow[from=1-1, to=2-1]
	\arrow[from=1-2, to=1-3]
	\arrow[from=1-2, to=2-2]
	\arrow[from=1-3, to=2-3]
	\arrow[from=2-1, to=2-2]
	\arrow[from=2-1, to=3-1]
	\arrow[from=2-2, to=2-3]
	\arrow[from=2-2, to=3-2]
	\arrow[from=2-3, to=3-3]
	\arrow[from=3-1, to=3-2]
	\arrow[from=3-2, to=3-3]
\end{tikzcd}\]
where all rows and columns are cofibre sequences.  Since $C_1 \simeq \mathbb{1}\oplus \Sigma(\overline{S}/\overline{R}) $ and $\colim_i C_i = \Ba(S)/\Ba(R)$, it's enough to check 
\begin{itemize}[$\ast$]
    \item $\pi_{g,*}(C_{i-1}\to C_{i})$ is an isomorphism for $i>2$, $*\leqslant d+1$;  
    \item $\pi_{g,*}(C_1\to C_2)$ is an isomorphism when $*<d+2$ and an epimorphism when $*=d+2$.
\end{itemize}
Since $\pi_{0,*}(\overline{R})=\pi_{0,*}(\overline{S})=0$, the rightmost cofibre sequence evaluated at $\pi_{*,0}$ all vanishes. By assumption $\pi_{<g, \leqslant d}(\overline{S}/\overline{R})=0$, so by Künneth $\pi_{<g+i-1, \leqslant d+i}(\Sigma^i(\overline{S}/\overline{R})^{\otimes i})=0$. 
\end{prf}

\begin{cor}\label{cor: connectivity passing to bar}
    In the situation of \cref{lem: hurewicz}:
    \begin{enumerate}[(i)]
       \item If $\pi_{g',d'}(S,R)=0$ whenever $g'\leqslant g, d'\leqslant d$, then $\pi_{j,i}(\mathrm{Bar}(S),\mathrm{Bar}(R))=0$ whenever $j\leqslant g$, $i\leqslant d+1$.
       \item If $R\to S$ is $(g,d)$-connected, then $\mathrm{Bar}(R\to S)$ is $(g,d+1)$-connected.\hfill $\blacksquare$
    \end{enumerate}
\end{cor}

\section{Tools for computing the stability Hopf algebra}\label{sec: tools}

Let $\fancycat{G}$ be a symmetric monoidal groupoid and  $r:\fancycat{G}\to \Z$ be a symmetric monoidal functor, \cref{const: algebra from symmetric monoidal groupoid} associates to this data  an $\E_\infty$-algebra $A\defeq \mathrm{B}\fancycat{G}_{\F_2}$ in graded $\F_2$-modules.
The point of view provided by \autocite{RW25} is that, if $A$  is connected and satisfies the standard connectivity estimate, then \stepcounter{thm}
\[\label{pov}\leqnomode \tag{\thethm}\begin{tabular}{c}
     homological stability of  $(\fancycat{G},r)$  \\
     with coefficients in  $\F_2$
\end{tabular}
\approx  \begin{tabular}{c}
     vanishing lines of iterated cofibres \\
     of the form $A/\alpha_1,\cdots ,\alpha_r$ 
\end{tabular}
\]
where each $\alpha_i$ is a self-map of $A/\alpha_1,\dots,\alpha_{i-1}$. Under this viewpoint, stability information of $R$ of slope $<1$  can be detected on a discrete Hopf algebra $\Delta_A$.

\begin{defn}[\autocite{RW25}]\label{def: stability Hopf alg}
    Let $n>0$ and $A$ a connected $\E_n$-algebra in $\fancycat{D}(\bbk)^{\Z}$ satisfying (SCE). Then the \textit{stability Hopf algebra} of $A$ is the coconnective truncation (with respect to the diagonal $t$-structure) of $\Ba A$:
    \[\Delta_A\defeq\tau^{\mathrm{diag}}_{\leqslant0}\Ba A\in \Coalg^{\mathrm{aug}}_{\E_1}(\Alg^{\mathrm{aug}}_{\E_{n-1}}((\fancycat{D}(\F_2)^{\Z})^{\heartsuit}). \]
     When $n\geqslant 2$, we identify $\Delta_A $ with a positively graded, connected Hopf algebra structure on $\bigoplus_{n\geqslant0} \pi_{n,n}(\Ba A)$. When $n>2$ this is  commutative.   When $n=1$, $\Delta_A$ is only a coalgebra.  
\end{defn}
Applying $\Cobar$ to the truncation map $\Ba A\to \Delta_A$ gives us an $\E_n$-map \stepcounter{thm}
\[\label{map: truncation map}\tag{\thethm}\leqnomode t_A:A\cong \Cobar( \Ba A)\to a\defeq \Cobar(\Delta_A).\]
 \autocite[Theorem 7.1,  Theorem 8.2, and Proposition 8.5]{RW25} tells us how to transfer stability information between $A$ and $\Delta_A$ using the map $t_A$.  Specialisations of these results that will be used later are collected in the following theorem.

\begin{thm}[\autocite{RW25}]\label{basechange thm}
 We say that a map $A\to A'$ of $\E_\infty$-algebras in $\fancycat{D}(\F_2)$ and a constant $\theta>0$ satisfies property \stepcounter{thm} (\thethm) if 
   for every $A$-module $M$,  $0<\lambda\leqslant \theta$, and any constant $\mu$, 
        \[\pi_{g,d}(M) =0\quad \text{for }d<\lambda g+\mu\iff\pi_{g,d}(M\otimes_AA')\quad\text{for }d<\lambda g+\mu.\] 
    \begin{enumerate}[(i),before*=\leavevmode\vspace{-0.17\baselineskip}]
        \item Let $A$ be a connected $\E_\infty$-algebra in $\fancycat{D}(\F_2)$ satisfying (SCE). Then $t_A:A\to a$ satisfies (\thethm) with $\theta=1$.
     \item Let $A\to B$ be  a map of connected, graded Hopf algebras over $\F_2$.
     \begin{enumerate}
         \item If  $A\to B$  is surjective with kernel supported in gradings $\geqslant N$, then $\Cobar(A)\to \Cobar(B)$ satisfies (\thethm) with $\theta = \frac{N-1}{N}$.
         \item  If  $A\to B$  is injective with cokernel supported in gradings $\geqslant N$, then $\Cobar(A)\to \Cobar(B)$ satisfies (\thethm) with $\theta = \frac{N-2}{N}$.
     \end{enumerate}
    \end{enumerate}
\end{thm}

The result of this section is the following  partial manual (\cref{thm: summary theorem of diag hopf}) for computing $\Delta_A$ out of an $\E_n$-cell structure of $A$. That is, we describe how attaching a cell changes the stability Hopf algebra, in terms of the bar operation (\cref{def: bar operation}) on the attaching map. An informal explanation of the $[-]$-construction that appears in the statement will be given immediately after; the precise definition take some work and will occupy \S\ref{sec: Toda bracket}.

\begin{thm}\label{thm: summary theorem of diag hopf} Let $n\geqslant2$, $A$ be a connected, augmented $\E_n$-algebra in  $\fancycat{D}(\F_2)^{\Z}$ satisfying (SCE). Let $x\in \pi_{g,d}(A)$ be a homotopy class, and $A'\defeq A/\!/^{g,d+1}_{\E_n}x$ be formed by attaching an $\E_n$-cell along $x$. Suppose $d\geqslant g-2>0$, so that $A'\defeq A/\!/^{g,d+1}_{\E_n}x$ also satisfies (SCE), then:
\begin{enumerate}[(i)]
    \item When $d>g$, the natural map $A\to A'$ induces an isomorphism $\Delta_A\to \Delta_{A'}$ of Hopf algebra.
    \item When $d+1 = g$, $\Delta_A\to \Delta_{A'}$ realises $\Delta_{A'}$ as the quotient algebra $\Delta_{A}/(\overline{x})$.
    \item Suppose $d+2=g$ and  $\overline{x}=0$. Let $\mu_A:A\otimes A$ be the multiplication on $A$. Given an expression of the form \[x = (\mu_A)_*(q(x)),\quad  q(x)=\sum_i a_i\otimes b_i\]  with  the bidegrees  $\abs{a_i}=(g_i,d_i)$ and $\abs{b_i} = (g_i',d_i')$ satisfying $d_i +1=g_i>0$ and $d_i'+1 = g_i'>0$, there is an associated 
        homotopy class of maps\[[q(x)]:\Sigma \Sigma \bbS^{g,d}\to \Ba(A')\] inducing an equivalence of $\E_{n-1}$-algebras
        \[\Ba(A)\sqcup^{\E_{n-1}}\E_{n-1}(\bbS^{g,d+2}\{[q(x)]\})\xrightarrow[]{\cong} \Ba (A').\]
        In particular, if we view $[q(x)]$ as an element of $\Delta_{A'}$ in grading $g$, then 
        \begin{align*}
          \Delta_{A}\ast \F_2\<[q(x)]> &\to \Delta_{A'} &\text{if\;}  n=2\\
          \Delta_A\otimes \F_2[[q(x)]]&\to \Delta_{A'} &\text{if\;} n>2
        \end{align*}
        are isomorphisms of Hopf algebras. Here $\ast$ means free product of associative  $\F_2$ algebras, and $\otimes$ means tensor product of commutative $\F_2$ algebras. In both cases, the coproduct of  $[q(x)]$ is 
        \[{\psi}([q(x)]) = 1\otimes [q(x)]+[q(x)]\otimes 1+  \sum_i \overline{a_i}\otimes \overline{b_i}.\]
        If $q(x)$ is of the form $a\otimes b$, then we write $[a\mid b]$ for  $[q(x)]$.
\end{enumerate}
\end{thm}
 Roughly speaking, a decomposition $q(x)$ provides a reason for $\overline{x}$ to be zero in $\Ba(A)$, and hence in $\Ba(A')$. On the other hand, the cell we attached is a reason for $x$ to be zero in $A'$, and thus another reason for $\overline{x}$ to be zero in $\Ba(A')$. By gluing these two null-homotopies one obtain $[q(x)]$.

 \subsection{The $[-]$-construction}\label{sec: Toda bracket}
Let $A$ be a connected (and hence canonically augmented) $\E_n$ algebra $A$ in $\fancycat{D}(\F_2)^{\Z}$ satisfying (SCE). We write $\overline{A}$ for the augmentation ideal of $A$ and identify $\pi_{*,*}(\overline{A})$ as the summand of $\pi_{*,*}(A)$ on those classes that vanishes under the augmentation map. 
Let's fix a point in the space of binary operations for $\E_1\hookrightarrow\E_n$, and let $\mu_A:A\otimes A\to A$ be the multiplication determined by this point \autocite[\nopp 5.3.3.1]{HA}.  Write 
 $\mu_{\overline{A}}:\overline{A}\otimes \overline{A}\to \overline{A}$ be the restriction of $\mu_A$ using the decomposition of $A\otimes A$ induced by  $A\cong \mathbb{1}\oplus \overline{ A}$.  Since $\mu_A$ is functorial in $A$, so is $\mu_{\overline{A}}$. 
\begin{defn}\label{def: quadratic decomp}
    A \textit{decomposition} $q(x)$ of a class $x\in \pi_{g,d}(\overline{A})$ is a factorisation  in $\fancycat{E}$  of a representative of $x$  through $\mu_{\overline{A}}$. I.e. it is the data of a diagram
    \[\begin{tikzcd}
        \bbS^{g,d}\arrow[r]\arrow[rr, bend right =30,"x"] &\overline{A}\otimes \overline{A} \arrow[r, "\mu_{\overline{A}}"]  & \overline{A}
    \end{tikzcd}\] 
    in $\fancycat{D}(\F_2)^{\Z}$. We will also denote the map $\bbS^{g,d}\to \overline{A}\otimes \overline{A}$ by $q(x)$. 
\end{defn}
Since the $\Ba(-)$ is \say{same as $\E_1$-indecomposables up to a shift}, it should kill all the decomposables. Indeed, the following lemma gives us a null-homotopy of 
\[\Sigma \overline{A}^{\otimes 2}\xrightarrow[]{\Sigma\mu_{\overline{A}}}\Sigma\overline{A}\to \Ba (A),\]
functorial in $A$. As in \S\ref{sec: the bar operation}, we write $\sk_i(A) = \colim_{\Delta^{\mathrm{op}}_{\leqslant i}}A^{\otimes i}$ and $\overline{\sk_i(A)} =\mathrm{cofib}(\mathbb{1}\to \sk_i(A))$

\begin{lem}\label{lem: admitting qd implies null}
    Let $F$ denote the fibre of $\overline{\sk_1(A)}\to \overline{\sk_2(A)}$. Then the map $F\to \overline{\sk_1(A)}$ is equivalent to the suspended multiplication  $\Sigma \mu_{\overline{A}}:\Sigma \overline{A}^{\otimes 2}\to \Sigma\overline{A}$. 
\end{lem}

\begin{defn}
    Let $\varepsilon_{\mathrm{can.}A}$ denote the null-homotopy of $i_{\overline{A}}\Sigma\mu_{\overline{A}}:\Sigma\overline{A}\otimes\overline{A}\to \Ba (A)$ obtained by composing the null-homotopy of $\Sigma \overline{A}^{\otimes 2}\to \overline{sk_2}(A)$ coming from the fibre sequence with the map $\overline{\sk_2}(A)\to \Ba (A)$. Given a decomposition $q(x)$, \stepcounter{thm}
    \[\label{eq: null homotopy of q(x)}\leqnomode\tag{\thethm}\varepsilon_{q(x)}\defeq \Sigma q(x).\varepsilon_{\mathrm{can.}A}\]
    is a null-homotopy of $\overline{x}:\Sigma \bbS^{g,d}\to\Ba A$.
\end{defn}

\begin{prf}[of \cref{lem: admitting qd implies null}]
    Let $\iota:\Delta^{\mathrm{op}}_{\leqslant 1}\hookrightarrow \Delta^{\mathrm{op}}_{\leqslant 2}$ denote the inclusion, $D_i:\Delta^{\mathrm{op}}_{\leqslant i}\to \fancycat{D}(\F_2)^{\Z}$ be the diagram obtained by splitting off the unit in the truncated simplicial diagram $A^{\otimes \bullet}$. Then $\overline{\sk_1(A)}\cong \colim D_1 \cong \colim \iota_!D_1$, and the map of diagrams $\iota_!D_1\to D_2$ is described as below.
    \[\begin{tikzcd}[ sep = large]
    0 \arrow[r] \ar[equal]{d}& \overline{R} \arrow[l, yshift =1ex]\arrow[l,yshift = -1ex]\arrow[r,yshift = -1ex]\arrow[r,yshift = 1ex]  \ar[equal]{d}& (\mathbb{1}\otimes\overline{R})\oplus (\overline{R}\otimes \mathbb{1}) \arrow[l]\arrow[l, yshift= 2 ex] \arrow[l,  yshift =-2ex]\arrow[d]\\
    0 \arrow[r]  & \overline{R} \arrow[l, yshift =1ex]\arrow[l,yshift = -1ex]\arrow[r,yshift = -2ex]\arrow[r,yshift = 2ex] & (\mathbb{1}\otimes \overline{R}) \oplus( \overline{R}\otimes\mathbb{1}) \oplus (\overline{R}\otimes \overline{R} )\arrow[l,"(\ast)"{description}]\arrow[l, yshift= 3ex] \arrow[l,  yshift =-3ex]
\end{tikzcd}\]
In the top row the three face maps $(\mathbb{1}\otimes\overline{R})\oplus (\overline{R}\otimes \mathbb{1})\to \overline{R}$ are projection in to the first factor, the zero map, and projection to the second factor; and the degeneracy maps $\overline{R}\to (\mathbb{1}\otimes \overline{R})\oplus (\overline{R}\otimes 1)$ are inclusion into first and second factors. In the bottom row, the unlabelled face and degeneracy maps in the second column are again given by inclusions and projections. The map $(\ast)$ is projection into the $\overline{R}\otimes \overline{R}$ summand followed by multiplication. The rightmost vertical map is inclusion of summands. Since $\fancycat{D}(\F_2)^{\Z}$ is stable, $F$ can be computed as the colimit of pointwise fibre $D_1\to D_2$.
\end{prf}
We are now ready to describe the $[-]$-construction. 
\begin{construction}\label{defn: toda [q(x)]} 
 Since $\Ba:\Alg_{\E_n}^{\mathrm{aug}}(\fancycat{D}(\F_2)^{\Z})\to\Alg_{\E_{n-1}}^{\mathrm{aug}}(\fancycat{D}(\F_2)^{\Z}) $ is a left adjoint, it sends the cell-attachment pushout square defining $A'$ to a pushout square \stepcounter{thm}
\[\leqnomode \tag{\thethm}\label{sqaure: bar pushout}\begin{tikzcd}
    \E_{n-1}(\Sigma \bbS^{g,d}) \arrow[r]\arrow[d] & \Ba A\arrow[d]\\
    \E_{n-1}(0) \arrow[r] & \Ba(A')
\end{tikzcd}\]
in (augmented) $\E_{n-1}$-algebras. In particular it adjoints to a square
\stepcounter{thm}
\[\leqnomode \tag{\thethm}\label{sqaure: toda ingredient 2}\begin{tikzcd}
    \Sigma \bbS^{g,d} \arrow[r, "\overline{x}"]\arrow[d] & \Ba A\arrow[d]\\
    0 \arrow[r] & \Ba(A')
\end{tikzcd}\]
in $\fancycat{D}(\F_2)^{\Z}$. By gluing  (\ref{sqaure: toda ingredient 2})  with   $\epsilon_{q(x)}$  along $\overline{x}$, we obtain the following diagram in $\fancycat{D}(\F_2)^{\Z}$\stepcounter{thm}
\[\leqnomode \tag{\thethm}\label{sqaure: house}\begin{tikzcd}[row sep = small]
 & 0 \arrow[dr, bend left =10]&\\
    \Sigma \bbS^{g,d} \arrow[ur, bend left =10]\arrow[rr,"\overline{x}"]\arrow[d] & &  \Ba A\arrow[d]\\
    0 \arrow[rr] &&  \Ba(A')
\end{tikzcd}\]
Let \[[q(x)]:\Sigma \Sigma \bbS^{g,d}\to \Ba (A')\] be the map in $\fancycat{D}(\F_2)^{\Z}$ given by taking colimit of the $0 \leftarrow \Sigma \bbS^{g,d}\to 0$ part  of
    (\ref{sqaure: house}).
\end{construction}

\begin{exmp}\label{exmp: coproduct of bracket}
    Suppose $A = \E_1(\bigoplus_{i=1}^k\bbS^{g_i,d_i}\{a_i\}\oplus \bbS^{g_i',d_i'}\{b_i\})$ with $d_i=g_i-1>0$, $ d_i'=g_i'-1>0$, and $g_i+g_i' = g_1+g_1'$. If $x$ be determined by $q(x) = \sum_{i=1}^k a_i\otimes b_i$, then $\Ba(A') = \Delta_{A'}$, whose underlying vector space is $\F_2\{1\}\oplus \F_2\{ \{\overline{a_i},\overline{b_i}\}_i\}\oplus \F_2\{[q(x)]\}$, and
    \begin{claim}
        The reduced coproduct of $[q(x)]\in \Delta_{A'}$ is $\sum_{i=1}^k\overline{a_i}\otimes\overline{b_i}$.
    \end{claim}
    \begin{prfoc}
        By a change of grading, i.e. left Kan extend along the symmetric monoidal functor $+:\Z\times \Z\to \Z$,  we may assume \[A = \E_1(\bigoplus_{i=1}^k\bbS^{g_i,d_i,0}\{a_i\}\oplus \bbS^{g_i',0,d_i'}\{b_i\})\] is in fact bigraded. Then for grading reasons the coproduct of $q[x]$ is of the form 
        \[\psi([q(x)]) =1\otimes [q(x)]+ [q(x)]\otimes 1+\sum_{i,j}\lambda_{ij}\overline{a_i}\otimes \overline{b_j},\quad \lambda_{ij}\in \F_2.\]
        The key realisation is that  this coproduct corresponds to a  relation in the homotopy group of $\Cobar(\Delta_{A'})$. And as $\Ba A'$ is diagonally coconnective for our $A'$, by  \cref{prop: bar cobar restrict to equiv}  the map $A'\to \Cobar (\Delta_{A'})$ is an equivalence, so the relation can be checked on $A'$. More precisely, the homotopy group \[\pi_{g,d,d'}(\Cobar (\Delta_{A'})) = \Cotor^{g-d,g-d'}_{\Delta_{A'}}(\F_2,\F_2)(g)\]
        can be computed using the cobar complex of the coalgebra $\Delta_{A'}$. As $\overline{a_i},\overline{b_i}\in \Delta_{A'}$ are primitive, they define 
        classes $[\overline{a_i}], [\overline{b_i}]$ in the homotopy groups of $\Cobar (\Delta_{A'})$, which corresponds to $a_i,b_i$ under the equivalence $A'\cong \Cobar (\Delta_{A'})$. The coproduct $\psi([q(x)])$ imposes the relation $\sum_{i,j}\lambda_{ij} [\overline{a_i}][\overline{b_i}] = 0$ in $\pi_{g_1+g_1',d_1,d_1'}(\Cobar(\Delta_{A'}))$, so $\sum_{i,j}\lambda_{ij}a_ib_j =0$ in $\pi_{g_1+g_1',d_1,d_1'}(A')$. But we know the only non-trivial relation in this tridegree is the relation $x=0$.
    \end{prfoc}
\end{exmp}

\subsection{Formal properties of the truncation functor} We will need \cref{cor: truncation preserves colimit} below for the proof of the main theorem. Put $\fancycat{E}\defeq \fancycat{D}(\F_2)^{\Z}$. First we recall

\begin{lem}\label{lem: truncation is left adjoint}
\autocite[\nopp 2.2.1.9 (i)]{HA} The restriction of $\tau^{\mathrm{diag}}$ to diagonally connective objects $\tau^{\mathrm{diag}}_{\leqslant 0}:\fancycat{E}_{\geqslant0}\to \fancycat{E}^{\heartsuit}$ is left adjoint to inclusion. In particular, it preserves colimits.\hfill $\blacksquare$
\end{lem}

Since $\tau^{\mathrm{diag}}_{\leqslant 0}$ is strong symmetric monoidal and the inclusion $\iota:\fancycat{E}^{\heartsuit}\hookrightarrow\fancycat{E}_{\geqslant 0}$ is lax, they both lifts to map between algebras; sometimes it will be convenient to given them another name, say $\overline{\tau}^{\mathrm{diag}}_{\leqslant 0}$ and $\overline{\iota}$, to distinguish them from their counterpart on the underlying.

\begin{lem}
    $\overline{\tau}^{\mathrm{diag}}_{\leqslant 0}:\Alg_{\E_{n-1}}(\fancycat{E}_{\geqslant 0})\to \Alg_{\E_{n-1}}(\fancycat{E}^{\heartsuit})$ is the left adjoint of $\overline{\iota}$.
\end{lem}

\begin{prf}
To compress the notation, let's denote $\overline{\tau}^{\mathrm{diag}}_{\leqslant0}$ by $\overline{\tau}$ and write $\Alg$ for $\Alg_{\E_{n-1}}$ temporarily.  
Let $\alpha:\id_{\fancycat{E}_{\geqslant 0}}\Rightarrow\iota\tau $ be the unit of the adjunction in \cref{lem: truncation is left adjoint}, then by  the proof of \autocite[\nopp2.2.1.9(i)]{HA} $\alpha$ lifts to a natural transformation $\alpha:\id_{\Alg(\fancycat{E}_{\geqslant 0})}\Rightarrow\iota\tau$, where now $\iota$ and $\tau$ are the induced functor of algebras. We need to check for every $A\in \Alg(\fancycat{E}_{\geqslant 0})$, $B\in \Alg(\fancycat{E}^{\heartsuit})$, \stepcounter{thm}
\[\leqnomode\tag{\thethm}\label{eq: check adjunction}\Map_{\Alg(\fancycat{E}^{\heartsuit})}(\overline{\tau} A, B)\to \Map_{\Alg(\fancycat{E}_{\geqslant 0})}(\overline{\iota}\overline{\tau} A, \iota B)\xrightarrow[]{-\circ \alpha_A} \Map_{\Alg(\fancycat{E}_{\geqslant 0})}( A, \overline{\iota} B)\] is an equivalence. We will need
\begin{claim}
    If $\E_{n-1}^{\heartsuit}$ denotes the free $\E_{n-1}$ algebra functor in the heart (namely the symmetric power if $n>2$ and the tensor power if $n=2$). Then the natural transformation $\E_{n-1}^{\heartsuit}\tau\Rightarrow \tau \E_{n-1}$ of functors $\fancycat{E}_{\geqslant 0}\to \Alg(\fancycat{E}^{\heartsuit})$ is an equivalence.
\end{claim}
\begin{prfoc}
    On objects $\tau$ just take the diagonal homotopy groups. The claim follows from the description of homotopy groups of a free algebra.
\end{prfoc}
Suppose $A =  \E_{n-1}(X)$ is a free algebra, then as $\alpha$ underlies the unit of the adjunction of \cref{lem: truncation is left adjoint}, the map  (\ref{eq: check adjunction})  is equivalent to the composite of the following sequence of equivalences on mapping spaces induced by adjunction.
\begin{align*}
    \Map_{\Alg(\fancycat{E}^{\heartsuit})}(\overline{\tau }A,B)&\cong \Map_{\Alg(\fancycat{E}^{\heartsuit})}(\E_{n-1}^{\heartsuit}(\tau X),B)  & \text{by the claim}\\
    &\cong \Map_{\fancycat{E}^{\heartsuit}}(\tau X, B) \\
    &\cong \Map_{\fancycat{E}_{\geqslant 0}}(X, \iota B) &  \text{by \cref{lem: truncation is left adjoint}}\\
    &\cong \Map_{\Alg(\fancycat{E}_{\geqslant 0})}(R, \overline{\iota} B).
\end{align*}
For general $A$, we can use the simplicial resolution $\Ba(\E_{n-1}, \E_{n-1}, -)$ by free algebras in $\fancycat{E}_{\geqslant0}$. This realises (\ref{eq: check adjunction}) as a colimit of a map of simplicial spaces, where every edge in this map of diagrams is an equivalence.
\end{prf}

\begin{cor}\label{cor: truncation preserves colimit}
    $\overline{\tau}^{\mathrm{diag}}_{\leqslant 0}:\Alg_{\E_{n-1}}(\fancycat{E}_{\geqslant 0})\to \Alg_{\E_{n-1}}(\fancycat{E}^{\heartsuit})$ preserves colimits.\hfill $\blacksquare$
\end{cor}

 \subsection{Proof of \cref{thm: summary theorem of diag hopf}}
 
(i) and (ii) follows from \cref{cor: truncation preserves colimit} by applying the truncation functor $\overline{\tau}^{\mathrm{diag}}_{\leqslant 0}$ to the $\E_{n-1}$-pushout diagram  (\ref{sqaure: bar pushout}). 

Now suppose we are in the situation of (iii).
Consider the adjoint of  (\ref{sqaure: house}) to $\E_{n-1}$-algebras 
\stepcounter{thm}
\[\leqnomode \tag{\thethm}\label{sqaure: alg house}\begin{tikzcd}[sep = small]
 & \E_{n-1}(0) \arrow[dr, bend left =10]&\\
   \E_{n-1}( \Sigma \bbS^{g,d}) \arrow[ur, bend left =10]\arrow[rr]\arrow[d] & &  \Ba A\arrow[d]\\
    \E_{n-1}(0) \arrow[rr] &&  \Ba(A')
\end{tikzcd}\]
By construction of  (\ref{sqaure: house}), the square here is a pushout is $\Alg_{\E_{n-1}}(\fancycat{D}(\F_2)^{\Z})$. So by the pushout pasting law we obtain an an equivalence in $\E_{n-1}$-algebras from $\Sigma_{\E_{n-1}}\E_{n-1}(\Sigma \bbS^{g,d})\sqcup^{\E_{n-1}}\Ba A$ to $\Ba(A')$. The expression of $\psi([q(x)])$ is due to \cref{exmp: coproduct of bracket} and naturality.\qedd

\section{Building the detector}\label{sec: building the detector}
As explained in the introduction, homological stability of $\GL_n(\F_2)$ with $\F_2$-coefficient is controlled by the stability Hopf algebra $\Delta_{\CGL(\F_2)}$ (\cref{def: stability Hopf alg})  of the $\E_\infty$ $\F_2$-module $\CGL(\F_2)$ of\cref{const: algebra from symmetric monoidal groupoid}. An alternative type of application of the stability Hopf algebra construction  is  detecting periodic families. Suppose $R$ is connected and satisfies (SCE), then we have the $\E_\infty$ map $t_R:R\to r$ of (\ref{map: truncation map}), $r = \Cobar(\Delta_R)$. If one can find a Hopf algebra map of the form $\Delta_R\to \Gamma$ such that the target $\Gamma$ has computable Cotor groups, then one can try to detect homotopy classes via basechange along the $\E_\infty$-map $R\to \Cobar (\Gamma)$.   Our $\Gamma$ will be the truncated dual Steenrod algebra $A(1)_*$, and in \S\ref{sec: cgl to a(1)} we describe the map $\Delta_{\CGL(\F_2)}\to A(1)_*$. For this we need a description of $\Delta_{\CGL(\F_2)}$ in small gradings, which can be found in \S\ref{sec: stability Hopf alg of CGL}. where we describe $\Delta_{\CGL(\F_2)}$ up to grading $5$ by finding a cell approximation for $\CGL$ in low degrees and using tools developed in \S\ref{sec: tools}. In \S\ref{sec: homotopy groups of a(1) and a(1)/h_{10}} we compute the homotopy groups of $a(1)\defeq \Cobar(A(1)_*)$  and the cofibre $a(1)/h_{10}$. Our supply of classes in $\CGL(\F_2)/\sigma$ comes from the endomorphism that we will build in \S\ref{sec: self map on CGL/sigma}.
 \textit{Throughout this section  $\CGL\defeq \CGL(\F_2)$.}

\subsection{Stability Hopf algebra of $\CGL(\F_2)$}\label{sec: stability Hopf alg of CGL}

Our main goal is the following result.
\begin{thm}\label{thm: stability Hopf alg F_2 up to grading 5}Let $\delta \defeq [\sigma\mid Q^1(\sigma)]\in \Delta_{\CGL}(3)$ and $\rho\defeq [\sigma \mid Q^2(\sigma)]\in \Delta_{\CGL}(4)$ be the class associated to the  decomposition $\sigma \otimes Q^1(\sigma)$ and $\sigma^2 \otimes Q^2(\sigma)$ respectively (c.f. \cref{thm: summary theorem of diag hopf} (iii)), \cref{fig: basis for CGL} describes a basis for $\Delta_{\CGL}$ up to grading $5$.
    \begin{figure}[h]
        \centering
        {\begin{tikzcd}
	n & 1 & 2 & 3 & 4 & 5 \\
	{\text{basis of }\Delta_{\CGL}(n)} & {\overline{\sigma}} & {\overline{\sigma}^2} & {\overline{\sigma}^3} & {\overline{\sigma}^4} & {\overline{\sigma}^5} \\
	&&& \delta & {\overline{\sigma}\delta} & {\overline{\sigma}^2\delta} \\
	&&&& \rho & {\overline{\sigma}\rho}
	\arrow[dotted, no head, from=1-1, to=1-2]
	\arrow[dotted, no head, from=1-1, to=2-1]
	\arrow[dotted, no head, from=1-2, to=1-3]
	\arrow[dotted, no head, from=1-2, to=2-2]
	\arrow[dotted, no head, from=1-3, to=1-4]
	\arrow[dotted, no head, from=1-3, to=2-3]
	\arrow[dotted, no head, from=1-4, to=1-5]
	\arrow[dotted, no head, from=1-4, to=2-4]
	\arrow[dotted, no head, from=1-5, to=1-6]
	\arrow[dotted, no head, from=1-5, to=2-5]
	\arrow[dotted, no head, from=2-1, to=2-2]
	\arrow[dotted, no head, from=2-2, to=2-3]
	\arrow[dotted, no head, from=2-3, to=2-4]
	\arrow[dotted, no head, from=2-4, to=2-5]
	\arrow[dotted, no head, from=2-4, to=3-4]
	\arrow[dotted, no head, from=2-5, to=2-6]
	\arrow[dotted, no head, from=2-5, to=3-5]
	\arrow[dotted, no head, from=2-6, to=1-6]
	\arrow[dotted, no head, from=3-4, to=3-5]
	\arrow[dotted, no head, from=3-5, to=3-6]
	\arrow[dotted, no head, from=3-5, to=4-5]
	\arrow[dotted, no head, from=3-6, to=2-6]
	\arrow[dotted, no head, from=4-5, to=4-6]
	\arrow[dotted, no head, from=4-6, to=3-6]
\end{tikzcd}}
        \caption{$\bigoplus_{n\leqslant5}\Delta_{\CGL}(n)$}
        \label{fig: basis for CGL}
    \end{figure}
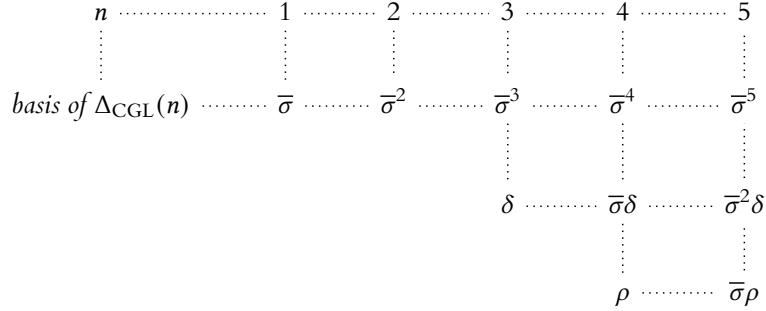
Here $\overline{\sigma}$ is the bar operation (\cref{def: bar operation}) on $\sigma$.
The coproduct structure in this range is determined by 
\begin{align*}
    \psi(\overline{\sigma})&= 1\otimes \overline{\sigma} + \overline{\sigma}\otimes 1\\
    \psi(\delta)& = 1\otimes \delta + \overline{\sigma}\otimes \overline{\sigma}^2 +\delta\otimes 1\\
    \psi(\rho)&= 1\otimes \rho + \rho\otimes 1.
\end{align*}
For future's convenience, we spell out the remaining coproducts in grading $\leqslant 4$: $\overline{\sigma}^2,\overline{\sigma}^4$ are primitive, $\psi(\overline{\sigma}^3) =1\otimes \overline{\sigma}^3 + \overline{\sigma}\otimes \overline{\sigma}^2 + \overline{\sigma}^3\otimes 1$, and $\psi(\overline{\sigma}\delta) =1\otimes \overline{\sigma}\delta + \overline{\sigma}\otimes \overline{\sigma}^3 + \delta\otimes \overline{\sigma} + \overline{\sigma}\otimes \delta + \overline{\sigma}^2 \otimes \overline{\sigma}^2 + \overline{\sigma}\delta\otimes 1$.
\end{thm}

We start by  building an $\E_\infty$-cell complex $X_5$ together with a $(5,4)$-connected $\E_\infty$-map $f_5:X_5\to \CGL$. 
Start with a representative $f_0:X_0\defeq \E_\infty(\bbS^{1,0}\sigma)\to \CGL$ of the class $\sigma\in \pi_{1,0}(\CGL)$. Suppose we have built $f_i:X_i\to \CGL$, then: If $\pi_{*,*}(f)$ is an isomorphism in the range $\{(g,d)\mid g\leqslant 5,d\leqslant 3\}$, we stop, and it will happen that this terminating $f_i$ is already $(5,4)$-connected. Otherwise, we look at the lowest row $d$ then the lowest column $g$ at which $\pi_{g,d}(f_i)$ fails to be an isomorphism. There are two cases. 

(1):  
If the failure is due to $\ker\pi_{g,d}(f_i)\neq 0$, then pick a non-zero $x_i$ from the kernel. We take  \[f_{i+1}:X_{i+1}\defeq X_i/\!/_{\E_\infty}^{g,d+1}x_i\to \CGL.\]
to be be $\E_\infty$-map induced by any representative 
$\E_\infty(\bbS^{g,d}\{x_i\})\to X_i$  of $x$ and any null homotopy of its composite with $f_i$.
    We then compute $\pi_{*,*}(X_{i+1})$ via the cell-attachment spectral sequence $\{E_r(\fil^{\mathrm{cell}}X_{i+1})\}_r$ of  \autocite[\S10.3.2]{GKRW18a}. As a trigraded algebra \[E_1(\fil^{\mathrm{cell}}X_{i+1})\cong \pi_{*,*,*}(0_!X_i)\otimes \pi_{*,*,*}(1_!\bbS^{g,d+1}\{\beta_i\}),\]where $i_!X$ of an object $X$ is the graded object that is $X$ in grading $i$ and zero elsewhere, and $\beta_i$ is the cell such that  $d^1\beta_i =x_i$. In this calculation it will always be the case that, in the range $\{(g,d)\mid g\leqslant 5, d\leqslant 4\}$, the $E_2$ page coincides with the $E_\infty$ page  where both be identified with the algebra quotient $\pi_{*,*,*}(0_!X_i)/(0_!x_i)$. This convenience is essentially because the Dyer-Lashof of $\beta_i$ (including $\beta_i^2$) has grading $>5$.

    (2): If $\im \pi_{g,d}(f_i)\neq \pi_{g,d}(\CGL)$,
 we pick a vector space splitting $\pi_{g,d}(\CGL) \cong \im \pi_{g,d}(f_i)\oplus V_i$, a basis $\{x'_{ij}\}$ of $V_i$, and take 
    \[f_{i+1}:X_{i+1}\defeq X_i\sqcup^{\E_\infty}\E_\infty(\bigoplus_{x'_{ij}}\bbS^{g,d}x_i')\to \CGL.\]
    Since finite coproduct in $\E_\infty$-algebras is be computed by the tensor product, we can read off $\pi_{*,*}(X_{i+1})$. Our choices of the $x_i$ and $x_{ij}'$ in each step is list below.
\begin{center}
\begin{tabular}{l l l}
   \textbf{Step 1. }   &  $x_0 = \sigma Q^1(\sigma)$ & $f_1:X_1 = \E_\infty(\bbS^{1,0}\sigma)/\!/^{3,2}_{\E_\infty}\sigma Q^1(\sigma)\to \CGL$.  \\
     \textbf{Step 2. }& $x_1 = \sigma^2 Q^2(\sigma)$ & $ f_2:X_2  = X_1/\!/^{4,3}_{\E_\infty}\sigma^2 Q^2(\sigma)\to \CGL$.\\
     \textbf{Step 3.}   & $x_2 = \sigma Q^3(\sigma)$ & $ f_3:X_3  = X_1/\!/^{3,4}_{\E_\infty}\sigma  Q^3(\sigma)\to \CGL$.\\
     \textbf{Step 4.} &  $x_3' = \nu_1,\nu_2$   & $ f_4:X_4 =X_3\sqcup^{\E_\infty}\E_\infty(\bbS^{3,3}\nu_1\oplus \bbS^{3.3}\nu_2)\to\CGL$.\\
     \textbf{Step 5.}  & $x_4= \sigma \nu_1,\sigma \nu_2$ & $ f_5:X_5 =X_4/\!/_{\E_\infty}^{4,4}(\sigma\nu_1,\sigma\nu_2)\to \CGL$. 
\end{tabular}
\end{center}

From step 3 onwards, all the cells attached have slope $\geqslant 1$. So \cref{thm: summary theorem of diag hopf} (i) tells us that $X_2\to X_5$ induces an isomorphism on the stability Hopf algebra, and hence \cref{thm: stability Hopf alg F_2 up to grading 5} will follow from the following two propositions.

\begin{prop}
    Let $\overline{\sigma}, \delta,\rho\in \Delta_{X_2}$ be defined analogously to their counterpart in \cref{thm: stability Hopf alg F_2 up to grading 5}. Then \[\Delta_{X_2} = \F_2[\sigma, \delta,\rho],\] with $\sigma,\rho$ both primitive, and  $ \psi(\delta)= 1\otimes \delta + \overline{\sigma}\otimes \overline{\sigma}^2 +\delta\otimes 1$.
\end{prop}
\begin{prf}
    Apply \cref{thm: summary theorem of diag hopf} (iii) iteratively. To see $\rho$ is primitive we use \cref{thm: bar and operations} (iii).
\end{prf}

\begin{prop}\label{prop: connectivity of Delta_f_2}
    The Hopf algebra map $\Delta_{f_2}:\Delta_{X_2}\to \Delta_{\CGL}$  is a bijection in grading $\leqslant 5$.
\end{prop}
\begin{prf} Applying the Hurewicz theorem in the form of \cref{cor: connectivity passing to bar} to the $(5,4)$-connected map $f_5:X_5\to \CGL$ we see that $\Delta_{f_5}$ is $5$-connected. Since $\Delta_{X_2}\to \Delta_{X_5}$ is an isomorphism, $\Delta_{f_2}$ is also $5$-connected. To see $\Delta_{f_2}$ at grading 5 is in fact a bijection, we compute the coproduct
\[\psi(\Delta_{f_2}(\kappa))  = \Delta_{f_2}(\psi (\kappa))\in \Delta_{\CGL}\]
for a generic element $\kappa \defeq \lambda_1\overline{\sigma}^5 + \lambda_2\overline{\sigma}^2\delta + \lambda_3\overline{\sigma}\rho$ ($\lambda_i\in \F_2$) in $\Delta_{X_2}(5)$ and use injectivity of $\Delta_{f_2}$ in lower grading.\end{prf}

\subsection{The map $\cgl\to a(1)$.}\label{sec: Delta BGL -> a(1)}\label{sec: cgl to a(1)}

The dual Steenrod algebra $A_*$ is a commutative Hopf algebra with underlying algebra $\F_2[\xi_1,\xi_2,\cdots ]$ and comultiplication specified by $\psi(\xi_n) = \sum_{i=0}^n\xi_{n-i}^{2^i}\otimes \xi_i$. For each $n\geqslant0$, the algebra quotient $A(n)_*\defeq  A_*/(\xi_1^{2^{n+1}}, \cdots ,\xi_{n+1}^2,\xi_{n+2},\xi_{n+3},\cdots )$ is in fact a quotient Hopf algebra. In this section we construct a map of dual Hopf algebras $\varphi: A(1)_*^{\vee}\to \Delta_{\CGL}^\vee$. Our detector is  then the $\E_\infty $-map\footnote{Here we made identifications $\Delta_{\CGL}^{\vee\vee}\cong \Delta_{\CGL}$ and $A(1)_*^{\vee\vee}\cong A(1)_*$. Both Hopf algebras are finite dimensional in each degree. For $\Delta_{\CGL}$, this  comes from identifying $\Delta_{\CGL}(g) \cong H_0(\GL_g(\F_2);\St^{\E_1}(g)\otimes\F_2)$ (\autocite[\S1.6]{RW25}, \autocite[189]{GKRW18a}) as the group homology of a finite group. The Hopf algebra $A(1)_*$ has total dimension $8$. } \[\Cobar(\varphi^\vee):\mathrm{cgl}\to a(1),\]  where $\cgl\defeq \Cobar(\Delta_{\CGL})$ and $a(1)\defeq \Cobar(A(1)_*)$.
\begin{rem}
The lower $A(n)$'s have meanings in homotopy theory: for $n= 0,1,2$, $A(n)$ computes stable $H\F_2$ cohomology operations relative to $H\F_2$, $\mathrm{ko}$ and $\mathrm{tmf}$ respectively. However we suspect that the appearance of $A(1)$ in our situation is  merely a coincidence.
\end{rem}

The dual $A(1)_*^\vee$ is the sub Hopf algebra $A(1) = \<\Sq^1,\Sq^2>$ of the Steenrod algebra.   As an algebra $A(1)$ is the quotient of the free associative algebra on symbols $\Sq^1$, $\Sq^2$ by the two sided ideal on $\Sq^{1,1}$ and $\Sq^{2,2}+\Sq^{1,2,1}$. We now specify elements $s_1,s_2\in \Delta_{\CGL}^\vee$ such that  \[\varphi:\Sq^i\in A(1)\mapsto s_i\in \Delta_{\CGL}^\vee\] determines a map of Hopf algebras.
Recall from (\cref{thm: stability Hopf alg F_2 up to grading 5}) that $B\defeq \{\overline{\sigma}, \overline{\sigma}^2,\overline{\sigma}^3,\delta,\overline{\sigma}^4,\overline{\sigma}\delta, \rho\}$ is a basis for $\bigoplus_{n=1}^4\Delta_{\CGL}(n)$.  Let  $B^\vee = \{b^\vee\mid b\in B\}$ denote the dual basis. We take $s_1,s_2$ to be $(\overline{\sigma}^1)^{\vee},(\overline{\sigma}^2)^{\vee}$.

\begin{lem} \cref{fig: basis for dual hopf alg} describes a basis for $\Delta_{\CGL}^\vee(n)$, $n\leqslant 4$. Moreover, the element $s_1$ is primitive, and \[\psi(s_2)=1\otimes s_2 + s_1\otimes s_1 +s_2\otimes 1.\]

\begin{figure}[h]
    \centering
{\begin{tikzcd}
	{n} & 1 & 2 & 3 & 4 \\
	{\text{basis of }\Delta_{\CGL}^\vee(n)} & {s_1} & {s_2} & {s_1s_2=(\overline{\sigma}^3)^\vee+ \delta^\vee} & {(\overline{\sigma}^4)^\vee} \\
	&& {(s_1^2=0)} & {s_2s_1=(\overline{\sigma}^3)^\vee} & {s_2^2 = s_1s_2s_1= (\overline{\sigma}\delta)^\vee} \\
	&&&& {\rho^\vee}
	\arrow[dotted, no head, from=1-1, to=1-2]
	\arrow[dotted, no head, from=1-1, to=2-1]
	\arrow[dotted, no head, from=1-2, to=1-3]
	\arrow[dotted, no head, from=1-2, to=2-2]
	\arrow[dotted, no head, from=1-3, to=1-4]
	\arrow[dotted, no head, from=1-3, to=2-3]
	\arrow[dotted, no head, from=1-4, to=1-5]
	\arrow[dotted, no head, from=1-4, to=2-4]
	\arrow[dotted, no head, from=1-5, to=2-5]
	\arrow[dotted, no head, from=2-1, to=2-2]
	\arrow[dotted, no head, from=2-2, to=2-3]
	\arrow[dotted, no head, from=2-4, to=2-3]
	\arrow[dotted, no head, from=2-4, to=2-5]
	\arrow[dotted, no head, from=2-4, to=3-4]
	\arrow[dotted, no head, from=2-5, to=3-5]
	\arrow[dotted, no head, from=3-4, to=3-5]
	\arrow[dotted, no head, from=3-5, to=4-5]
\end{tikzcd}}
    \caption{$\bigoplus_{n\leqslant 4}\Delta_{\CGL}^\vee(n)$}
    \label{fig: basis for dual hopf alg}
\end{figure}
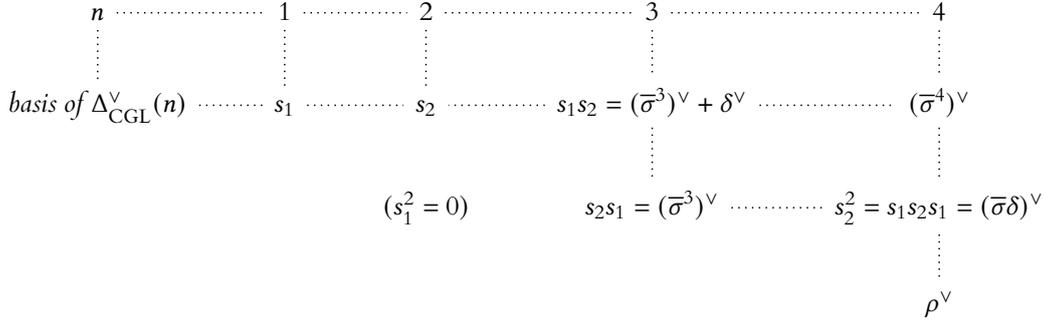

\end{lem}
\begin{prf}
    Since taking linear dual exchanges primitives with indecomposables, $s_1$ is primitive and $s_2  = (\overline{\sigma}^2)^\vee$ is not. The only possibly non-zero reduced coproduct for $s_1$ is $s_1\otimes s_1$. This also implies that $s_1^2=0$ as $s_1^2$ is primitive and $\Delta_{\CGL}(2)^\vee$ is $1$-dimensional. To see the identities in the $n=3$ column, let $\<-,->:\Delta_{\CGL}^\vee\otimes \Delta_{\CGL}\to \F_2$ be the evaluation map. One has  
\setlength{\tabcolsep}{12pt}
\[\begin{tabular}{l l}
    $\<s_2s_1, \overline{\sigma}^3>=1$ & $\<s_2s_1,\delta>= 0$ \\
   $\<s_1s_2+s_2s_1,\overline{\sigma}^3>= 0$  & $\<s_1s_2+s_2s_1,\delta>= 1$ 
\end{tabular}\]
which implies $s_2s_1= (\overline{\sigma}^3)^\vee$ and $s_1s_2+s_2s_1 = \delta^\vee$. For the $n=4$ column, both $s_2^2$ and $s_1s_2s_1= s_1(s_2s_1) = s_1(\overline{\sigma}^3)^\vee$ evaluates to $1$ against $\overline{\sigma}\delta$ and $0$ against $\overline{\sigma}^4, \rho$, so $s_2^2=s_1s_2s_1=(\overline{\sigma}\delta)^\vee$. 
\end{prf}

The relations $s_1^2=0$ and $s_2^2=s_1s_2s_1$ in $\Delta_{\CGL}^\vee$ shows that $\varphi$ defines a map of algebras.  To see $\varphi$ is a map of Hopf algebras, we check that $\psi\varphi = (\varphi\otimes\varphi) \psi$ holds on the algebra generators $\Sq^1$ and $\Sq^2$. This is indeed the case by our formula for $\psi(s_1)$ and $\psi(s_2)$.

\begin{rem}\label{rem: description of phi dual}Since
    $\varphi^\vee:\Delta_{\CGL}\to A(1)_* = \F_2[\xi_1,\xi_2]/(\xi_1^4,\xi_2^2)$ is a bijection in grading $1$ and preserves coproduct, we must have
    \[\varphi^\vee(\overline{\sigma}) = \xi_1,\qquad \quad \varphi^\vee(\delta)  = \xi_2+\xi_1^3.\]
\end{rem}

\subsection{Homotopy groups of $a(1)$ and $a(1)/h_{10}$.} \label{sec: homotopy groups of a(1) and a(1)/h_{10}}Suppose $A = \mathbb{k} \oplus \overline{A}$ is a connected\footnote{By this we mean that the unit map $\mathbb{1}\to A$ is an isomorphism in gradings $\leqslant0$. In particular, $A$ is canonically augmented and the augmentation ideal $\overline{A} $ is supported in positive gradings.} graded Hopf algebra over a field $\bbk$, then the homotopy group $a\defeq\Cobar(A)$ is the Cotor groups of $A$: \stepcounter{thm}
\[\tag{\thethm}\label{eq: homotopy group of cobar = cotor}\leqnomode\pi_{g,d}(a) = \Cotor_{g-d}^{A}(\mathbb{k},\mathbb{k})(g).\]
In particular, $\pi_{*,*}(a)$ can be accessed via the classical cobar complex. When $A = A(1)_*$, the primitive elements $\zeta_1, \zeta_1^2$ and $\zeta_2$ defines cycles $\underline{\zeta_i^{2^j}}$ in the cobar complex. We denote the corresponding class in $\pi_{*,*}(a(1))$ by $h_{ij}$. We will see that $h_{10}$ correspond to the the class $\sigma$ under the map $\CGL\to \cgl\to a(1)$. The purpose of this section is to describe $\pi_{*,*}(a(1))$ with its algebra structure and compute $\pi_{*,*}(a(1)/h_{10})$ as a $\pi_{*,*}(a(1))$-module.

\subsubsection{The Ivanovsky spectral sequence}The filtration
\[(\fil^{\mathrm{aug}}A)_p\defeq \begin{cases*}
    A & for  $p\geqslant 0$\\
    (\overline{A})^{-p}  & for $p<1$
\end{cases*}\]
on $A$ induces a filtration on the cobar complex complex on $A$. Suppose $A$ is finite-dimensional in each grading, then we obtain a spectral sequence converging strongly to  $\Cotor^A(\mathbb{k},\mathbb{k})$. This spectral sequence was introduced by  L. N. Ivanovsky \autocite{Iva64} with $A$ being the dual Steenrod algebra\footnote{The May spectral sequence is obtained by filtering the bar complex of the Steenrod algebra by augmentation ideals.}; a good reference is \autocite{BS94}.

This spectral sequence has a natural $W_{n-1}$-structure, with  $n=\infty$ if $A$ is commutative and $n=1$ otherwise, given by the 
 identification (\ref{eq: homotopy group of cobar = cotor}): We can view $\fil^{\mathrm{aug}}A$ as a discrete $\E_1$-coalgebra-$\E_n$-algebra in filtered objects,  so its Cobar \[\fil^{\mathrm{aug}}a\defeq \Cobar (\fil^{\mathrm{aug}}A)\] is an $\E_{n+1}$-filtration on $a$.  This filtered object is called the \textit{underived canonical multiplicative filtration} in \autocite[\S8.4]{RW25}. The associated spectral sequence $\{E^r(\fil^{\mathrm{aug}}a)\}_r$ is the Ivanovsky spectral sequence.

 \begin{exmp}\label{exmp: E^1 of a(1)}
     When $A = A(1)_*$, the $E^1$ page is the algebra $\F_2[h_{10}, h_{11},h_{20}]$, and the
     Dyer-Lashof operations on the $E^1$ page are such that $h_{ij} = Q_1^{\circ j}(h_{i0})$; see the discussion in \autocite[\S8.4]{RW25}.
 \end{exmp}

Unwinding the filtration on the cobar complex of $A$, we see that the differentials in $\{E^r(\fil^{\mathrm{aug}}a)\}_r$ can be imported from the coproducts on $A$ as follows. Suppose $x\in \fil^{\mathrm{aug}}_r A$ ($r<0$) has reduced coproduct \[\overline{\psi}(a) = \sum x'_i\otimes x''_i \in \bigcup_{p+q\leqslant -n+r} \fil^{\mathrm{aug}}_pA\otimes \fil^{\mathrm{aug}}_qA,\] then $x\in E^1(\fil^{\mathrm{aug}}a)$ survives to  the $E^n$ page and $d^n([a]) $ is represented by the image of $\sum x_i'x_i''$ mod $\fil^{\mathrm{aug}}_{r-n-1}A$. We will compute $\pi_{*,*}(a(1))$ using this spectral sequence, and compute $\pi_{*,*}(a(1)/h_{10})$ via the spectral sequence associated to the cofibre of $\fil^{\mathrm{aug}}a(1)$ by a filtered lift of $h_{10}$.

\subsubsection{The homotopy groups of $a(1)$} 
\begin{figure}[h]
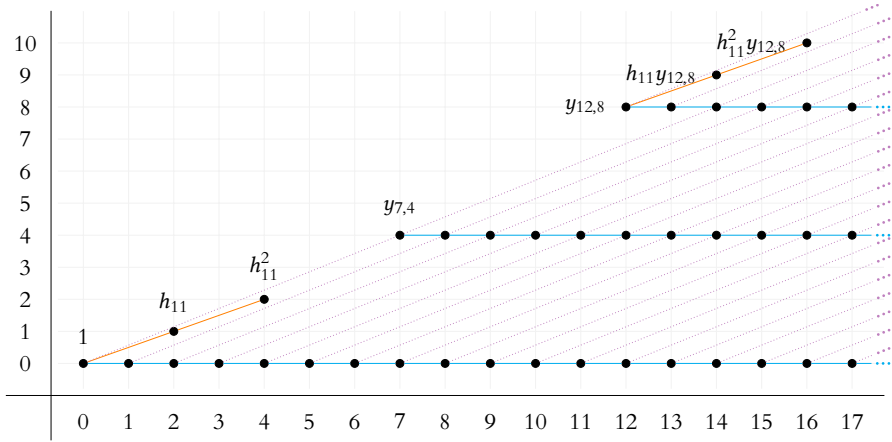

\begin{center}
\adjustbox{scale=0.85}{
 \NewSseqGroup\bottomgroup {} {
\class(0,0)
\class(2,1)
\class(4,2)
\structline[orange](0,0)(2,1)
\structline[orange](2,1)(4,2)
}
 \NewSseqGroup\horizontalgroup {} {
 \foreach \x in {0,...,30}  {
\class(\x,0)
}
 \foreach \x in {0,...,17} {
\structline[cyan](\x,0)(\x+1,0)
}
}

\begin{sseqpage}[ classes = fill,  xscale = 0.7,  yscale = 0.5, grid = go, grid color =gray!10, x range = {0}{17}, y range = {0}{10}, x axis tail = 0.7cm, y axis tail = 0.7cm, x tick gap = 0.4cm, y tick gap = 0.4cm ]
\bottomgroup(0,0)
\bottomgroup(12,8)
\horizontalgroup (1,0)
\structline[cyan] (0,0)(1,0)
\horizontalgroup (7,4)
\horizontalgroup (13,8)
\horizontalgroup (19,12)
\structline[cyan](12,8)(13,8)
\foreach \x  in {0,...,10} {
\structline[violet!50, densely dotted](7+\x,4)(14+\x,8)
}
\foreach \x  in {0,...,17} {
\structline[violet!50, densely dotted](0+\x,0)(7+\x,4)
}
\foreach \x  in {0,...,6} {
\structline[violet!50, densely dotted](12+\x,8)(19+\x,12)
}
\classoptions["1" above](0,0)
\classoptions["h_{11}" above ](2,1)
\classoptions["h_{11}^2" above](4,2)
% \classoptions["h_{10}^{16}" above] (16,0)
% \classoptions["h_{10}^{9}y_{7,4}" below] (16,4)
\classoptions["y_{7,4}" above] (7,4)
% \classoptions["h_{10}y_{7,4}" below] (8,4)
\classoptions["y_{12,8}" left] (12,8)
\classoptions["h_{11}y_{12,8}" left] (14,9)
\classoptions["h_{11}^2y_{12,8}" left] (16,10)
% \classoptions["h_{10}^5y_{12,8}" above] (17,8)
% \classoptions["h_{10}y_{7,4}^2" above ] (15,8)

\end{sseqpage}}
\end{center}
    \caption{$\pi_{g,d}(a(1))$ for $g\leqslant 17, d\leqslant 10$. Blue / orange / dotted purple structure line $\defeq h_{10}$ / $ h_{11}$ /  $y_{7,4}$-extension.}
    \label{fig: homotopy groups of a(1)}
\end{figure}

Throughout this part we put $E^r \defeq E^r(\fil^{\mathrm{aug}}a(1))$. Recall the description  $E^1_{*,*,*} = \F_2[h_{10}, h_{11}, h_{20}]$ from \cref{exmp: E^1 of a(1)}. 
\begin{thm}\label{thm: homotopy groups of a(1)}
    The spectral sequence $\{E^r\}_r$ collapses at $E^5$, and \stepcounter{thm}
    \[\leqnomode\tag{\thethm}\label{eq: alg of a(1)}E^\infty_{*,*,*}=E^5_{*,*,*} = \frac{\F_2[h_{10}, \;h_{11},\;y_{7,4} , \;y_{12,8}]}{(h_{10}h_{11}, \;h_{11}^3,\; h_{11}y_{7,4},\; y_{7,4}^2 -h_{10}^2 y_{12,8})}\]
    where $y_{7,4}\in E^5 = E^\infty$ is represented by $h_{10}h_{20}^2\in E^1_{7,4,-3}$, and $y_{12,8}$ by $h_{20}^4\in E^1_{12,8,-4}$. Since every $E^\infty_{g,d,*}$ is $1$-dimensional, the filtration on  $\pi_{*,*}(a(1))$ is trivial, and hence $\pi_{*,*}(a(1))$
    is canonically identified with the $E^\infty$-page.
\end{thm}
\cref{fig: homotopy groups of a(1)} describes $\pi_{*,*}(a(1))$ in low degrees. Note that $\pi_{*,*}(a(1))$ is $(12,8)$-periodic: additively it is the free $\F_2[y_{12,8}]$-module on $\F_2[h_{10}, h_{11},y_{74}]/(h_{10}h_{11},h_{11}^3, y_{7,4}^2)$.

We begin by collecting some differentials in this spectral sequence. Part (i) below follows from observing that $d^r:E^r_{g,d,f}\to E^r_{g,d-1,f-r}$ lowers $g+f$ by $r$, and $\abs{h_{10}} =(1,0,-1) $, $\abs{h_{20}} = (3,2,-1)$, $\abs{h_{11}} = (2,1,-2)$ all have even $g+f$.
The reduced coproduct $\overline{\psi}(\zeta_2) = \zeta_1\otimes \zeta_1^2$ gives us (iii) below. One can obtain (ii) by either noting that $\zeta_1$ is primitive or by considering degrees. The rest are consequences of the the $W_\infty$-structure on the spectral sequence  \autocite[\S16.6]{GKRW18a}.
\begin{lem}\label{lem: differentials}
    \begin{enumerate}[(i),before*=\leavevmode\vspace{-0.17\baselineskip}]
    \item There are no non-trivial differentials $d^r$ for odd $r$.
    \item $h_{10}\in E^1_{1,0,-1}$ is a permanent cycle.
    \item  $h_{20}\in E^1_{3,2,-1}$ survives to $E^2$, and $d^2(h_{20}) = h_{10}h_{11}$.
          \item $h_{20}^2\in E^1_{6,4,-2}$ survives to $E^4$ and $d^4([h_{20}^2]) = h_{11}^3$.
    \item  $h_{20}^4\in E^1$, and hence every $h_{20}^{4k}\in E^1$, $k\geqslant0$,  is a permanent cycle. 
    \item For every $k\geqslant 0$ and $1\leqslant i\leqslant 3$,  the class $h_{20}^{4k+i}\in E^1$ survives to $E^{3+(-1)^i}$ with
    \begin{align*}
        d^2(h_{20}^{4k+1})&=h_{10}h_{11}\cdot h_{20}^{4k},  & i=1\\
        d^4(h_{20}^{4k+2})&=h_{11}^3\cdot h_{20}^{4k},& i=2\\
        d^2(h_{20}^{4k+3})&=h_{10}h_{11}\cdot h_{20}^{4k+2}, & i=3.
    \end{align*}
    \item $h_{10}h_{20}^2\in E^1_{7,4,-3}$ survives to page $4$, and  $d^4(h_{10}h_{20}^2 ) =[h_{10}h_{11}^3]=0$.\hfill $\blacksquare$
    \end{enumerate}
\end{lem}

\begin{lem}\label{lem: E_infty a(1)} 
$E^3=\F_2[h_{10}, h_{11}, h_{20}^2]/(h_{10}h_{11})$.
\end{lem}
\begin{prf} By \cref{lem: differentials}$E^1 = E^2$ and the  $d^2$-differential is given by  
    \[d^2(h_{10}^ih_{11}^{j}h_{20}^k) =\begin{cases}
    0 &, k \equiv 0\;\mathrm{mod} 2\\
    h_{10}^{i+1}h_{11}^{j+1}h_{20}^{4l+\epsilon} &, k=4l+\epsilon+1, \;\epsilon\in \{0,1\}\\
\end{cases}.\]
Now the result follows from 
\begin{claim}
    For every $g,d,f$, one has $\dim E^2_{g,d,f} = \dim E^1_{g,d,f}\leqslant 1$. 
\end{claim}
\begin{prfoc}
Suppose $h_{10}^ah_{11}^bh_{20}^c=h_{10}^{a'}h_{11}^{b'}h_{20}^{c'}$. Then 
\[a+2b+3c = a'+2b'+3c',\quad b+2c = b'+2c',\quad -a-2b-c= -a'-2b'-c'.\]
by equating the tridegree, from which we obtain $c = c'$, then $b =b'$, then $a = a'$.
\end{prfoc}
\end{prf}

\begin{prf}[of \cref{thm: homotopy groups of a(1)}]
    The monomials $h_{10}^ih_{20}^{2k}$ and $h_{11}^jh_{20}^{2k}$, $i,k\geqslant 0$, $j=1,2$, form a basis for the $E^4$-page. The $d^4$-differential with respect to this basis is given by 
    \begin{align*}
    d^4(h_{10}^ih_{20}^{4k})&=0\\
    d^4(h_{11}^jh_{20}^{4k})&=0\\
    d^4(h_{11}^jh_{20}^{4k+2})& = h_{11}^{j+3}h_{20}^{4k}\\
    d^4(h_{10}^ih_{20}^{4k+2})&=h_{10}^ih_{11}^3h_{20}^{4k}=0
\end{align*}
Since $E^4$ is a quotient of $E^2$, the fact that $\dim E^2_{g,d,f}\leqslant 1$, as we saw in the proof of \cref{lem: E_infty a(1)}, tells us $\dim E^4_{g,d,f}\leqslant 1$. In particular there cannot be two non-trivial differential starting at the same tridegree, and hence we get the desired description of $E^5_{*,*,*}$.  To finish, we need to show that there are  no non-trivial differentials $d^r$ when $r\geqslant 5$, and we only need to check on the generators of (\ref{eq: alg of a(1)}). Since $h_{10}$ is a permanent cycle (\cref{lem: differentials}(ii)), by the $W_\infty$-structure on the spectral sequence so is $h_{11} = Q^1(h_{10})$. For homological degree reason  $d^2(y_{7,4})$ and $d^r(y_{12,8})$  must be zero (compare \cref{fig: homotopy groups of a(1)}).
\end{prf}

\begin{rem}
    $\pi_{*,*}(a(1))$ is the $E_2$ page of the 2-local Adams spectral sequence of the connective real $K$-theory $\mathrm{ko} $.
\end{rem}

\subsubsection{The homotopy groups of $a(1)/h_{10}$.} Let $\overline{E}^r = E^r(\fil^{\mathrm{aug}}a(1)/\widetilde{h_{10}})$, where $\widetilde{h_{10}}\in \pi_{1,0,-1}(\fil^{\mathrm{aug}}a(1))$ is defined by $h_{10}$ in the Bockstein spectral sequence. 
By the $\{E^r\}_r$-module structure on $\{\overline{E}^r\}_r$, the class $[h_{20}]\in \overline{E}^1$ is a permanent cycle. 
Since $\overline{E}^1= E^1/(h_{10})$, it follows from the claim in the proof of \cref{lem: E_infty a(1)} that every $\overline{E}^1_{g,d,*}$ and hence every $\overline{E}^r_{g,d,*}$ is either zero or 1-dimensional, in particular $[h_{20}]$ corresponds to a class $z_{3,2}\in \pi_{3,2}(a(1)/h_{10})$. Let $z_{0,0}\in \pi_{0,0}(a(1)/h_{10})$ be the image of $1$. 

\begin{thm}\label{thm: homotopy of a(1)/h10}
    Let $\mathscr{L}\subseteq \pi_{*,*}(a(1)/h_{10})$ be the subspace spanned by the following set of classes \[z_{0,0},\quad h_{11}.z_{0,0},\quad h_{11}^2.z_{0,0},\quad  z_{3,2},\quad h_{11}.z_{3,2},\quad h_{11}^2.z_{3,2}.\]
    Then the following describes $\pi_{*,*}(a(1)/h_{10})$ as a $\pi_{*,*}(a(1))$-module. See also \cref{fig:lightening flash}.
    \begin{enumerate}[(i)]
        \item The classes listed above are all non-zero, and \[\pi_{*,*}(a(1)/h_{10}) = \F_2[y_{12,8}].\mathscr{L}\] as a left $\F_2[h_{11},y_{12,8}]/(h_{11}^3)\subseteq \pi_{*,*}(a(1))$-module.
        \item $h_{10}.z_{3,2}$ is non-zero, so in fact $h_{10}.z_{3,2} = h_{11}^2.z_{0,0}$. 
    \end{enumerate}
\end{thm}
\begin{figure}[h]
\centering
{
        \NewSseqGroup\mygroup {} {
\class(0,0)
\class(2,1)
\class(4,2)
\structline[orange](0,0)(2,1)
\structline[orange](2,1)(4,2)
}
\begin{sseqpage}[ classes = fill,  xscale = 0.9,  yscale = 0.8, grid = go, grid color =gray!10]
\mygroup(0,0)
\mygroup(3,2)
\classoptions["z_{0,0}" right](0,0)
\classoptions["h_{11}.z_{0,0}" right](2,1)
\classoptions["{h_{11}^2.z_{0,0} = h_{10}.z_{3,2}}" right](4,2)
\classoptions["z_{3,2}" left](3,2)
\classoptions["h_{11}.z_{3,2}" left](5,3)
\classoptions["h_{11}^2. z_{3,2}" left](7,4)
\structline[cyan](3,2)(4,2)
\end{sseqpage}}
    \caption{The \textit{lightning flash} $\mathscr{L}$: $\pi_{g,d}(a(1)/h_{10}$, $g\leqslant7$, $d\leqslant 4$. $\pi_{*,*}(a(1)/h_{10}) =\F_2[y_{12,8}].\mathscr{L} $.}
    \label{fig:lightening flash}
\end{figure}

\begin{prf} By the observation we made at the beginning of this part, $\pi_{g,d}(a(1)/h_{10})= \overline{E}^\infty_{g,d,*}$. 
Using the module structure we see that $\overline{E}^1 = \cdots = \overline{E}^4$, and
\begin{align*}
    \overline{d}^r([h_{11}^ih_{20}^{4k}])&=0\quad\text{ for all }r>0,\\
    \overline{d}^4([h_{11}^ih_{20}^{4k+2}]) & =[ h_{11}^{i+3}h_{20}^{4k}],  & i,k\geqslant 0.
\end{align*}
We also have that $\overline{d}^2([h_{20}])=0$. In fact $\overline{d}^r([h_{20}])=0$ for all $r$ because the target lands in tridegree $(3,1,*)$, and there is nothing there even in $\overline{E}^1$. So the remaining $\overline{d}^4$ differentials are
\begin{align*}
    \overline{d}^4([h_{11}^ih_{20}^{4k+1}]) &=\overline{d}^4 (h_{11}^ih_{20}^{4k}.[h_{20}])=0\\
    \overline{d}^4([h_{11}^i h_{20}^{4k+3}]) & = \overline{d}^4 (h_{11}^ih_{20}^{4k+2}.[h_{20}]) = [h_{11}^{i+3}h_{20}^{4k+1}], & k\geqslant 0.
\end{align*}
To conclude (i) we only need $\overline{d}^r([h_{11}^ih_{20}^{4k+1}]) =0$ for $r>4$ as well. The bidegree $(12k+2i+3,8k+i+2)$ of $[h_{11}^ih_{20}^{4k+1}]$ has  $2d-g =4k+1$. The target of the differential has $2d-g = 4k-1$. The other potentially non-zero bidegrees are the  $(12k'+2i', 8k'+i')$, which have $2d-g  =4k'$.

For (ii), we employ a standard trick for detecting hidden extension. Suppose for a contradiction that $h_{10}.z_{3,2}=0$. Then as $h_{11}\cdot h_{10}=0\in \pi_{*,*}(a(1))$, we can form the Toda bracket $\<h_{11}, h_{11}, z_{3,2}>$ in the setting of \autocite[\S3.3]{BK24} (with $(X'',X',X)=(a(1),a(1),a(1)/\sigma)$).  
Moss' convergence theorem in the form of \autocite[Theorem 1.1]{BK24} applied to $(\mathfrak{a},\mathfrak{a}', \mathfrak{a}'') = (h_{11}, h_{10}, [h_{20}])$ and $r = 2$ tells us that the Massey product $\<h_{11}, h_{10}, [h_{20}]>\subset E^2(\fil^{\mathrm{aug}}a(1)/\widetilde{\sigma})\defeq \overline{E}^2$ contains a permanent cycle that converges some class in the Toda bracket. But  $\<h_{11}, h_{10}, [h_{20}]> =\{h_{20}.[h_{20}]\}$ and $h_{20}.[h_{20}]\in \overline{E}^2$ is not a permanent cycle: it survives to the page $4$ where $\overline{d}^4(h_{20}.[h_{20}]) =\overline{d}^4([h_{20}^2]) = [h_{11}^3]\neq 0$. Note the crossing differential hypothesis \autocite[(2.5)]{BK24} we need to check amounts to 
\begin{itemize}[$\ast$]
    \item All elements in $\overline{E}^4_{3,2,0}, \overline{E}^3_{3,2,-1}, \overline{E}^3_{4,3,0}$ are permanent cycles.
\end{itemize}
But in fact all of those vector spaces are zero. 
    
\end{prf}

\subsection{The self map on the stabilisation cofibre}\label{sec: self map on CGL/sigma}
So far we have collected various \say{$\E_\infty$-approximations} of $\CGL$: the bottom part
 $f_1:X_1\to \CGL$ of a minimal cell structure of $\CGL$ (\S\ref{sec: stability Hopf alg of CGL}), the {truncation map} $\CGL\to \cgl$ (\ref{map: truncation map}), and $\cgl\to a(1)$ constructed as cobar of  a map of  Hopf algebras (\S\ref{sec: Delta BGL -> a(1)}). Putting the set together we obtain an $\E_\infty$-map\stepcounter{thm}
 \[\tag{\thethm}\label{eq: X_1 -> a(1)}\Psi_{X_1}:X_1\to a(1)\]
 with $\Psi_{X_1*}(\sigma) = h_{10}$  (\cref{lem: map induced on SS}). So if we view $a(1)$ as a right $X_1$-module via $\Psi_{X_1}$, then there is an equivalence  $a(1)/h_{10}\simeq a(1)\otimes_{X_1}X_1/\sigma$ of left $a(1)$-modules. The goal of this section is to construct a left $\fil^{\E_\infty}X_1$-module endomorphism $\widetilde{\varphi_{X_1}}$ of tridegree $(12,8,-4)$ with the property that the basechange of the $X_1$-module map $\tau^{-1}\widetilde{\varphi_{X_1}}$ to $a(1)$ is equivalent to  acting by  $y_{12,8}\in \pi_{12,8}(a(1))$.

 The strategy is to take a well-chosen class in $\tensor*[^{\mathrm{BS}}]{E}{^\infty}(\fil^{\E_\infty }X_1/\widetilde{\sigma})$ and pick a filtered class $\widetilde{\omega_{12,8}}$ that is detected by it. We check  the filtered homotopy group at the tridegree of $\widetilde{\sigma}\widetilde{\omega_{12,8}}$ vanishes, so the left module map induced by any such $\widetilde{\omega_{12,8}}$ descend to the cofibre $\fil^{\E_\infty}X_1/\widetilde{\sigma}$. 

\subsubsection{Digression: Normal sequences}\label{sec: normal sequence}
The homotopy group of a free $\E_\infty$-algebra on  a graded $\F_2$-module $X$ is the free commutative algebra generated by those $Q^I(x)$, where $x$ ranges over a basis of $\pi_{*,*}(X)$ and $I$ ranges over all the admissible sequences for $x$ satisfying certain excess condition \autocite[181]{GKRW18a}; we say such $I$ is a  \textit{normal sequence} for $x$.  Explicitly, this means either
\begin{itemize}[$\ast$]
    \item $I=\varnothing$, and we interpret $Q^\varnothing$ as the identity operation; or 
    \item $I = (s)$ with $s>\abs{x}$; or
    \item $I=(s_1,\dots,s_r)$, $r\geqslant 2$, is admissible, meaning $s_{i-1}\leqslant 2s_i$ for $2\leqslant i\leqslant r$, and \[e(I)\defeq s_1-(s_2+\cdots + s_r)>\abs{x}.\]
\end{itemize} 
The normality condition is easier when expressed in terms of the \say{lower-indexed} Dyer-Lashof operations. They are related to the upper-indexed ones by the formula $Q^s(x) = Q_{s-\abs{x}}(x)$, where $\abs{x}$ denote the homological degree. So given a non-empty sequence $I$, one has $Q^I(x) =Q_{I'}(x)$, where 
\[I' = (s_1',\cdots ,s_r'),\quad s_i'  \defeq  s_i-\sum_{j>i}s_j - \abs{x}.\]
We say $I'$ is normal if $I$ is. We make the following observations.
\begin{obs}
    Suppose $I$ has length $r>1$. Let $I_{>1}$ denote the truncated sequence $(s_2,\dots,s_r)$. Then
    \begin{enumerate}[(i)]
        \item $(I_{>1})' = (s_2',\dots, s_r')$, so we can denote this sequence by $I'_{>1}$.
        \item If $I$ is normal, then so is $I_{>1}$. 
        \item If $I'$ is normal, then so is  $I'_{>1}$.
    \end{enumerate}
\end{obs}
\begin{lem}\label{lem: subsript normal seq}
    A non-empty $I =(s_1,\dots,s_r)$ is normal for $x$ iff  $I'=(s_1',\dots,s_r')$ is such that \stepcounter{thm}
\[\leqnomode \label{normality condition}\tag{\thethm}0<s_1'\leqslant \cdots  \leqslant s_r'.\]
\end{lem}
\begin{prf}
We induct on the length of $I$. When $I=  (s)$, then the condition for $I$ being normal is precisely $s' = s-\abs{x}>0$. Now suppose $I$ has length $r>1$. Suppose $I$ is normal, then by the observation above so is $I'_{>1}$. So we only need to check that $0<s_1'\leqslant s_2$.  The part $s_1' = e(I)-\abs{x} >0$ is the  excess condition for $I$, and that $s_1'-s_2' = s_1-2s_2\leqslant0$ is part the admissibility for $I$. Conversely, suppose (\ref{normality condition}) is satisfied. As saw above $s_1'>0$ is precisely the excess condition, so we only need to check admissibility.
By inductive hypothesis $I'_{>1} $ is normal, leaving us to verify $s_1\leqslant 2s_2$. Again, we just saw that $s_1-2s_2 =s_1'-s_2'$.
\end{prf}

\subsubsection{Construction of the filtered endomorphism $\widetilde{\varphi_{X_1}}$}\label{sec: self map on fil X} Throughout this section, let $\{E^r,d^r\}_r$ and $\{\overline{E}^r,\overline{d}^r\}_r$ denote the spectral sequence associated to $\fil^{\E_\infty}X_1$ and its left module $\fil^{\E_\infty}X_1/\widetilde{\sigma}$. Let $\beta$ be the $(3,2)$-cell in 
\[\gr \fil^{\E_\infty}X_1 = \E_\infty((-1)_!\bbS^{1,0}\{\sigma\}\oplus (-1)_!\bbS^{3,2}\{\beta\}).\]
For a class $x\in E^1_{*,*,*}$, we write $[x]$ for its image in the quotient algebra $\overline{E}^1_{*,*,*} = E^1_{*,*,*}/(\sigma)$. Then 
the majority of work goes into \cref{lem: E infty 12 8} and \cref{lem: E_infty 13 8} below. The first lemma describes the $E^\infty$-page of the Bockstein spectral sequence of $\fil^{\E_\infty}X_1/\widetilde{\sigma}$ in $(g,d,f)$-degree $(12,8,-4)$, in particular the class $[\beta^4]$ defines, potentially non-uniquely, a filtered homotopy class $\widetilde{\omega_{12,8}}$. 
The second lemma tells us any representative of $\widetilde{\omega_{12,8}}$ factors through the cofibre $\Sigma^{12,8,-4}\fil^{\E_\infty}X_1/\widetilde{\sigma}$, giving us the desired left module self-map $\widetilde{\varphi}_{X_1}$.
The proofs shall be delayed to the end of the section. 
\begin{lem}\label{lem: E infty 12 8}
    $\tensor*[^{\mathrm{BS}}]{\overline{E}}{^{\infty}_{12,8,-4,*}}$ is spanned by $[\beta^4]$ and $\tau^8[Q_{1,1}(\sigma)^2 Q_1(\sigma)^2]$. 
\end{lem}

\begin{defn}
    We let $\widetilde{\omega_{12,8}}\in \pi_{12,8,-4}(\fil^{\E_\infty}X_1/\widetilde{\sigma})$ be any class detected by the class $[\beta^4]$ in \cref{lem: E infty 12 8}.  
\end{defn}
The ambiguity of $\widetilde{\omega_{12,8}}$ is exactly the filtered homotopy class defined by  $\tau^8[Q_{1,1}(\sigma)^2 Q_1(\sigma)^2]$. The next lemma tells us that this ambiguity is immaterial.
 \begin{lem}\label{lem: E_infty 13 8}
 The group $\pi_{13,8,-5}(\fil^{\E_\infty}X_1/\widetilde{\sigma})$ is zero. In particular, $\widetilde{\sigma}.\widetilde{\omega_{12,8}}=0$ regardless of our choice of $\widetilde{\omega_{12,8}}$.
 \end{lem}
\begin{defn}\label{ref: def of base self map}
 Take any representative $\widetilde{\omega_{12,8}}: \Sigma^{12,8,-4}\fil^{\E_\infty}X_1\to \fil^{\E_\infty}X_1/\widetilde{\sigma}$. By \cref{lem: E_infty 13 8} above the composite $\widetilde{\omega_{12,8}}\circ \widetilde{\sigma}$ is null. Let
    \[\widetilde{\varphi_{X_1}}:\Sigma^{12,8,-4}\fil^{\E_\infty}X_1/\widetilde{\sigma}\to\fil^{\E_\infty}X_1/\widetilde{\sigma}\]be the left $\fil^{\E_\infty}X_1$-module endomorphism obtained via any such null-homotopy.
\end{defn}

\begin{prf}[of \cref{lem: E infty 12 8}]
We start by noting the the spectral sequences $E^r$ and $\overline{E}^r$ are very similar to the cell-attachment spectral sequence of the algebra $\overline{\mathbf{A}}$ and that for $\overline{\mathbf{A}}/\sigma$ in Section 6.2 of \autocite[]{GKRW18b}. Since the differentials are analysed extensively there, we will assert differentials in our settings without the details. The class called $\tau$ there corresponds to our $\beta$ here.  The result is  a consequence of 
\begin{claim}
$\overline{E}^1_{12,8,*}$ is spanned by the classes in \cref{tab:basis foe E^1_{12,8,*}}.
\begin{table}[h]
    \centering
    \begin{tabular}{lll}
    name &  class  & multiplicative filtration \\
    \hline 
     & $[\beta^4]$  &$-4$\\
     $\chi_2$ & $[Q_1(\sigma)^2 \beta ^2Q_{2}(\sigma)]$ & $-6$\\
     $\chi_2'$ &  $[Q_1(\sigma)^2Q_{1,1}(\sigma)^2]$ & $-12$ \\
     $\chi_3$   & $[Q_1(\sigma)^3 Q_1(\beta)]$  &$-8$\\
     $\chi_3'$ & $[Q_1(\sigma)^3Q_{1,1}(\sigma)Q_2(\sigma)]$  & $-12$\\
     $\chi_4$  & $[Q_2(\sigma)^2 Q_1(\sigma)^4]$  & $-12$\\
      $\chi_5$ & $[Q_1(\sigma)^5Q_3(\sigma)]$   & $-12$
\end{tabular}
    \caption{Basis for $\overline{E}^1_{12,8,*}$.}
    \label{tab:basis foe E^1_{12,8,*}}
\end{table}
    \end{claim}
The claim tells us $\tensor*[^{\mathrm{BS}}]{\overline{E}}{^1_{12,8,-4,*}}$ is spanned by $[\beta^4]$, $\tau^2\chi_2$, 
 $\tau^8\chi'_2$, $\tau^4\chi_3$, $\tau^8\chi_3'$, $\tau^8\chi_4$, $\tau^8\chi_5$. The reslut then follows from the case analysis below.
 \begin{itemize}[$\ast$]
 \item  In $\tau$-filtration $-2$, the class $\tau^2 \chi_2$ supports a non-trivial differential $\tensor*[^{\mathrm{BS}}]{\overline{d}}{_2}(\tau^2\chi_2)  =\tau^6[ Q_1(\sigma)^5 Q_2(\sigma)]$. 
     \item In $\tau$-filtration $-4$, there is only $\tau^4\chi_3$, which is the target of $ \tensor*[^{\mathrm{BS}}]{\overline{d}}{^4}(Q_1(\sigma)Q_2(\sigma)^2[\beta^2])= \tau^4\chi_3$.
     \item In $\tau$-filtration $-8$, $\tau^8\chi_2$ is a  permanent cycle and all of $\tau^8\chi_3',\tau^8\chi_4,\tau^8\chi_5$ are target of differentials
      \begin{align*}
     \tensor*[^{\mathrm{BS}}]{\overline{d}}{^4}(\tau^4Q_{1,1}Q_2(\sigma).[\beta^2]) &= \tau^8\chi_3'\\
       \tensor*[^{\mathrm{BS}}]{\overline{d}}{^4}(\tau^4Q_2(\sigma)^2Q_1(\sigma).[\beta^2])& = \tau^8\chi_4\\
         \tensor*[^{\mathrm{BS}}]{\overline{d}}{^4}(\tau^4 Q_1(\sigma)^2Q_3(\sigma).[\beta^2]) & = \tau^8\chi_5
 \end{align*}
 So the $\tensor*[^{\mathrm{BS}}]{\overline{E}}{^{\infty}_{12,8,-4,-8}}$ is spanned by  $\tau^8\chi_2$. 
 \end{itemize}

\begin{prfoc}Let $S$ be the set of $Q_I(x)$ with $x\in \{\sigma, \beta\} $ and $I$ ranges over the non-empty normal sequences for $x$. Let $M(S)$ be the set of 
 monomials generated by $S$ and identify $\overline{E}^1_{*,*,*}$ with the vector space $\F_2\{M(S)\}$ spanned by $M(S)$ via the bijection \[\F_2\{M(S)\}\hookrightarrow \F_2[\{Q^I(x)\}_{I,x}] = E^1_{*,*,*}\to \overline{E}^1_{*,*,*}\]
 which, in symbols, sends $Q_I(x)$  to $[Q_I(x)]$. For each $r\geqslant 0$, define
 \[A(S)^{(r)}\defeq \{x\in M(S)_{12,8,*}\mid x\text{ is divisible by }Q_1(\sigma)^{r}\text{ but not by }Q_1(\sigma)^{r+1}\}.\]
We now show  that $A(S)^{(0)} = \{\beta^4\}$, $A(S)^{(i)} =\{\chi_i,\chi_i'\}$ for $i =2,3$, $A(S)^{(i)} = \{\chi_i\}$ for $i = 4,5$,  and all other $A(S)^{(r)}$ are empty. Firstly, $A(S)^{(r)} = \varnothing$ for $r \geqslant 6$ simply the grading of $Q_1(\sigma)^r$ is $2r$, and all element of $S$ have positive grading. To compute  other $A(S)^{(r)}$, it is convenient to organise $S$ by the length of $I$: Let $L_rS$ denote the set of element o the form $Q_I(\sigma), Q_I(\beta)$, with $I$ having length $r$, then one can write down  the following table of  $L_rS\cap \{\text{ elements of homological degree} \leqslant 8\}$ using  \cref{lem: subsript normal seq}.
\begin{table}[h]
    \centering
    \begin{tabular}{l|ll l l }
    $r=0$  & $(3,2)\;\;\beta$  & \\
    $r=1$ & $(2,i)\;\;Q_i(\sigma)$, $1\leqslant i\leqslant 8$ & $(6,4+j)\;\;Q_j(\beta)$, $1\leqslant j\leqslant 4$\\
    $r=2$  & $(4,3)\;\;Q_{1,1}(\sigma)$ &  $(4,5)\;\;Q_{1,2}(\sigma)$ & $(4,7)\;\;Q_{1,3}(\sigma)$ & $(4,6)\;Q_{2,2}(\sigma)$ \\
    & $(4,8)\; Q_{2,3}(\sigma)$ & $(4,9)\; Q_{1,4}(\sigma)$ & $(4,9)\; Q_{3,3}(\sigma)$
\end{tabular}
    \caption{$L_rS\cap \{\text{ elements of homological degree} \leqslant 8\}$. The classes are labelled by its  grading and homological degree.}
    \label{tab: LrS}
\end{table}
If we define the  \textit{slope} of an element $x\in M(S)$ of tridegree $(g,d,f)$  to be $\frac{d}{g}$. Then 
\[\mathrm{slope}(Q_{(s_1,s_2,\dots, s_n)}(x)) = \frac{s_1 + 2s_2 + \cdots + 2^{n-1}s_n + 2^n d(x)}{2^n g(x)}.\]
In particular, when $r>0$, $\min \mathrm{slope}(L_rS) = \frac{2^r-1}{2^r}$, strictly achieved at $Q_1^{\circ r}(\sigma)$. For an interval $I\subseteq [0,1]$, let $S'_I$  denote the set of element in $S'\defeq S\setminus \{Q_1(\sigma)\}$ of slope $\lambda\in I$. Let $M(S')$ be the set of non-zero monomials on $S'$. Since $A(S)^{(r)} \subseteq Q_1(\sigma)^r.\F_2\{M(S')\}$, our strategy is to  compute $A(S)^{(r)}$ by analysing slope of elements in $S'_{[0,1]}$. A useful consequence of the above discussion of minimum slope is
\[S'_{[0,\frac{2}{3})}= \varnothing, \quad \quad S'_{[0,\frac{2}{3}]} = \{\beta\},\quad \quad S'_{(\frac{2}{3}, \frac{3}{4}]} = \{Q_{1,1}(\sigma)\},\]
from which it follows that $A(S)^{(0)} =\{\beta^4\}$.
\begin{itemize}[$\ast$]
    \item If $x\in A(S)^{(1)}$, then $x = Q_1(\sigma)x'$ for some $x'\in M(S')_{10,7,*}$. Since $\frac{7}{10}< \frac{3}{4}$ and $S'_{[0,\frac{2}{3})}= \varnothing$, $A(S)^{(1)} = \varnothing$. 
    \item 
If $x\in A(S)^{(2)}$,  then $x  =Q_1(\sigma)^2x'$ for  some $x'\in M(S')_{8,6,*}$. If $y$ is $\beta$-divisible,  then as $\beta$ is the only element in $S'$ of odd grading, $y$ must be $\beta^2$-divisible, so $y =\beta^2 z $ with $z \in M(S')_{2,2}(S') = \{Q_2(\sigma)\}$.  If $y$ is not $\beta$-divisible, then $\frac{6}{8} =\frac{3}{4}$ tells us $y=Q_{1,1}(\sigma)^2$.
\item  If $x\in A(S)^{(3)}$, then $x = Q_1(\sigma)^3y$ for some $y\in M(S')_{6,5,*}$. As before, if $y$ is $\beta$-divisible, then $y$ is $\beta^2$-divisible, but this cannot happen since $M(S')_{0,1,*} =\varnothing$. Now, $\min\mathrm{slope}(L_rS)$ exceeds $\frac{5}{6}$ when $r\geqslant 3$, so $y$ has to contain  at least a factor in $\ (L_1S\cup L_2S)\cap S'$. If $y$ as no factor in $L_2S$, then we must have $y=Q_1(\beta)$, for otherwise $y$ is a  monomials of grading $6$ on $Q_i(\sigma)$, $i\geqslant 2$, and all monomial has homological degree $>5$.  If $y$ has a factor in $L_2S\cap S'$, then this can only be $Q_{1,1}(\sigma)$ ($M(S')_{2,1} = \varnothing$), so $y = Q_{1,1}(\sigma) Q_2(\sigma)$.
\item  If $x \in A(S)^{(4)}$, then $x= Q_1(\sigma)^4 y$, $y\in M(S')_{4,4} = \{Q_2(\sigma)^2\}$.
\item Finally, if $x\in A(S)^{(5)}$ then $x =Q_1(\sigma)^5 y$ for some $y\in M(S')_{2,3} = \{Q_3(\sigma)\}$. 
\end{itemize}
\end{prfoc}
\end{prf}

 \begin{prf}[of \cref{lem: E_infty 13 8}]
 We first show that $\overline{E}^1_{13,8,*}$ is spanned by the classses in \cref{tab:E^1_{13,8,*}}.
 \begin{table}[h]
     \centering
\begin{tabular}{lll}
    name &  class  & multiplicative filtration \\
    \hline 
     $\varpi_2$& $[Q^1(\sigma)^2\beta^3]$  &$-7$\\
     $\varpi_3$ & $[Q_1(\sigma)^3Q_{1,1}(\sigma)\beta  ]$ & $-11$\\
     $\varpi_4$ &  $[Q_1(\sigma)^4Q_{2}(\sigma)\beta]$ & $-11$ 
\end{tabular}
     \caption{Basis of $\overline{E}^1_{13,8,*}$.}
     \label{tab:E^1_{13,8,*}}
 \end{table}
 We use a similar argument as in the proof of   \cref{lem: E infty 12 8}; We also adopt the same set of notations used there, except this time we set $A(S)^{(r)}$ to consist of exactly $Q_1(\sigma)^r$-divisible elements in $M(S)$ of filtration $(13,8,*)$.
\begin{claim}
    $A(S)^{(r)} = \{\varpi_r\}$ for $r =2,3,4$, and all other $A(S)^{(r)}$ are empty. 
\end{claim}
\begin{prfoc}\begin{itemize}[$\ast$]
    \item Suppose $x  = \beta ^i y\in A(S)^{(0)}$ for some $i\geqslant 0$, then $\mathrm{slope}(y) = \frac{8-2i}{13-3i}$, which is $<\frac{3}{4}$ for all $i<7$. So by iteratively applying the fact $S'_{[0,\frac{3}{4})} = S_{\frac{2}{3}}' = \{\beta\}$,  we see any $x\in A(S)^{(0)}$ must be $\beta^7$-divisible, which cannot happen because the grading of  $\beta^8$ is too large.
    \item Suppose $x\in A(S)^{(1)}$, then $ x= Q_1(\sigma )y$ for some $y\in M(S')_{11,7,*}$. As above, the inequality $\frac{11-2i}{7-3i}<\frac{3}{4}$ for $i<5$ shows that $y$ must be $\beta^5$-divisible, which cannot happen.
    \item Suppose $x\in A(S)^{(2)}$, so $x = Q_1(\sigma)^2 y$ for some $y\in M(S')_{9,6,*}$. Since $\frac{6}{9}= \frac{2}{3}= \min \mathrm{slope}(S')$, one can only have $y = \beta ^3$.
    \item If $x\in A(S)^{(3)}$ then $x = Q_1(\sigma)^3y$ for some $y\in M(S')_{7,5,*}$. Since $\beta$ is the only element of $S'$ with odd grading, $y = \beta z$ for some $z\in M(S')_{4,3,*}$. But $ M(S')_{4,3,*}=\{Q_{1,1}(\sigma)\}$: element in $S'$ of length $j$ has grading $\geqslant 2^j$, and we wrote down the list of element of length $\leqslant 2$ with homological grading $\leqslant 8$ in  \cref{tab: LrS}. 
    \item If $x\in A(S)^{(4)}$, then $x = Q_1(\sigma)^4 y$ for some $y\in M(S')_{5,4,*}$. Again, consideration on parity of the grading shows $y =\beta z$ for some $z\in M(S')_{2,2,*}$, and  $M(S')_{2,2,*} =\{Q_2(\sigma)\}$.
    \item $A(S)^{(5)} = Q_1(\sigma)^5\cdot M(S')_{3,3,*}$ and $A(S)^{(6)} = Q_1(\sigma)^6 \cdot M(S')_{1,2,*}$. But both $ M(S')_{3,3,*}$ and $M(S')_{1,2,*}$ are empty.
    \item The grading of $Q_1(\sigma)^r$ exceeds $13$ when $r\geqslant 7$, so $A(S)^{(r)} =\varnothing$ for $r\geqslant  7$.
\end{itemize}
\end{prfoc}
By the claim, $\tensor*[^{\mathrm{BS}}]{\overline{E}}{^1}_{13,8,-5,*}$  is spanned by $\tau^2\varpi_2$, $\tau^6\varpi_3$, $\tau^6\varpi_4$. The first class survives to page 4 and supports a non-trivial differential $\tensor*[^{\mathrm{BS}}]{\overline{d}}{^4}(\tau^2\varpi_1) = \tau^6[\beta Q^1(\sigma)^5]$ there. The other two are permanent cycles, but they are also target of differentials
\begin{align*}
    \tensor*[^{\mathrm{BS}}]{\overline{d}}{^4}(\tau^2Q_{1,1}(\sigma)[\beta^3]) &= \tau^6\varpi_3\\
    \tensor*[^{\mathrm{BS}}]{\overline{d}}{^4}(\tau^2 Q_1(\sigma) Q_2(\sigma)[\beta^3])&=\tau^6\varpi_4.
\end{align*} So $\tensor*[^{\mathrm{BS}}]{\overline{E}}{^\infty}_{13,8,-5,*}$ is zero, and the lemma follows.
 \end{prf}

\subsubsection{Basechange to $a(1)$.}\label{sec: basechange to a(1)} 

Our aim is to describe the basechange of the endomorphism $\varphi_{X_1}\defeq \tau^{-1}\widetilde{\varphi_{X_1}}$ along the $\E_\infty$-map $\Psi_{X_1}:X_1\to a(1)$ of (\ref{eq: X_1 -> a(1)}), see \cref{prop: basechang as expected} below. To do so, we refine $\Psi_{X_1}$ to a filtered $\E_\infty$-map
\[\widetilde{\Theta}\defeq u_{a(1)}\circ \fil^{\E_\infty}(\Psi_{X_1}):\fil^{\E_\infty}X_1\to \fil^{\mathrm{aug}}a(1),\]
where $u_{a(1)}:\fil^{\E_\infty}a(1)\to \fil^{\mathrm{aug}}a(1)$ is obtained from a universal property of the canonical multiplicative filtration; refer the reader to \autocite[\S8.4]{RW25} for  details.

\begin{lem}\label{lem: map induced on SS}
    Let $\{E^r(\widetilde{\Theta})\}_r$ be the map of spectral sequences induced by $\widetilde{\Theta}$. Then
    \begin{enumerate}[(i)]
        \item $E^1(\widetilde{\Theta}) (\sigma) = h_{10}\in E^1(\fil^{\mathrm{aug}}a(1))$, 
        \item $E^1(\widetilde{\Theta})(\beta)  = h_{20}\in  E^1(\fil^{\mathrm{aug}}a(1))$
    \end{enumerate}
 where the $\sigma$, $\beta\in E^1(\fil^{\E_\infty}X_1)$ are as in \S\ref{sec: self map on fil X}. In particular, $\Psi_{X_1} \cong \tau^{-1}\widetilde{\Theta}$ sends $\sigma\in \pi_{1,0}(X_1)$ to $h_{10}\in \pi_{1,0}(a(1))$, and hence $a(1)\otimes_{X_1}X_1/\sigma\simeq a(1)/h_{10}$ as left $a(1)$-modules.
\end{lem}

\begin{prf} For part (i), we first remind ourself that the official definition definition of $\Psi_{X_1}$ is the composite
\[X_1\xrightarrow[]{f_1}\CGL\xrightarrow[]{t_{\CGL}}\cgl\xrightarrow[]{\Cobar(\varphi^\vee)}a(1),\]
where $f_1$ was the inclusion of bottom cells, $t_{\CGL}$ is the truncation $\E_\infty$-map (\ref{map: truncation map}), and  $\varphi^\vee$ was constructed in \S\ref{sec: cgl to a(1)}. Consider the Hopf algebra map 
\[\varphi':\Delta_{X_1} = \F_2[\overline{\sigma},\delta]\xrightarrow[]{} \Delta_{\CGL}\xrightarrow[]{\varphi^\vee} A(1)_*\]
where first map sends the generators to the element in $\Delta_{\CGL}$ with the same name. Then
by naturality of the truncation map $t_{(-)}$, $\Psi_{X_1}$ is equivalent  to
\[X_1\xrightarrow[]{t_{X_1}}x_1\xrightarrow[]{\Cobar(\varphi')}a(1)\]
 and thus it suffices to show
\begin{enumerate}[(a)]
    \item $u_{x_1}\fil^{\E_\infty}(t_{X_1}):\fil^{\E_\infty}X_1\to \fil^{\mathrm{aug}}x_1$ sends $\sigma\in E^1(\fil^{\E_\infty}(X_1))$ to $\underline{\overline{\sigma}}\in E^1(\fil^{\mathrm{aug}}x_1)$. 
    \item $E^1(\fil^{\mathrm{aug}}x_1)\to E^1(\fil^{\mathrm{aug}}a(1))$ induced by $\varphi'$ sends $\underline{\overline{\sigma}}$ to $h_{10}$.
\end{enumerate}
Here $\underline{\overline{\sigma}}$ is defined by the cycle in the cobar complex given by the primitive element $\overline{\sigma}\in \Delta_{X_1}$. In particular,  (b) is just a consequence of $\varphi'(\overline{\sigma}) = \xi_1$. By naturality, we can reduce part (a) to 
\begin{itemize}
    \item[(a')] Let $R =\E_{\infty}(\bbS^{1,0}\{\sigma\})$, then $E^1(u_r\circ \fil^{\E_{\infty}}(t_R))(\sigma) = \underline{\overline{\sigma}}$.
\end{itemize}
Since both $\sigma$ and $\overline{\underline{\sigma}}$ are permanent cycles that corresponds to the generator of their respective $\pi_{1,0} \cong \F_2$, it suffices to show that $\tau^{-1}(u_r\circ \fil^{\E_\infty}(t_R))$
is an isomorphism on $\pi_{1,0}$. Since $\tau^{-1}u_r$ is an equivalence, it remains to consider $\tau^{-1}\fil^{\E_\infty}t_R\simeq t_R$. By  \autocite[Lemma A.11]{Bur22}, $t_R$ is an equivalence through grading $1$ iff $\Ba(t_R)$ is. As $R$ is connected, both $\Ba R$ and $\tau^{\mathrm{diag}}_{\leqslant0}\Ba R$ are simply connected, so $\Ba \Cobar \simeq \id$  at them. As a consequence,
 $\Ba(t_R)$ is equivalent to the the truncation map $\Ba R\to \tau^{\mathrm{diag}}_{\leqslant0}\Ba R$ with respect to the diagonal $t$-structure, which is indeed an equivalence through grading $1$. 

For  (ii), simply note that $h_{20}$ is the unique class in $E^1_{3,2,-1}(\fil^{\mathrm{aug}}a(1))$ that survives to $E^2$ and whose image under $d^2$ is $h_{10}h_{11}$. On the other hand, $d^1(E^1(\widetilde{\Theta})(\beta)) = E^1(\widetilde{\Theta})(d^1\beta)=0$ and  
\[d^2(E^2(\widetilde{\Theta})(\beta))=E^2(\widetilde{\Theta})(d^2\beta) =E^2(\widetilde{\Theta})(\sigma Q_1(\sigma)) =h_{10}h_{11}. \]
The last equality uses (i) and that $\widetilde{\Theta}$ is $\E_\infty$.
 \end{prf}

 \begin{prop}\label{prop: basechang as expected}
     Under the identification above, $a(1)\otimes_{X_1}{\varphi_{X_1}}$ is equivalent to multiplication by $y_{12,8}$ as a left module endomorphism.
 \end{prop}
\begin{prf}
    Write $[-,-]$ for the homotopy class of left $a(1)$-module maps. 
    The cofibre sequence defining $a(1)/h_{10}$ induces a long exact sequence with a segment of the form
    \[\pi_{13,9}(a(1)/h_{10})\to [\Sigma^{12,8}a(1)/h_{10}, a(1)/h_{10}]\to \pi_{12,8}(a(1)/h_{10})\to \pi_{13,8}(a(1)/h_{10}).\]
  By \cref{thm: homotopy of a(1)/h10} the first and last term vanishes, giving us  \[[\Sigma^{12,8}a(1)/h_{10}, a(1)/h_{10}]\cong \pi_{12,8}(a(1)/h_{10})\cong \F_2.\] Now it is enough to show $a(1)\otimes_{X_1}{\varphi_{X_1}}$ is not null. Well, the image of $(a(1)\otimes_{X_1}{\varphi_{X_1}})_*(z_{0,0})$  is detected by the non-zero element 
 $E^\infty(\widetilde{\Theta})(\beta^4) = h_{20}^4$.
\end{prf}

\section{Applications of the detector}\label{sec: applications of the detector}
In \S\ref{sec: detecting periodic} we state and prove the precise version of \cref{intro thm: rel} and \cref{intro thm: abs} in the introduction. In \S\ref{sec: snake lemma} we show that,  up to filtration, the relative classes $\alpha_{i1},\alpha_{i2}$ do not contribute new information to the absolute homology groups.
The situation can be summarised by \cref{fig: flash in CGL} below. \textit{As in the previous section, $\CGL\defeq \CGL(\F_2)$.}

\begin{figure}[h]
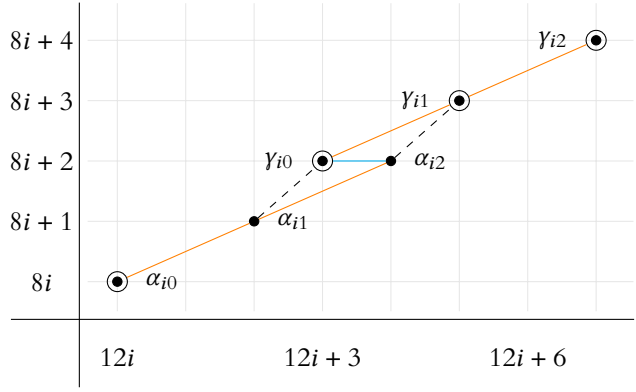

\centering
{
\begin{sseqpage}[ classes = fill,  xscale = 0.9,  yscale = 0.8, grid = go, grid color =gray!20, x range = {0}{7}, y range = {0}{4}, x tick handler = {
\ifnum#1 = 0\relax
12i
\else
\ifnum#1 = 3\relax
% \vphantom is fragile so we \protect it
\protect\vphantom{2}12i+#1
\else
\ifnum#1 = 6\relax
% \vphantom is fragile so we \protect it
\protect\vphantom{2}12i+#1
\else
\relax 
\fi
\fi
\fi
}, y tick handler = {
\ifnum#1 = 0\relax
8i
\else
\relax 8i+#1\fi
}]
\mygroup(0,0)
\mygroup(3,2)
\classoptions["\alpha_{i0}" right,circlen = 2,fill](0,0)
\classoptions["\alpha_{i1}" right](2,1)
\classoptions["\alpha_{i2}" right](4,2)
\classoptions["\gamma_{i0}" left,circlen = 2,fill](3,2)
\classoptions["\gamma_{i1}" left,circlen = 2,fill](5,3)
\classoptions["\gamma_{i2}" left,circlen = 2,fill](7,4)
\structline[cyan](3,2)(4,2)
\structline[dashed](3,2)(2,1)
\structline[dashed](5,3)(4,2)
\end{sseqpage}}
    \caption{The non-zero classes $\alpha_{ij},\gamma_{ij}$ ($i\geqslant 0$) of \cref{prop: BGL/sigma classes}; the circled classes survives under the boundary map $\partial^{\CGL}_*$; the dashed line indicates $q^{\CGL}_*\partial^{\CGL}_*(\gamma_{ij})\equiv \alpha_{i(j+1)}$ in $E^\infty(\fil^{\E_\infty}\CGL/\widetilde{\sigma})$. The orange / blue line indicates a $Q^1(\sigma)$ / $\sigma$-extension.}
    \label{fig: flash in CGL}
\end{figure}

\subsection{Detecting periodic families in $\CGL(\F_2)/\sigma$}\label{sec: detecting periodic}

Let \[\Psi:\CGL\to a(1)\] be the composite of $t_{\CGL}:\CGL\to \cgl$  (\ref{map: truncation map}) with the detection map $\cgl\to a(1)$  that we built in \S\ref{sec: cgl to a(1)}. As $\Psi f_1=\Psi_{X_1}$, so by \cref{lem: map induced on SS} $\Psi_*(\sigma) = h_{10}$ and hence $\Psi$ induces a map of cofibre sequences

\[\begin{tikzcd}
        \Sigma^{1,0}\CGL \ar["\sigma"]{r}\ar["\Sigma^{1,0}\Psi"]{d}& \CGL\ar["q^{\CGL}"]{r} \ar["\Psi"]{d}& \CGL/\sigma  \ar["\partial^{\CGL}"]{r} \ar["\overline{\Psi}"]{d}& \Sigma^{1,1}\CGL\ar["\Sigma^{1,1}\Psi"]{d}\\
        \Sigma^{1,0}a(1) \ar["h_{10}"]{r} & a(1)\ar["q^{a(1)}"]{r} &a(1)/h_{10} \ar["\partial^{a(1)}"]{r} & \Sigma^{1,1}a(1)
    \end{tikzcd}\]
An inspection of the long exacts sequence  shows that the boundary map \[\pi_{3,2}(\CGL/\sigma)\to \pi_{2,1}(\CGL) = \F_2\{Q^1(\sigma)\}\] is an isomorphism; let $\omega_{3,2}\in \pi_{3,2}(\CGL/\sigma)$ be the class corresponding to $Q^1(\sigma)$. We write
 $\omega_{0,0}\defeq q^{\CGL}_*(1)\in \pi_{0,0}(\CGL/\sigma)$. Lastly, let $\varphi:\Sigma^{12,8}\CGL/\sigma\to \CGL/\sigma$ denote the left $\CGL$-module endomorphism $\CGL\otimes_{X_1}{\varphi_{X_1}}$.

\begin{prop}\label{prop: BGL/sigma classes}
    For every $i\geqslant 0$, $j\in \{0,1,2\}$,  the classes\stepcounter{thm}
    \begin{align*}
        \alpha_{ij}&\defeq \varphi^i_*(Q^1(\sigma)^j.\omega_{0,0})\in \pi_{12i+2j, 8i+j}(\CGL/\sigma)\\
        \gamma_{ij}&\defeq \varphi^i_*(Q^1(\sigma)^j.\omega_{3,2})\in \pi_{12i+2j+3, 8i+j+2}(\CGL/\sigma)\leqnomode\tag{\thethm}\label{classes: in cofibre}
    \end{align*}
    are such that
    \begin{align*}
        \Psi_*(\alpha_{ij})  &= y_{12,8}^ih_{11}^j.z_{0,0}\\
        \Psi_*(\gamma_{ij}) &= y_{12,8}^ih_{11}^j.z_{3,2}
    \end{align*}
   In particular, $\alpha_{ij}, \gamma_{ij}$ are non-zero. 
\end{prop}

\begin{prf}
 Using the $\pi_{*,*}(a(1))$-module structure of $\pi_{*,*}(a(1)/h_{10})$ given in \cref{thm: homotopy of a(1)/h10} and the fact that $\Psi_*(\sigma)= h_{10}$, it suffices to check $\Psi_*(\alpha_{00}) =z_{0,0}$ and $\Psi_*(\gamma_{00}) = z_{3,2}$. This can be done by using the naturality of long exact sequence. \end{prf}

\begin{rem}\label{rem: stability on X/sigma} By construction $f_{1*}:\pi_{g,d}(X_1)\to\pi_{g,d}(\CGL)$ is an isomorphism for $g\leqslant 3$, $d\leqslant 2$, so $\omega_{3,2}\in \pi_{3,2}(\CGL/\sigma)$ corresponds uniquely to a class in $\pi_{3,2}(X_1/\sigma)$ which we will also denote by $\omega_{3,2}$. And of course, the image of $1$ in $\pi_{0,0}(X_1/\sigma)$ is sent to $\omega_{0,0}\in \pi_{0,0}(\CGL/\sigma)$ under $f_1$, so we can also call this class $\omega_{0,0}$. Now formula (\ref{classes: in cofibre}) with ${\varphi_{X_1}}$ in place of $\varphi$ defines classes in $X_1/\sigma$ that are mapped to $\alpha_{ij},\gamma_{ij}$ under $f_1$. So \cref{prop: BGL/sigma classes} also tells us that $X_1/\sigma$ does not have a slope $>\frac{2}{3}$ vanishing line. Applying \autocite[Theorem 8.1 (ii)]{RW25} to $t_{X_1}:X_1\to x_1$ we see that $x_1/\sigma$ also does not have a slope $>\frac{2}{3}$ vanishing line. 
\end{rem}

We now show that the top part and the tips of the lightning flash survives under the boundary map. 

\begin{thm}\label{thm: absolute clases}
Let $\alpha_{ij},\gamma_{ij}\in \pi_{*,*}(\CGL/\sigma)$, $i\geqslant 0, j\in\{0,1,2\}$ be the classes in \cref{prop: BGL/sigma classes}. 
\begin{enumerate}[(i)]
    \item For every $i\geqslant0$, $j\in \{0,1\}$, the class \[{u}_{ij}\defeq \partial^{\CGL}_*( \gamma_{ij})\in \pi_{12i+2j+2, 8i+j+1}(\CGL)\]
is non-zero.
    \item The classes
    \begin{align*}
        {s}_i&\defeq \partial^{\CGL}_*(\alpha_{i0})\in \pi_{12i-1,8i-1}(\CGL)&,i>0\\
        u_{i2}&\defeq \partial^{\CGL}_*(\gamma_{i2})\in \pi_{12i+6,8i+3}(\CGL) &, i\geqslant 0
    \end{align*}
    are non-zero.
\end{enumerate}
\end{thm}

\begin{prf}
     For (i), it suffices to check $\overline{\Psi}_*(u_{ij})$, which by naturality equals  $\partial^{a(1)}(y^i_{12,8}h_{11}^i.z_{3,2})$, are non-zero. They are indeed non-zero because  in the long exact sequence for $a(1)$, the vector spaces to both side of the boundary map vanishes, and the target of that boundary map is non-zero. 

For (ii), suppose $u_{i2}=0$ for a contradiction. Then $\gamma_{i2} =q^{\CGL}_*(\gamma_i')$ for some class $\gamma_i'\in \pi_{12i+7,8i+4}(\CGL)$. $\Psi_*(\gamma_i')$ is non-zero because $q_{a(1)*}\Psi_*(\gamma_i')=\overline{\Psi}_*(\gamma_{i2})\neq0$ by \cref{prop: BGL/sigma classes}, so $\Psi_*(\gamma_i')$ must equal the only non-zero class in that bidegree, namely $y_{12,8}^iy_{7,4}$ (c.f. (\ref{eq: alg of a(1)}), \cref{fig: homotopy groups of a(1)}) which is $h_{10}$-torsion free. But every class in $\pi_{*,>0}(\CGL)$ is $\sigma$-torsion\footnote{The stable homology $H_d(\colim_n(\GL_n(\F_2));\F_2)$ vanishes for $d>0$.}. Similarly, if $s_i=0$, then $\alpha_{i0} = q^{\CGL}_*(\alpha_i')$ for some class $\alpha_i'$. As before one must have $\Psi_*(\alpha_i') =y_{12,8}^i$, which is $h_{10}$-torsion free.\end{prf}

\subsection{A snake lemma for filtered objects}\label{sec: snake lemma}
Our main result for this section is \cref{cor: collapsing lightning flash}, which says that for $\alpha_{i1}$ and $\alpha_{i2}$ are related to $\gamma_{i0}$ and respectively $\gamma_{i1}$ by a \say{$\sigma$-Bockstein} {up to certain filtration}. Due to the possibly non-trivial choice involved in  building  the endomorphism $\varphi_{\CGL}$, this is the best we can hope for. \cref{cor: collapsing lightning flash} itself does not shed new light on $H_*(\GL_n(\F_2);\F_2)$; rather we hope that the method of proof will be interesting in its own right.

 For $X\in \Fil(\fancycat{D}(\F_2)^{\Z})$, we write 
\[\delta^\infty_1:X/\tau\to \Sigma^{0,1,1}X\] for the boundary map of the cofibre sequence $\Sigma^{0,0,1}X\xrightarrow[]{\tau} X\xrightarrow[]{\mathrm{mod}\tau}X/\tau$. We learnt this definition from \autocite[Section 2.3]{CDN24}. Morally speaking, $\delta_1^\infty$ encodes all the differentials in $\{E^r(X)\}_r$.

Given a cofibre sequence  $X\xrightarrow[]{f}Y\xrightarrow[]{g}Z$ in $\Fil(\fancycat{D}(\F_2)^{\mathrm{gr}})$ with boundary map $Z\xrightarrow[]{\partial}\Sigma^{0,1,0}X$, we propose a way to to compute the boundary map on certain subset of $\pi_{*,*,*}(Z)$. 
Consider the following commutative diagram \stepcounter{thm}
\[\leqnomode\tag{\thethm}\label{eq: big diagram}\begin{tikzcd}
    X\arrow[r, "f"]\arrow[d,"\mathrm{mod}\tau"]  & Y\arrow[r,"g"] \arrow[d,"\mathrm{mod}\tau"] & Z\arrow[d,"\mathrm{mod}\tau"] \arrow[r,"\partial"]& \Sigma^{0,1,0}X\\
    X/\tau \arrow[r,"f/\tau"]\arrow[d,"\delta_1^\infty"]& Y/\tau  \arrow[r,"g/\tau"]\arrow[d,"\delta_1^\infty"]& Z/\tau \arrow[d,"\delta_1^\infty"]\\
    \Sigma^{0,1,1}X \arrow[r, "\Sigma^{0,1,1}f"]\arrow[d,"\Sigma^{0,1,0}\tau"]&\Sigma^{0,1,1}Y \arrow[r, "\Sigma^{0,1,1}g"]& \Sigma^{0,1,1}Z\\
    \Sigma^{0,1,0}X
\end{tikzcd}\]
whose rows and columns are cofibre sequences. Suppose $\widetilde{z}\in \pi_{g,d,f}(Z)$ has \[(\partial/\tau \circ \mathrm{mod}\tau)_*(z) =0,\]
then we can pick a representative $\widetilde{z}:\bbS^{g,d,f}\to Z$ \textit{and} a null-homotopy of $\partial/\tau \circ \mathrm{mod}\tau\circ \widetilde{z}$, giving rise to a $y:\bbS^{g,d,f}\to Y/\tau$. The bottom right homotopy square together with the right vertical cofibre sequence provides a null-homotopy of $\Sigma^{0,1,1}g\circ \delta_1^{\infty}\circ y$, giving us a map $\widetilde{x}:\bbS^{g,d,f}\to \Sigma^{0,1,1}X$. 

\begin{lem}
    In the situation above, \[\partial_*(z) = \tau \widetilde{x}\in \pi_{g,d,f}(\Sigma^{0,1,0}X)\] regardless of those choices we made. In particular, we can specify the map $y$ by a class $y\in \pi_{g,d,f}(Y/\tau)$ such that $(g/\tau)_*(y) = (\mathrm{mod}\tau)_*(z)$.
\end{lem}
\begin{prf}
The input data provides a map of diagrams from 
\[\begin{tikzcd}
&0 \arrow[d] \arrow[r]& \bbS^{g,d,f} \arrow[r,"\partial^{\mathrm{hor}}"]\arrow[d,equal] & \bbS^{g,d,f}\\
     0 \arrow[r]\arrow[d]& \bbS^{g,d,f}\arrow[r, equal]  \arrow[d,equal]& \bbS^{g,d,f}\arrow[d] \\
     \bbS^{g,d,f} \arrow[r,equal]\arrow[d,"\partial^{\mathrm{vert}}"]  & \bbS^{g,d,f} \arrow[r] &0 \\
  \bbS^{g,d,f}
\end{tikzcd}\]
to (\ref{eq: big diagram}). The rows and columns of this diagram are cofibre sequences. Since we are over $\F_2$,  the path components of the mapping space from $\bbS^{g,d,f}$ to itself
consist only of that of the zero map and that of the identity map, in particular  $\partial^{\mathrm{hor}}$ and $\partial^{\mathrm{vert}}$ are equivalent. Now the lemma follows from naturality.
\end{prf}

\begin{thm}
    Let $\widetilde{\gamma_{i\epsilon}}\defeq \widetilde{\varphi_{X_1}}^i(Q^1(\widetilde{\sigma})^\epsilon.\widetilde{\omega_{3,2}})$, $\epsilon= 0,1$,  and $\widetilde{\alpha_{i(\epsilon+1)}}\defeq \widetilde{\varphi_{X_1}}^i(Q^1(\sigma)^{\epsilon+1}\widetilde{\omega_{0,0}})$. \footnote{Here $X_1$ is the 2-cell complex $\E_\infty(\sigma)/\!/_{\E_\infty}\sigma Q^1(\sigma)$ built in \S\ref{sec: stability Hopf alg of CGL}. The name classes here are defined in \S\ref{sec: self map on fil X}. $X$ in the proof will refer to the filtered object $\Sigma^{1,0,-1}\fil^{\E_\infty}X_1$ as said.  We apologise for this symbol overload.}Consider the cofibre sequence 
    \[ \fil^{\E_\infty}X_1\xrightarrow{\widetilde{q}} \fil^{\E_\infty}X_1/\widetilde{\sigma}\xrightarrow[]{\widetilde{\partial}}\Sigma^{1,1,-1}\fil^{\E_\infty}X_1.\]
Then ${\widetilde{q}}_*{\widetilde{\partial}}_*(\widetilde{\gamma_{i\epsilon}})$ and $\tau\widetilde{\alpha_{i(\epsilon+1)}}$ have the same image in $\tensor*[^{\mathrm{BS}}]{E}{^\infty}(\fil^{\E_\infty}X/\widetilde{\sigma})$. 
\end{thm}
\begin{prf}
     Take $X\to Y\to Z$ to be 
     \[\Sigma^{1,0,-1}\fil^{\E_\infty}X\xrightarrow[]{\widetilde{\sigma}} \fil^{\E_\infty}X \xrightarrow[]{\widetilde{q}} \fil^{\E_\infty}X/\widetilde{\sigma}\]
     and
     \[\widetilde{z}= \widetilde{\omega_{12,8}}^i\widetilde{\omega_{3,2}}Q^1(\widetilde{\sigma})^{\epsilon},\quad y = \beta^{4i+1}Q^1(\sigma)^{\epsilon}\in \pi_{*,*,*}(\fil^{\E_\infty }X_1/\tau) = \tensor*[^{\mathrm{BS}}]{E}{^1_{*,*,*,0}}(\fil^{\E_\infty}X_1) \]
  in the above discussion. 
\begin{claim}[1]
    $(\mathrm{mod}\tau)_*\delta_1^\infty(y) = \sigma Q^1(\sigma)^{\epsilon+1}\beta^{4i}\in \pi_{*,*-1,*-1}(Y/\tau)$.
\end{claim}
\begin{prfoc}[1]
    Let  $\delta_1^2$ be the boundary map in the cofibre sequence
$Y/\tau\xrightarrow[]{\tau}Y/\tau^2 \xrightarrow[]{\mathrm{mod}\tau}Y/\tau$. 
Then there is a commutative triangle 
\[\begin{tikzcd}
    Y/\tau \arrow[r,"\delta_1^\infty"]\arrow[dr, "\delta_1^2"{below}]& \Sigma^{0,1,1}Y\arrow[d, "\mathrm{mod}\tau"]\\
     & \Sigma^{0,1,1}Y/\tau
\end{tikzcd}.\]
By the truncated omnibus \autocite[Theorem 2.21]{CDN24},  $\delta_1^2(y)$ is a representative of $d^2(y)$. There is only one such representative, namely the proposed class.
\end{prfoc}
\begin{claim}[2]
    Regardless of the choice involved in obtaining $\widetilde{x}$, we have
    \[(\mathrm{mod}\tau)_*(\widetilde{x}) = Q^1(\sigma)^{\epsilon+1}\beta^{4i}\in \pi_{*-1,*,*+1}(X/\tau) = \tensor*[^{\mathrm{BS}}]{E}{^1_{*-1,*,*+1}}(\fil^{\E_\infty}X_1)\]
\end{claim}
\begin{prfoc}We have
\[\sigma (\mathrm{mod}\tau)_*(\widetilde{x}) = (\mathrm{mod}\tau)_*(\widetilde{\sigma}\widetilde{x}) = (\mathrm{mod}\tau)_*(\delta_1^\infty(y))=\sigma Q^1(\sigma)^{\epsilon+1}\beta^{4i}\]
where the first equality is the naturality of $\mathrm{mod}\tau$, the second equality follows from the definition of $\widetilde{x}$, and the last equality is Claim (1).  Since the Bockstein $E_1$-page of $\fil^{\E_\infty}X_1$ is $\sigma$-torsion free, the claim follows.
\end{prfoc}
By Claim (2) $(\mathrm{mod}\tau)_*(q_*\widetilde{x})=(\widetilde{q}/\tau)_*(\mathrm{mod}\tau)_*=[Q_1(\sigma)^{\epsilon+1}\beta^{4i}]$, so \[q_*\widetilde{x}= \widetilde{\alpha_{i(\epsilon+1)}} +\tau \widetilde{\alpha}'\] for some $\widetilde{\alpha}'$. By the discussion at the start of the section $\widetilde{\partial}_*(\widetilde{\gamma_{i\epsilon}}) =\tau\widetilde{x}$, so 
\[\widetilde{q}_*\widetilde{\partial}_*(\widetilde{\gamma_{i\epsilon}}) =  \tau \widetilde{\alpha_{i(\epsilon+1)}} +\tau^2 \widetilde{\alpha}'.\]
This proves the result because the filtration on the filtered homotopy groups by the Bockstein spectral sequence is the $\tau$-adic filtration.
\end{prf}

By inverting $\tau$, we deduce

\begin{cor}\label{cor: collapsing lightning flash} 
    The classes $q^{\CGL}_*(u_{i\epsilon})$ and $\alpha_{i(\epsilon+1)}$  have the same image in $E^\infty_{*,*}(\fil^{\E_\infty}\CGL/\widetilde{\sigma})$. \hfill $\blacksquare$
\end{cor}

\section{Consequences for $\GL_n(\fancycat{O})$ and $\Aut(F_n)$}\label{sec: stability for Z}
Let $\fancycat{O}$ be a Dedekind domain that has $\Z/2$ as a quotient (e.g. $\fancycat{O} =\Z$). 
In this section prove some \say{$\frac{2}{3}$-stability} consequences for  $\GL_n(\fancycat{O})$ and $\Aut(F_n)$ (\cref{thm: transporting 2/3} and \cref{prop: GLnZ and AutFn}). First we need some setups for $\Aut(F_n)$.

Let $\cat{G}_{\Aut(F)}$ be the 1-groupoid with object non-negative integers, morphisms given by the automorphism group $\Aut(F_n)$ of free groups on $n$ letters. Let  $\cat{G}_{\Aut(F)}$ be equipped with the symmetric monoidal structure given by addition. The functor $r:\cat{G}_{\Aut(F)}\to \Z$ that given by inclusion on  objects is symmetric monoidal. Applying \cref{const: algebra from symmetric monoidal groupoid} we obtain an  $\E_\infty$-algebra  $\mathrm{CAut}(F)\defeq (\mathrm{B}\cat{G}_{\Aut(F)})_{\F_2}$ whose homotopy group  at grading $n$ computes the group homology $H_*(\Aut(F_n);\F_2)$. The homomorphisms $\Aut(F_n)\to \GL_n(\Z)$ induced by $\Z^{\oplus n} = F_n^{\mathrm{ab}}$ induces an $\E_\infty$-map $\mathrm{CAut}(F)\to \CGL(\Z)$, fitting into\stepcounter{thm}
\[\mathrm{CAut}(F)\to \CGL(\Z)\to\CGL(\F_2)\to a(1).\leqnomode\tag{\thethm}\label{eq: sequence from autFn}\]
Since $\Aut(F_1) = \GL_1(\Z)$, the stabilisation class $\sigma\in \pi_{1,0}(\CGL(\F_2))$ lifts all the way to $\mathrm{CAut}(F)$, such that the homotopy groups of $\mathrm{CAut}(F)/\sigma$ at grading $n$ computes $H_*(\Aut(F_n),\Aut(F_{n-1});\F_2)$. We would like to lift the endomorphism $\varphi$ to one on $\mathrm{CAut}(F)/\sigma$. However, since $\sigma Q^1(\sigma)$ is non-zero in $\mathrm{CAut}(F)$,  we cannot map $X_1$ to $\mathrm{CAut}(F)$. But 
\begin{lem}\label{lem: connectivity of AutF to CGLZ}
     $\pi_{n,d}(\mathrm{CAut}(F))\to \pi_{n,d}(\CGL(\Z))$ is an isomorphism when $n\leqslant 3$, $d\leqslant 1$. In particular, if we let $\theta_1\in \pi_{1,1}(\mathrm{CAut}(F))$ be 
     the class $\theta_1\in \pi_{1,1}(\CGL(\Z))$, then $\sigma^2\theta_1 = \sigma Q^1(\sigma)\in \pi_{3,1}(\mathrm{CAut}(F))$ and hence we have an $\E_\infty$-map
     \stepcounter{thm}
\[g_1:Y_1\defeq \E_\infty(\bbS^{1,0}\{\sigma\}\oplus \bbS^{1,1}\{\theta_1\})/\!/_{\E_\infty}(\sigma Q^1(\sigma)-\sigma^2\theta_1)\to \mathrm{CAut}(F).\leqnomode\tag{\thethm}\label{eq: g1}\]
\end{lem}
\begin{prf}
    The $d=0$ case is clear. For $d=1$, it's enough to check $\Aut(\F_n)\to \GL_n(\Z)$ is an isomorphism on abelianisation for $n =1,2,3$. For $n=1$, $\Aut(F_1) $ is just $\GL_1(\Z)$. For $n=2$, there is a well-known isomorphism $\mathrm{Out}(F_2)\cong \GL_2(\Z)$. Using the explicit presentation of of $\Aut(F_3)$ provided by \autocite[Corollary 1]{AFV07} one can compute $\Aut(F_3)^{\mathrm{ab}}\cong \Z/2$, which is also abstractly isomorphic to $\GL_3(\Z)^{\ab}$. We are done because $\Aut(F_3)^{\ab}\to \GL_3(\Z)^{\ab}$  is surjective.
\end{prf}
By \cref{thm: summary theorem of diag hopf} (iii) and  \cref{thm: bar and operations} (iii), the
$\E_\infty$-map $c:Y_1\to X_1$ given by collapsing the cell $\theta_1$ induces an isomorphism on the stability Hopf algebra.

\begin{thm}\label{thm: transporting 2/3}
    Let $A$ is be a connected $\E_\infty$-algebra in $\fancycat{D}(\F_2)^{\Z}$ satisfying $(SCE)$. Suppose there are $\E_\infty$ maps  
    \[Y_1\xrightarrow[]{i}A\xrightarrow[]{r}\CGL(\F_2)\]
    such that $\Delta_{ri}\circ (\Delta_c)^{-1}=\Delta_{f_1}:\Delta_{X_1}\to \Delta_{\CGL(\F_2)}$. Then $A/i_*\sigma$ does not have vanishing lines of slope $>\frac{2}{3}$.
\end{thm}
\begin{prf}Let $x_1,y_1, a$ denote respectively the $\E_\infty$ algebra $\Cobar(\Delta_{X_1})$, $\Cobar(\Delta_{Y_1})$, $\Cobar(\Delta_A)$. Since $y_1\xrightarrow[]{\cong}x_1$ is an equivalence, 
the endomorphism $\varphi_{x_1}\defeq \Cobar(\Delta_{\varphi_{X_1}})$ defines an endomorphism  on $y_1/\sigma$, which we call $\varphi_{y_1}$. By assumption, $\Cobar(\Delta_{ri})$ is equivalent to $\Cobar(\Delta_{f_1})$ under $x_1\cong  y_1$. In particular the composite with $\cgl\to a(1)$ (constructed in \S\ref{sec: cgl to a(1)}) is equivalent to $x_1\to a(1)$ and as a consequence, the basechange of  $\varphi_{y_1}$ along $y_1\to a(1)$ is equivalent to the basechange of $\varphi_{X_1}$ under $X_1\to a(1)$, which we saw is non nilpotent (\cref{prop: basechang as expected}). So $\varphi_{y_1}$ basechange to a non-nilpotent endomorphism on $a\otimes_AA/i_*(\sigma)$. Thus $a\otimes_AA/i_*(\sigma)$ does not have a slope $>\frac{2}{3}$ vanishing line. By \cref{basechange thm} (i), neither does $A/i_*(\sigma)$.
\end{prf}

\begin{prop}\label{prop: GLnZ and AutFn} Both
    \begin{enumerate}[(i)]
        \item $A =\CGL(\fancycat{O})$, where $\fancycat{O}$ is a Dedekind domain having $\Z/2$ as a quotient, and
        \item $A = \mathrm{CAut}(F)$ 
    \end{enumerate}
    satisfies the assumptions of \cref{thm: transporting 2/3}.
\end{prop}
\begin{prf} First we note that 
$\CGL(R)$ and $\mathrm{CAut}(F)$ are connected: If $r:\fancycat{G}\to \Z$ in \cref{const: algebra from symmetric monoidal groupoid} factors through the non-negative integers $\mathbb{N}$ and is such that the fibre at $0\in \Z$ is connected, then the resulting $\E_\infty$-algebra is connected. Asking for $R=\CGL(\Z), \mathrm{CAut}(F)$  to satisfy (SCE) is equivalent as asking the associated \textit{$\E_1$-splitting complexes} $S^{\E_1}(n)$ to be $(n-2)$-connected for every $n$, which are in turn the content of \autocite[Theorem 1.1]{Cha80} for $A = \CGL(R)$ and \autocite[Corollary 4.5]{Hep20} for $A =\mathrm{CAut}(F) $. We explain the details in  \S\ref{sec: splitting complex} below.
For the factorisation:
\begin{enumerate}[(i)]
    \item  The relation $\sigma^2\theta_1 = \sigma Q^1(\sigma)$ in $\CGL(\Z)$ defines an $\E_\infty$-map $Y_1\to \CGL(\Z)$. For $A = \CGL(\fancycat{O})$, we take $i$ to be the composite of $Y_1\to \CGL(\Z)$ with the map $\CGL(\Z)\to A$ and $r$ to be the map induced by $\fancycat{O}\to \Z$. 
    \item For $A = \mathrm{CAut}(F)$ we take $i$ to be the map (\ref{eq: g1}) and $r$ to be the one induced by $\Aut(F_n)\to \GL_n(\Z)\to \GL_n(\F_2)$. 
\end{enumerate}
The criterion on the Hopf algebras amounts to having $\overline{\sigma},\delta\in \Delta_{Y_1}$ maps to $\overline{\sigma},\delta\in \Delta_{\CGL(\F_2)}$. Since in $\Delta_{\CGL(\F_2)}$ the element $\overline{\sigma}$ is primitive and $\psi(\delta) =1\otimes\delta+\overline{\sigma}\otimes \overline{\sigma}^2 + \delta\otimes 1$, by considering coproducts and using the algebra structure on $\overline{\sigma} $
it is enough to show $\Delta_{Y_1}\to \Delta_{\CGL}(\F_2)$ is an isomorphism in gradings $1$. By \cref{cor: connectivity passing to bar}  we only need to show $Y_1\to \CGL(\F_2)$ is $(1,0)$-connected. This is indeed the case.
\end{prf}

\begin{prop}\label{prop: GLnZ AutFn stab}
    Let $A$ be either $\CGL(\Z)$ or $\mathrm{CAut}(F)$. Then $\pi_{n,d}(A/\sigma)=0$ for $d<\frac{2}{3}(n-1)$.
\end{prop}
\begin{prf} Since the map $Y_1\to A\to \CGL(\F_2)$ in \cref{thm: transporting 2/3} on stability Hopf algebras  is equivalent to the one induced by  $X_1\to \CGL(\F_2)$, which we know is an isomorphism through gradings $\leqslant 5$, $\Delta_{R}\to \Delta_{\CGL(\F_2)}$ is surjective through gradings $\leqslant 5$. On the other hand, $\Delta_{R}(n) = \Delta_{Y_1}(n)$ for $n\leqslant 2$ by \cref{cor: connectivity passing to bar} applied to \cref{prop: dyer-lashof basis for BGL(Z)}, \cref{lem: connectivity of AutF to CGLZ}. Thus  $\Delta_{R}\to \Delta_{\CGL(\F_2)}$  is injective through gradings $\leqslant 2$. Now by \cref{basechange thm} (ii) the map $r\to \cgl(\F_2)$ reflects vanishing lines of slope $\leqslant \min \{\frac{3-1}{3}, \frac{6-2}{6}\} = \frac{2}{3}$, and the conclusion follows from the known vanishing line \autocite[Theorem B]{GKRW18b} for $\CGL(\F_2)/\sigma$. 
\end{prf}

\subsubsection{The $\E_1$ splitting complexes are highly connected.}\label{sec: splitting complex} Let $\cat{G}$ be a symmetric monoidal $1$-groupoid with underlying object the natural numbers and the symmetric monoidal structure on objects by addition. Let $R$ be the $\E_\infty$-algebra in $\fancycat{D}(\F_2)^{\Z}$ obtained by applying \cref{const: algebra from symmetric monoidal groupoid} to $r:\cat{G}\to \Z$ that is $\mathbb{N}\hookrightarrow\Z$ on object. Associated to $(\cat{G},r)$ there is a $G_n\defeq \Aut_{\cat{G}}(n)$-semisimplicial set $\mathcal{S}^{\E_1}_{\bullet}$ for each $n\geqslant 2$, whose $p$-simplices are
\[\mathcal{S}^{\E_1}_{p}= \bigsqcup_{\substack{n_0+\cdots + n_{p+1}=n\\n_i>0}}\frac{G_n}{G_{n_0}\times \cdots \times G_{n_p}}\]
and face maps given by merging the adjacent factors. \autocite[Definition 17.9, Remark 17.11]{GKRW18a}. Let $S^{\E_1}(n)$ denote its thick realisation and  $S^{\E_1}(n)/\!/G_n$  the homotopy orbit, then there are natural isomorphisms
\[\pi_{n,*}(\Ba R)\cong H_{n,*-2}(S^{\E_1}(n)/\!/G_n;\F_2)\]
\autocite[Definition 17.9, Corollary 17.5, Lemma 17.10]{GKRW18a}. In particular, if each $S^{\E_1}(n)$ is $(n-3)$-connected, then $R$ satisfies (SCE). 

When $R = \CGL(\fancycat{O})$, $\fancycat{O}$ a Dedekind domain, $\St^{\E_1}_R(n)$ is homeomorphic to the nerve of the poset $\fancycat{T}^{\E_1}(n)  \defeq \{(P,Q)\}$ of splittings  $\fancycat{O}^{\oplus n} = P\oplus Q$ with $0\subsetneq P,Q\subsetneq \fancycat{O}^{\oplus n}$. \autocite[Theorem 1.1]{Cha80} says $\abs{\fancycat{T}^{\E_1}(n)}$ is homotopy equivalent to wedge of $(n-2)$-spheres.

For $R  =\mathrm{CAut}(F)$, we claim that the semisimplicial set $\mathcal{S}^{\E_1}_{\bullet}$ is isomorphic to the semisimplicial nerve\footnote{in the sense of \autocite[Definition 3.3]{ERW19}} of the non-unital poset $SP_n^{\circ}$, obtained from the poset $SP_n$ of \autocite[Proposition 3.5]{Hep20} by removing all identity relations. 
Recall objects of $SP_n$ are ordered pair of non-trivial proper subgroups $(P,Q)$ such that $P\ast Q = F_n$,  and $(P,Q)\leqslant (P',Q')$ if there are proper subgroups $J_0,J_1,J_2$ of $F_n$ such that $(P,Q) = (J_0, J_1\ast J_2)$, $(P', Q') = (J_0\ast J_1, J_2)$. In particular, 
a length-$p$ chain in $SP_n^{\circ}$ is same as an ordered tuple $(H_0,\dots, H_{p+1})$ of proper non-zero subgroups of $F_n$ with $H_0\ast\cdots \ast H_{p+1} =F_n $ (via the map $\lambda$ in the proof of \autocite[Corollary 4.5]{Hep20}). Fix a free basis $e_1,\cdots, e_n$ of $F_n$. The claimed isomorphism sends $g \in\frac{G_n}{G_{n_0}\times \cdots \times G_{n_p}}$ to  $g.(\<e_1,\dots, e_{n_0}>,\dots, \<e_{n-n_p+1},\dots, e_{n}>)$.  Note the unitalisation\footnote{in the sense of \autocite[Definition 3.7]{ERW19}.} $(SP_n^{\circ})^+$  of $SP_n^{\circ}$ is precisely $SP_n$, and by \autocite[]{ERW19} ${\norm{SP_n^{\circ}}}\to {\norm{SP_n^{\circ +}}}$ is a weak homotopy equivalence. The thick realisation $\norm{SP_n}$ is homotopy equivalent to the thin realisation $\abs{SP_n}$ \autocite[Lemma 1.7]{ERW19}.  That  $\abs{SP_n}$ is homotopy equivalent to the wedge of $(n-2)$-spheres is precisely the content of \autocite[Corollary 4.5]{Hep20}.

 % actual content
\printbibliography

\end{document}